\documentclass[authoryear]{elsarticle}

\usepackage{marginnote}

\usepackage{multirow}
\usepackage[ruled]{algorithm2e}
\usepackage{graphicx}  
\usepackage[font=small,format=plain,labelfont=bf,up,textfont=it,up]{caption}

\usepackage{todonotes}

\usepackage{subfigure}
\usepackage{geometry}
\usepackage{amsmath}
\usepackage{amssymb}
\usepackage{amsfonts}
\usepackage{amsthm}
\usepackage{url}
\usepackage{fancyhdr}
\usepackage{float}

\usepackage{pgf,tikz}
\usetikzlibrary{shapes,arrows,calc,automata}

\usepackage{ctable}

\let\oldref\ref

\renewcommand{\ref}[1]{(\oldref{#1})}

\title{\textbf{Multi-Level Spectral Deferred Corrections Scheme \\ for the Shallow Water Equations \\ on the Rotating Sphere}}
\date{}

\author[lbl]{Fran\c cois P. Hamon\corref{cor1}}
\ead{fhamon@lbl.gov}
\author[uexeter,umunich]{Martin Schreiber}
\author[lbl2]{Michael L. Minion}

\cortext[cor1]{Corresponding author}
\address[lbl]{Center for Computational Sciences and Engineering, Lawrence Berkeley National Laboratory, Berkeley, USA}
\address[uexeter]{Department of Mathematics/Computer Science, University of Exeter, Exeter, United Kingdom}
\address[umunich]{Chair of Computer Architecture and Parallel Systems, Technical University of Munich, Germany}
\address[lbl2]{Department of Applied Mathematics, Lawrence Berkeley National Laboratory, Berkeley, USA}

\begin{document}

\begin{abstract} 
Efficient time integration schemes are necessary to capture the complex processes 
involved in atmospheric flows over long periods of time. In this work, we propose a high-order,
implicit-explicit numerical scheme that combines Multi-Level Spectral Deferred Corrections (MLSDC)
and the Spherical Harmonics (SH) transform to solve the wave-propagation problems arising from 
the shallow-water equations on the rotating
sphere.

The iterative temporal integration is based on a sequence of corrections distributed on
coupled space-time levels to perform a significant portion of the calculations on a coarse 
representation of the problem and hence to reduce the time-to-solution while preserving accuracy. 
In our scheme, referred to as MLSDC-SH, the spatial discretization plays a key role in the efficiency 
of MLSDC, since the SH basis allows for consistent transfer functions between space-time 
levels that preserve important physical properties of the solution.

We study the performance of the MLSDC-SH scheme with shallow-water test cases commonly 
used in numerical atmospheric modeling. We use this suite of test cases, which gradually adds
more complexity to the nonlinear system of governing partial differential equations, to perform
a detailed analysis of the accuracy of MLSDC-SH upon refinement in time. We illustrate the stability 
properties of MLSDC-SH and show that the proposed scheme achieves up to eighth-order convergence in
time. Finally, we study the conditions in which MLSDC-SH achieves its theoretical speedup, and we
show that it can significantly reduce the computational cost compared to single-level Spectral Deferred Corrections (SDC).
\end{abstract}

\begin{keyword} high-order time integration, multi-level spectral deferred corrections, 
implicit-explicit splitting, atmospheric flows, shallow-water equations on the rotating 
sphere, spherical harmonics
\end{keyword}

\maketitle

\section{\label{section_introduction}Introduction}

The numerical modeling of global atmospheric processes presents a challenging application area
requiring accurate time integration methods for the discretized governing partial differential equations. 
These complex processes operate on a wide range of time scales but often have to be simulated over long 
periods of time -- up to a hundred years for long-term paleoclimate studies --  which constitutes a 
challenge for the design of stable and efficient integration schemes.  One strategy for creating more efficient
temporal integration schemes for such systems is to employ a semi-implicit scheme that
allows larger time steps to be taken than with explicit methods at a cost that is less than that of fully 
implicit methods \citep{giraldo2005semi}.
A second strategy is to use a parallel-in-time strategy to solve multiple time steps concurrently on multiple processors. 
Examples of parallel-in-time methods include Parareal \citep{lions2001resolution}, the Parallel Full Approximation 
Scheme in Space and Time, (PFASST, \cite{emmett2012toward}), and MultiGrid Reduction in Time 
(MGRIT, \cite{falgout2014parallel}). In this work, we consider semi-implicit, iterative, multi-level 
temporal integration methods based on Spectral Deferred Corrections (SDC) that are easily extended to high-order 
and also serve as a first step toward constructing parallel-in-time integration methods for the atmospheric 
dynamics based on PFASST. 

SDC methods are first presented in \cite{dutt2000spectral} and consist in applying a sequence of 
low-order corrections -- referred to as sweeps -- to a provisional solution in order to achieve high-order accuracy. 
Single-level SDC schemes have been applied to a wide range of problems, including reacting flow simulation 
\citep{bourlioux2003high,layton2004conservative}, atmospheric modeling \citep{jia2013spectral}, particle motion 
in magnetic fields \citep{winkel2015high}, and radiative transport modeling \citep{crockatt2017arbitrary}. 
In \cite{jia2013spectral}, a fully implicit SDC scheme is combined with the Spectral Element Method (SEM) to solve the 
shallow-water equations on the rotating sphere.
The authors demonstrate that the SDC method can take larger stable time steps than competing explicit schemes
such as leapfrog, second-order Runge-Kutta methods, and 
implicit second-order Backward Differentiation Formula (BDF) method without loss of accuracy. 

The approach considered here for atmospheric simulations builds on the work of \cite{speck2015multi}, in which a 
Multi-Level Spectral Deferred Corrections (MLSDC) scheme is proposed to improve the efficiency of the SDC time 
integration process while preserving its high-order accuracy. MLSDC relies on the construction of coarse 
space-time representations -- referred to as  levels -- of the problem under consideration. The calculations
are then performed on this hierarchy of levels in a way that shifts a significant portion of the computational 
burden to the coarse levels. As in nonlinear multigrid methods, the space-time levels are coupled by the 
introduction of a Full Approximation Scheme (FAS) term in the collocation problems solved on coarse levels. 
With this multi-level approach, the iterative correction process requires fewer fine sweeps than the standard 
SDC scheme but still achieves fast convergence to the fixed point solution.
Synthetic numerical examples demonstrate the efficiency and accuracy of the MLSDC approach. 

The MLSDC approach is combined here with a spatial discretization based on the 
global Spherical Harmonics (SH) transform to solve the shallow-water equations on the rotating sphere. This study is 
relevant for practical applications since the SH transform is implemented in major forecasting systems
such as the Integrated Forecast System (IFS) at the European Centre for Medium-Range Weather Forecasts (ECMWF, \cite{wedi2013fast}) 
and the Global Spectral Model (GSM) at the Japan Meteorological Agency (JMA, \cite{kanamitsu1983description}). 
Using a highly accurate method in space significantly reduces the spatial discretization errors and allows us to focus
on the temporal integration. Our approach, referred to as MLSDC-SH, uses a temporal splitting 
in which only the stiff linear terms in the governing equations are treated implicitly, whereas less stiff terms are 
evaluated explicitly. Here, the word \textit{stiff} is used to denote the terms that limit the time step size of fully 
explicit schemes. The temporal integration scheme retains the main features of the multi-level algorithm presented in 
\cite{speck2015multi}, and takes full advantage of the structure of the spatial discretization to achieve efficiency. 
Specifically, we construct accurate interpolation and restriction functions between space-time levels by padding or 
truncating the spectral representation of the variables in the SH transform. In addition, the spherical harmonics 
combined with the implicit-explicit temporal splitting considered in this work circumvent the need for a global linear 
solver and rely on an efficient local solver for the implicit systems. 

We illustrate the properties of MLSDC-SH using a widely used suite of shallow-water test cases 
\citep{williamson1992standard,galewsky2004initial}. We start the numerical study
with a steady-state benchmark that highlights the connection between the magnitude of the spectral coefficients truncated
during coarsening and the convergence rate of MLSDC-SH upon refinement in time. Then we proceed to more challenging unsteady test cases to
show that MLSDC-SH is stable for large time steps and achieves up to eighth-order temporal convergence. Finally, we investigate
the conditions in which the proposed scheme achieves its theoretical speedup and we demonstrate that MLSDC-SH can reduce the 
computational cost compared to single-level SDC schemes.

In the remainder of the paper, we first introduce the system of governing equations in Section~\oldref{section_governing_equations}. 
Then, we briefly review the fundamentals of the spatial discretization based on the global SH 
transform in Section~\oldref{section_spatial_discretization}. In Section~\oldref{section_temporal_discretization},
we describe the implicit-explicit temporal integration scheme, with an emphasis on the Multi-Level
Spectral Deferred Corrections (MLSDC) scheme. Finally, in Section~\oldref{section_numerical_examples}, 
we present numerical examples on the sphere demonstrating the efficiency and accuracy of our approach.

\section{\label{section_governing_equations}Governing equations}

We consider the Shallow-Water Equations (SWE) on the rotating sphere. These equations capture the main 
horizontal effects present in the full atmospheric equations.
Well-defined test cases are available -- such as those considered in this work -- that relate the 
SWE to some key features of the full atmospheric equations. Hence, they provide a simplified assessment 
of the properties of temporal and spatial discretizations for atmospheric simulations on the rotating 
sphere. We use the vorticity-divergence formulation \citep{bourke1972efficient,hack1992description} in 
which the prognostic variables $\boldsymbol{U} = [ \Phi, \, \zeta, \, \delta ]^T$ are respectively the 
potential, $\Phi$, the vorticity, $\zeta$, and the divergence, $\delta$. Here, the vorticity and 
divergence state variables are used to overcome the singularities in the velocity field at the poles.

The system of governing partial differential equations is
\begin{align}
\frac{\partial \Phi'}{\partial t}  &= - \nabla \cdot ( \Phi' \boldsymbol{V} ) - \bar{\Phi} \delta + \nu \nabla^2 \Phi', \label{geopotential_equation} \\ 
\frac{\partial \zeta}{\partial t}  &= - \nabla \cdot ( \zeta + f ) \boldsymbol{V} + \nu \nabla^2 \zeta, \label{vorticity_equation} \\
\frac{\partial \delta}{\partial t} &= \boldsymbol{k} \cdot \nabla \times (\zeta + f) \boldsymbol{V} - \nabla^2 \bigg( \Phi + \frac{\boldsymbol{V} \cdot \boldsymbol{V}}{2} \bigg) + \nu \nabla^2 \delta, \label{divergence_equation}
\end{align}
where $\boldsymbol{k}$ is the outward radial unit vector. The average geopotential, $\bar{\Phi} = g \bar{h}$, is written as 
the product of the gravitational acceleration by the average height, and $\Phi'$ is defined as $\Phi' = \Phi - \bar{\Phi}$. 
The horizontal velocity vector is $\boldsymbol{V} \equiv \boldsymbol{i}u + \boldsymbol{j}v$, where $\boldsymbol{i}$ and $\boldsymbol{j}$
are the unit vectors in the eastward and northward directions, respectively. 
The Coriolis force is represented by $f = 2 \Omega \sin \phi$, where $\Omega$ is the angular rate of rotation, 
and $\phi$ is the latitude. The diffusion coefficient is denoted by $\nu$. Including a diffusion term in the governing equations 
is used in practice in atmospheric simulations to stabilize the flow dynamics and reduce the errors caused by nonlinearly interacting 
modes. Using the inviscid equations is not a viable option due to the extremely fast generation of small-scale features
\citep{galewsky2004initial}, in particular for global spectral methods using a collocated grid. 
For simplicity and reproducibility, we employ a second-order diffusion term with a diffusion coefficient set to 
$\nu = 1.0 \times 10^5 \, \text{m}^2.\text{s}^{-1}$ for all spatial resolutions as in \cite{galewsky2004initial}.
To express the velocities as a function of the prognostic variables, $\zeta$ and $\delta$, we first use the 
Helmholtz theorem which relates $\boldsymbol{V}$ to a scalar stream function, $\psi$, and a scalar velocity potential, $\chi$,
\begin{equation}
\boldsymbol{V} = \boldsymbol{k} \times \nabla \psi + \nabla \chi.
\label{helmholtz_theorem}
\end{equation}
Using the identities
\begin{align}
\zeta  &\equiv \boldsymbol{k} \cdot ( \nabla \times \boldsymbol{V} ), \label{identity_defining_zeta} \\
\delta &\equiv \nabla \cdot \boldsymbol{V}, \label{identity_defining_delta} 
\end{align}
the application of the curl and divergence operators to \ref{helmholtz_theorem} yields
$\zeta  =  \nabla^2 \psi$ and $\delta =  \nabla^2 \chi$.
The Laplacian operators can be efficiently inverted using the SH transform
to compute the stream function, $\psi$, and the velocity potential, $\chi$, as a function of $\zeta$ 
and $\delta$, as explained in Section~\oldref{section_spatial_discretization}. Equations 
\ref{geopotential_equation}, \ref{vorticity_equation}, and \ref{divergence_equation},
form the system that we would like to solve.

Next, we use the identities \ref{identity_defining_zeta} and \ref{identity_defining_delta} to split 
the right-hand side of \ref{geopotential_equation}, \ref{vorticity_equation}, 
and \ref{divergence_equation} into linear and nonlinear parts as follows
\begin{equation}
\frac{\partial \boldsymbol{U}}{\partial t} = \boldsymbol{\mathcal{L}}_G( \boldsymbol{U} )  + \boldsymbol{\mathcal{L}}_F ( \boldsymbol{U} ) 
                                           + \boldsymbol{\mathcal{N}} ( \boldsymbol{U} ).
\label{linear_nonlinear_decomposition}
\end{equation}
The first term in the right-hand side of \ref{linear_nonlinear_decomposition} represents the linear wave motion 
induced by gravitational forces and also includes the diffusion term
\begin{equation}
\boldsymbol{\mathcal{L}}_G ( \boldsymbol{U} ) \equiv [ - \bar{\Phi} \delta + \nu \nabla^2 \Phi', \, \, \nu \nabla^2 \zeta, \, \, - \nabla^2 \Phi + \nu \nabla^2 \delta]^T. \label{linear_wave_motion_induced_by_gravitational_forces} 
\end{equation}
The second term in the right-hand side of \ref{linear_nonlinear_decomposition} contains a linear harmonic 
oscillator on the velocity components that includes the Coriolis term 
\begin{align}
\boldsymbol{\mathcal{L}}_F ( \boldsymbol{U} ) \equiv [ 0, \, \, -f \delta - \boldsymbol{V} \cdot \nabla f, \, \, f \zeta + \boldsymbol{k} \cdot (\nabla f) \times \boldsymbol{V} ]^T. \label{linear_harmonic_oscillator_on_the_velocity_components}
\end{align}
The third term in the right-hand side of \ref{linear_nonlinear_decomposition} represents 
the nonlinear operators
\begin{equation}
\boldsymbol{\mathcal{N}} ( \boldsymbol{U} ) \equiv \left[ - \nabla \cdot ( \Phi' \boldsymbol{V} ), \, \, - \nabla \cdot ( \zeta \boldsymbol{V}), \, \,  \boldsymbol{k} \cdot \nabla \times ( \zeta \boldsymbol{V} ) - \nabla^2 \frac{ \boldsymbol{V} \cdot \boldsymbol{V} }{2} \right]^T. \label{nonlinear_operators}
\end{equation} 
In Section~\oldref{section_temporal_discretization}, this decomposition is used to define the  
temporal implicit-explicit splitting chosen based on the stiffness of the different terms. 
Next, the details of the spatial discretization of $\boldsymbol{\mathcal{L}}_F$, $\boldsymbol{\mathcal{L}}_G$, and 
$\boldsymbol{\mathcal{N}}$  are presented.

\section{\label{section_spatial_discretization}Spatial discretization}

This section presents an overview of the spatial discretization based on the global SH transform
applied to the system of governing equations. The global SH transform is a key feature 
of the multi-level scheme presented here since it allows for simple and accurate data transfer between 
different spatial levels. We will show with numerical examples in Section~\oldref{section_numerical_examples} 
that this is critical for the design of efficient MLSDC schemes. In the SH scheme, the 
representation of a function of longitude $\lambda$ and Gaussian latitude 
$\mu \equiv \sin( \phi )$, $\xi(\lambda,\mu)$, consists of a sum of spherical harmonic 
basis functions $P^r_s(\mu) e^{i r \lambda}$ weighted by the spectral coefficients 
$\xi^r_s$,
\begin{equation}
\xi(\lambda, \mu) = \sum^R_{r = -R} \sum^{S(r)}_{s = |r|} \xi^r_s P^r_s(\mu) e^{i r \lambda},
\label{spectral_to_physical_space}
\end{equation}
where the index $r$ (respectively, $s$) refers to the latitudinal (respectively, longitudinal)
mode. In \ref{spectral_to_physical_space}, $P^r_s$ is the normalized associated Legendre 
polynomial. Without loss of generality, we use a triangular truncation with $S(r) = R$.
In Section~\oldref{subsubsection_coarsening_strategy_and_transfer_functions}, we will explain that a coarse
representation of $\xi$ can be obtained by simply truncating the number of modes -- i.e., reducing $R$ 
and $S$ in \ref{spectral_to_physical_space} -- to construct a hierarchy of spatial levels with different degrees of coarsening
in MLSDC-SH. The transformation from physical to spectral space is achieved in two steps. The first step
consists in taking the discrete Fourier transform of $\xi(\lambda, \mu)$ in longitude -- i.e., over $\lambda$ --, 
defined as
\begin{equation}
\xi^r(\mu) = \frac{1}{I} \sum_{\iota=1}^I \xi(\lambda_\iota,\mu) e^{-i r \lambda_\iota}, 
\label{fourier_transform}
\end{equation}
where $I$ denotes the number of grid points in the longitudinal direction, located at longitudes
$\lambda_{\iota} = \frac{2 \pi \iota}{I}$. Then, in the second step, the application of the discrete Legendre transformation 
in latitude yields 
\begin{equation}
\xi^r_s = \sum_{j = 1}^J \xi^r(\mu_j) P^r_s(\mu_j) w_j.
\label{legendre_transform}
\end{equation}
In \ref{legendre_transform}, $J$ is the number of Gaussian latitudes $\mu_j$, chosen as the roots
of the Legendre polynomial of degree $J$, $P_J$, and $w_j$ denotes the Gaussian weight at latitude $\mu_j$.
This two-step global transform is applied to \ref{linear_nonlinear_decomposition} to obtain a system of coupled 
ordinary differential equations involving the prognostic variables in spectral space, 
$\boldsymbol{\Theta}^r_s = [ \Phi^r_s, \, \zeta^r_s, \, \delta^r_s ]$. Note that due to the symmetry of 
the spectral coefficients, it is sufficient to only include the indices $r \geq 0$. Hence,  for 
$r \in \{ 0, \dots, R \}$ and $s \in \{ r, \dots, R \}$, the equations are
\begin{equation}
\frac{\partial \boldsymbol{\Theta}^r_s}{\partial t} = (\boldsymbol{L}_{G})^r_s( \boldsymbol{\Theta} )
                                                   + (\boldsymbol{L}_{F})^r_s ( \boldsymbol{\Theta} )
                                                   + \boldsymbol{N}^r_s ( \boldsymbol{\Theta} ),
\label{semidiscrete_system}
\end{equation}
where $(\boldsymbol{L}_G)^r_s$, $(\boldsymbol{L}_F)^r_s$, and $\boldsymbol{N}^r_s$ are the discrete, spectral
representations of the operators defined in \ref{linear_wave_motion_induced_by_gravitational_forces}, 
\ref{linear_harmonic_oscillator_on_the_velocity_components}, and \ref{nonlinear_operators}. 
The state variable in spectral space, $\boldsymbol{\Theta}$, is defined as a vector 
of size $K = 3 R (R+1)/2$ as follows
\begin{equation}
\boldsymbol{\Theta} \equiv [ \boldsymbol{\Theta}^0_0, \, \boldsymbol{\Theta}^{1}_{0}, \dots, \, \boldsymbol{\Theta}^{R}_{R-1}, \, \boldsymbol{\Theta}^R_R ]^T.
\label{spectral_space}
\end{equation}
We 
refer to the work of \cite{hack1992description} for a thorough presentation of the 
scheme, including the full expression of the right-hand side of \ref{semidiscrete_system} in spectral
space. More details about an efficient implementation of the global SH transform
can be found in \cite{temperton1991scalar,rivier2002efficient}.
The implementation of the spherical harmonics transformation used in this work is based on the SHTns library developed by \cite{schaeffer2013efficient}.
Next, we proceed to the presentation of the
discretization in time based on spectral deferred corrections.

\section{\label{section_temporal_discretization}Temporal discretization}

\subsection{\label{subsection_temporal_splitting}Temporal splitting}

The choice of a temporal splitting for the right-hand side of \ref{semidiscrete_system} is one of
the key determinants of the performance of the scheme. Fully explicit schemes are based on 
inexpensive local updates but are limited by a severe stability restriction on the time step 
size. In the context of atmospheric modeling, this limitation is often caused by the presence
of fast waves (e.g., sound or gravity waves) propagating in the system.
Fully implicit 
schemes overcome the stability constraint on the time step size but rely on costly nonlinear 
global implicit solves to update all the degrees of freedom simultaneously 
\citep{evans2010accuracy,jia2013spectral,lott2015algorithmically}. 

Instead, implicit-explicit (IMEX) schemes only treat the stiff terms responsible for the 
propagation of the fast-moving waves implicitly, while the non-stiff terms that represent 
processes operating on a slower time scale are evaluated explicitly. This strategy reduces 
the cost of the implicit solves, relative to fully implicit solves, and allows for relatively large stable time steps.
A common IMEX approach employed in non-hydrostatic atmospheric modeling is based on dimensional 
splitting and implicitly discretizes only the terms involved in the (fast) vertical dynamics 
\citep{ullrich2012operator,durran2012implicit,weller2013runge,giraldo2013implicit,lock2014numerical,gardner2018implicit}.
Alternatively, the approach of \cite{robert1972implicit,giraldo2005semi} consists in linearizing
  the governing PDEs in the neighborhood of a reference state. The linearized piece is then discretized implicitly, and the term
  treated explicitly is obtained by subtracting the linearized piece from the nonlinear system.
  For the shallow-water equations, we directly discretize the fast linear
  terms on the right-hand side of \ref{semidiscrete_system} implicitly, while the other terms are evaluated explicitly.
Specifically, we investigate an IMEX scheme based on the following splitting, for $r \in \{ 0, \dots, R \}$ 
and $s \in \{ r, \dots, R \}$,
\begin{equation}
\frac{\partial \boldsymbol{\Theta}^r_s}{\partial t} = (\boldsymbol{F}_I)^r_s(\boldsymbol{\Theta}) 
                                              + (\boldsymbol{F}_E)^r_s(\boldsymbol{\Theta}),
\label{system_of_odes}
\end{equation}
in which the implicit right-hand side, $(\boldsymbol{F}_I)^r_s$, contains the terms 
representing linear wave motion induced by gravitational forces and the diffusion term. The explicit right-hand 
side, $(\boldsymbol{F}_E)^r_s$, contains the linear harmonic oscillator and the nonlinear terms. 
This temporal splitting leads to implicit and explicit right-hand sides defined as
\begin{align}
(\boldsymbol{F}_I)^r_s &\equiv (\boldsymbol{L}_G)^r_s, \label{system_of_odes_implicit_part} \\
(\boldsymbol{F}_E)^r_s &\equiv (\boldsymbol{L}_F)^r_s + \boldsymbol{N}^r_s. \label{system_of_odes_explicit_part}
\end{align}
This implicit-explicit approach greatly simplifies the solution strategy for the implicit systems
and circumvents the need for a global linear solver. The solution algorithm treats the geopotential separately 
from the divergence and vorticity variables. Since the Coriolis term and the nonlinear terms are treated 
explicitly, one can form a diagonal linear system in spectral space to update the geopotential, and then update locally the vorticity 
and divergence variables. This is explained in Section~\oldref{implicit_solver_for_sdc_and_mlsdc}.
We will investigate the stability and accuracy of the splitting with numerical examples in Section 
\oldref{section_numerical_examples}. Next, we describe the multi-level temporal integration scheme 
starting with the fundamentals of SDC.

\subsection{\label{subsection_implicit_explicit_spectral_deferred_correction}IMEX Spectral Deferred Corrections }

We start with a review of the fundamentals of the Spectral Deferred Corrections (SDC) scheme. SDC methods have been introduced in 
\cite{dutt2000spectral} and later extended to methods with different temporal splittings in
\cite{minion2003semi,bourlioux2003high,layton2004conservative}.  In \cite{minion2003semi}, an implicit-explicit SDC method is described and
referred to as {\it semi-implicit SDC} to contrast the method with subsequent {\it multi-implicit SDC} methods with multiple implicit terms
introduced in \cite{bourlioux2003high}.  Here we employ the more used term {\it IMEX} to refer to SDC methods with an implicit-explicit splitting.
The properties of IMEX SDC schemes for fast-wave slow-wave problems are analyzed in \cite{ruprecht2016spectral}.
We consider a system of coupled ODEs in the generic form
\begin{align}
\frac{\partial \boldsymbol{\Theta} }{\partial t} (t) &= \boldsymbol{F}_I \big( \boldsymbol{\Theta}(t) \big) + \boldsymbol{F}_E \big( \boldsymbol{\Theta}(t) \big), \qquad t \in [t^n,t^{n} + \Delta t], \\
\boldsymbol{\Theta}(t^n) &= \boldsymbol{\Theta}^n,
\end{align}
and its solution in integral form given by
\begin{equation}
  \boldsymbol{\Theta}(t) = \boldsymbol{\Theta}^n + \int^t_{t^n} (\boldsymbol{F}_I + \boldsymbol{F}_E) \big( \boldsymbol{\Theta}(a) \big)d a 
                         = \boldsymbol{\Theta}^n + \int^t_{t^n} \boldsymbol{F} \big( \boldsymbol{\Theta}(a) \big)d a,
  \label{integral_form}
\end{equation}
where $\boldsymbol{F}_I$ and $\boldsymbol{F}_E$ are the implicit and explicit right-hand sides, respectively, with 
$\boldsymbol{F} = \boldsymbol{F}_I + \boldsymbol{F}_E$, and $\boldsymbol{\Theta}(t)$ is the state variable at time $t$. 
In \ref{integral_form}, the integral is applied componentwise.  Denote by $\tilde{\boldsymbol{\Theta}}(t)$ an approximation 
of $\boldsymbol{\Theta}(t)$, and then define the correction $\boldsymbol{\Delta} \boldsymbol{\Theta}(t) = \boldsymbol{\Theta}(t) - \tilde{\boldsymbol{\Theta}}(t)$. 
The SDC scheme applied to the implicit-explicit temporal splitting described above iteratively improves the accuracy of the 
approximation based on a discretization of the update or correction equation
\begin{align}
\tilde{\boldsymbol{\Theta}}(t) + \boldsymbol{\Delta}\boldsymbol{\Theta}(t)   = \boldsymbol{\Theta}^{n} 
& + \int_{t^{n}}^{t} \big[ \boldsymbol{F}_E \big( \tilde{\boldsymbol{\Theta}}(a) + \boldsymbol{\Delta}\boldsymbol{\Theta}(a) \big) - \boldsymbol{F}_E \big( \tilde{\boldsymbol{\Theta}}(a) \big) \big] d a \nonumber \\
& + \int_{t^{n}}^{t} \big[ \boldsymbol{F}_I \big( \tilde{\boldsymbol{\Theta}}(a) + \boldsymbol{\Delta}\boldsymbol{\Theta}(a) \big) - \boldsymbol{F}_I \big( \tilde{\boldsymbol{\Theta}}(a) \big) \big] d a \nonumber \\
& + \int_{t^n}^{t} \boldsymbol{F} \big( \tilde{\boldsymbol{\Theta}}(a) \big)  d a,
    \label{standard_sdc_sweep_implicit}
\end{align}
where $\boldsymbol{\Theta}^{n}$ is the (known) state variable at the beginning of the time step. In the update equation
\ref{standard_sdc_sweep_implicit}, the last integral is computed with a high-order Gaussian quadrature rule.
However, the other integrals are approximated with simpler low-order quadrature rules.
We mention here that all the quadrature rules used in \ref{standard_sdc_sweep_implicit} are
based on a relatively small number of Gauss points -- up to five in this work -- compared to the quadrature
rule used in the discrete Legendre transform to obtain \ref{legendre_transform}.
Each pass of the discrete version of the update equation \ref{standard_sdc_sweep_implicit}, referred to as sweep, increases 
the formal order of accuracy by one until the order of accuracy of the quadrature applied to the third 
integral is reached \citep{hagstrom2007spectral,xia2007efficient,christlieb2009comments}. 

To discretize the update equation \ref{standard_sdc_sweep_implicit}, the correction algorithm 
uses a decomposition of the time interval $[t^n, t^{n+1}]$ into $M$ subintervals using $M+1$ temporal nodes, such that
\begin{equation}
  t^n \equiv t^{n,0} < t^{n,1} < \dots < t^{n,M} = t^n + \Delta t \equiv t^{n+1}.
\label{time_discretization_sdc_nodes}
\end{equation}
The points $t^{n,m}$ are chosen to correspond to Gaussian quadrature nodes. Throughout this paper, we use Gauss-Lobatto nodes. 
We use the shorthand notations $t^m = t^{n,m}$ and $\Delta t^m = t^{m+1} - t^m$. We denote by $\boldsymbol{\Theta}^{m+1,(k+1)}$ 
the approximate solution at node $m+1$ and at sweep $(k+1)$. The terms in the first integral of \ref{standard_sdc_sweep_implicit}
are treated explicitly, and therefore this integral is discretized with a forward Euler method. Conversely, the second integral 
in \ref{standard_sdc_sweep_implicit} is discretized implicitly.
The general form of the discrete version of equation \ref{standard_sdc_sweep_implicit} is then
\begin{align}
\boldsymbol{\Theta}^{m+1,(k+1)} = \boldsymbol{\Theta}^{n} 
                              &+ \Delta t \sum_{j = 1}^m \tilde{q}^E_{m+1,j} \big[ \boldsymbol{F}_E \big( \boldsymbol{\Theta}^{j,(k+1)} \big) 
                                                                                                - \boldsymbol{F}_E \big( \boldsymbol{\Theta}^{j,(k)}   \big) \big]  \nonumber \\
                              &+ \Delta t \sum_{j = 1}^{m+1} \tilde{q}^I_{m+1,j} \big[ \boldsymbol{F}_I \big( \boldsymbol{\Theta}^{j,(k+1)} \big) 
                                                                                                - \boldsymbol{F}_I \big( \boldsymbol{\Theta}^{j,(k)}   \big) \big] \nonumber \\
                              &+ \Delta t \sum_{j = 0}^{M} q_{m+1,j} \boldsymbol{F} \big( \boldsymbol{\Theta}^{j,(k)} \big). \label{update_equation_sdcq_discrete_form} 
\end{align}
In \ref{update_equation_sdcq_discrete_form},
the coefficients $\tilde{q}^E_{m+1,j}$ correspond to forward-Euler time stepping.
The coefficients $q_{m+1,j}$ correspond to the Lobatto IIIA optimal order collocation quadrature
\begin{equation}
q_{m+1,j} \equiv \frac{1}{\Delta t} \int_{t^{n,0}}^{t^{{n,m+1}}} L^j(a) da,
\label{weights}
\end{equation}
where $L^j$ denotes the $j^{\text{th}}$ Lagrange polynomial constructed using the SDC nodes \ref{time_discretization_sdc_nodes}.
We note that the formulation of the correction given in \ref{update_equation_sdcq_discrete_form} differs from that of \cite{jia2013spectral}
in two ways. First, our scheme is based on an implicit-explicit splitting, whereas that of \cite{jia2013spectral} is fully implicit.
Second, for the choice of the quadrature weights used in the discretization of the implicit correction integral, we adopt 
the approach of \cite{weiser2015faster}. Specifically, the weights $\tilde{q}^I_{m+1,j}$ in \ref{update_equation_sdcq_discrete_form}
are chosen to be the coefficients of the upper triangular matrix in the LU decomposition of $\boldsymbol{Q} = \{ q_{ij} \} \in \mathbb{R}^{(M+1)\times(M+1)}$,
while a diagonal matrix is used in \cite{jia2013spectral}.
This formulation leads to a faster convergence of the iterative process to the fixed-point solution and remains convergent even when
the underlying problem is stiff. We refer to \cite{weiser2015faster} for a proof, and to \cite{hamon2018concurrent} for numerical examples
illustrating the improved convergence.

Using these definitions, the integration scheme \ref{update_equation_sdcq_discrete_form} is effectively an iterative solution method for the 
collocation problem defined by
\begin{equation}
\boldsymbol{A} ( \vec{\boldsymbol{\Theta}} ) = \boldsymbol{1}_{M+1} \otimes \boldsymbol{\Theta}^{n,0}. \label{collocation_problem}
\end{equation}
The operator $\boldsymbol{A}$ is 
\begin{equation}
\boldsymbol{A} ( \vec{\boldsymbol{\Theta}} ) \equiv \vec{\boldsymbol{\Theta}} 
                                           - \Delta t ( \boldsymbol{Q} \otimes \boldsymbol{I}_{K} ) \vec{\boldsymbol{F}},
\end{equation}
where $\otimes$ denotes the Kronecker product and  
$\boldsymbol{I}_{K} \in \mathbb{R}^{K \times K}$ is the identity matrix. $\boldsymbol{1}_{M+1} \in \mathbb{R}^{M+1}$ is 
a vector of ones. Following the notation used in \cite{bolten2016multigrid}, the space-time vectors 
$\vec{\boldsymbol{\Theta}} \in \mathbb{C}^{(M+1)K}$ and $\vec{\boldsymbol{F}} \in \mathbb{C}^{(M+1)K}$ are such that
\begin{align}
\vec{\boldsymbol{\Theta}} &\equiv [\boldsymbol{\Theta}^{n,0}, \dots, \boldsymbol{\Theta}^{n,M}]^T,  \\
\vec{\boldsymbol{F}} &\equiv \vec{\boldsymbol{F}} ( \vec{\boldsymbol{\Theta}} ) = [\boldsymbol{F}( \boldsymbol{\Theta}^{n,0} ), \dots, \boldsymbol{F}( \boldsymbol{\Theta}^{n,M}) ]^T.
\end{align}
Next, we introduce the multi-level algorithm based on SDC that we will apply to the collocation problem \ref{collocation_problem}.

\subsection{\label{subsection_multi_level_spectral_deferred_correction}Multi-Level Spectral Deferred Corrections (MLSDC)}

Multi-Level Spectral Deferred Corrections (MLSDC) schemes are based on the idea of replacing some of the SDC iterations required
to converge to the collocation problem  \ref{collocation_problem} with SDC sweeps performed on a coarsened (and hence computationally cheaper) 
version of the problem. The solutions on different levels are coupled by the introduction of a Full Approximation Scheme (FAS) correction 
term explained below as in nonlinear multigrid methods. The combination of performing SDC sweeps on multiple space-time levels with a 
FAS correction term first appears as part of the PFASST method in \cite{emmett2012toward}.  The idea is generalized and analyzed in 
\cite{speck2015multi} showing how MLSDC can improve the efficiency for certain problems compared to single-level SDC methods.
In this work, only two-level MLSDC schemes are considered, and the study of MLSDC with three or more space-time levels to 
integrate the shallow-water equations is left for future work. 

\subsubsection{Full Approximation Scheme (FAS)}

We define two space-time levels to solve the collocation problem \ref{collocation_problem}, and we denote by $\ell = f$ (respectively, 
$\ell = c$) the fine level (respectively, the coarse level). We denote by $\vec{\boldsymbol{\Theta}}_{\ell} \in \mathbb{C}^{(M_{\ell}+1) K_{\ell}}$ 
and $\vec{\boldsymbol{F}}_{\ell} \in \mathbb{C}^{(M_{\ell}+1) K_{\ell}}$ the space-time vector and right-hand side at level $\ell$, respectively. 
The matrix $\boldsymbol{R}^{c}_{f} \in \mathbb{R}^{ (M_{c}+1)K_{c} \times  (M_{f} + 1)K_{f}}$ is the linear restriction operator
from the fine level to the coarse level. Here, $K_{\ell}$ represents the total number of spectral coefficients in \ref{spectral_space} on level $\ell$.
As in nonlinear multigrid methods \citep{brandt1977multi}, the coarse problem is modified by 
the introduction of a correction term, denoted by $\vec{\boldsymbol{\tau}}_{c}$, that couples the solutions at the two space-time levels. 
Specifically, the coarse problem reads
\begin{equation}
\boldsymbol{A}_{c} ( \vec{\boldsymbol{\Theta}}_{c} ) - \vec{\boldsymbol{\tau}}_{c} = \boldsymbol{1}_{M_{c}+1} \otimes \boldsymbol{\Theta}^{n,0}_{c},
\label{collocation_problem_level_l+1}
\end{equation}
where the FAS correction term at the coarse level is defined as
\begin{equation}
\vec{\boldsymbol{\tau}}_{c} \equiv \boldsymbol{A}_{c} ( \boldsymbol{R}^{c}_{f} \vec{\boldsymbol{\Theta}}_{f} ) 
                           - \boldsymbol{R}^{c}_{f} \boldsymbol{A}_{f} ( \vec{\boldsymbol{\Theta}}_{f} )  
                           + \boldsymbol{R}^{c}_{f} \vec{\boldsymbol{\tau}}_{f},
\label{tau_term_level_l+1}
\end{equation}
with, for the two-level case, $\vec{\boldsymbol{\tau}}_f = \boldsymbol{0}$ on the fine level. 
In \ref{collocation_problem_level_l+1}-\ref{tau_term_level_l+1}, the operator
$\boldsymbol{A}_{c}$ denotes an approximation of $\boldsymbol{A}$ at the coarse level. We note that 
\begin{align}
\boldsymbol{A}_{c} ( \boldsymbol{R}^{c}_{f} \vec{\boldsymbol{\Theta}}_{f} ) - \vec{\boldsymbol{\tau}}_{c} &= 
   \boldsymbol{A}_{c} ( \boldsymbol{R}^{c}_{f} \vec{\boldsymbol{\Theta}}_{f} ) 
 - \boldsymbol{A}_{c} ( \boldsymbol{R}^{c}_{f} \vec{\boldsymbol{\Theta}}_{f} ) 
 + \boldsymbol{R}^{c}_{f}  \boldsymbol{A}_{f} ( \vec{\boldsymbol{\Theta}}_{f} )
 - \boldsymbol{R}^{c}_{f}  \vec{\boldsymbol{\tau}}_{f}  \nonumber \\
&= \boldsymbol{R}^{c}_{f} \big( \boldsymbol{A}_{f} ( \vec{\boldsymbol{\Theta}}_{f} )- \vec{\boldsymbol{\tau}}_{f} \big) ,
\end{align}
which implies that the restriction of the fine solution, $\boldsymbol{R}^{c}_{f} \vec{\boldsymbol{\Theta}}_{f}$, is a solution of the coarse problem. 
On the coarse problem \ref{collocation_problem_level_l+1}, the modified SDC update for temporal node $m+1$ at sweep $(k+1)$ is
\begin{align}
\boldsymbol{\Theta}^{m+1,(k+1)}_{c} &= \boldsymbol{\Theta}^{n,0}_{c} \nonumber \\
                                &+ \Delta t \sum_{j = 1}^m (\tilde{q}^E_{m+1,j})_{c} \big[ \boldsymbol{F}_{E,c} \big( \boldsymbol{\Theta}^{j,(k+1)}_{c} \big) 
                                                                                                - \boldsymbol{F}_{E, c} \big( \boldsymbol{\Theta}^{j,(k)}_{c}   \big) \big] \nonumber \\
                                &+ \Delta t \sum_{j = 1}^{m+1} (\tilde{q}^I_{m+1,j})_{c}              \big[ \boldsymbol{F}_{I, c} \big( \boldsymbol{\Theta}^{j,(k+1)}_{c} \big) 
                                                                                                - \boldsymbol{F}_{I, c} \big( \boldsymbol{\Theta}^{j,(k)}_{c}   \big) \big] \nonumber \\
                                &+ \Delta t \sum_{j = 0}^{M} (q_{m+1,j})_c \boldsymbol{F}_c \big( \boldsymbol{\Theta}^{j,(k)}_c \big) 
                                + \boldsymbol{\tau}^{m+1,(k)}_{c} \label{modified_update_equation_sdcq_discrete_form} .
\end{align}

\subsubsection{IMEX MLSDC algorithm}

We are now ready to review the steps of the MLSDC algorithm of \cite{emmett2012toward,speck2015multi} for the case of two space-time levels. 
In this section, $\boldsymbol{\Theta}^{m, (k)}_{\ell}$ denotes the approximate solution at temporal node $m$, space-time 
level $\ell$, 
and sweep $(k)$. $\vec{\boldsymbol{\Theta}}^{(k)}_{\ell}$ is the space-time vector that contains the approximate solution 
at all temporal nodes on level $\ell$.
The vectors $\boldsymbol{F}^{m, (k)}_{\ell}$ and $\vec{\boldsymbol{F}}^{(k)}_{\ell}$ are defined 
analogously. The MLSDC iteration starts with an SDC sweep on the fine level. The iteration continues as in a V-cycle from the 
fine level to the coarse level, and then back to the fine level. The specifics of the MLSDC iteration with two space-time 
levels are detailed in Algorithm~\oldref{alg:mlsdc_iteration}.

\begin{algorithm}[H]
  \SetAlgoLined
  \caption{\label{alg:mlsdc_iteration}IMEX MLSDC iteration on two space-time levels denoted by ``coarse'' and ``fine''.}
  \BlankLine

 \vspace{-0.2cm}
 \KwData{Initial data $\boldsymbol{\Theta}^{0,(k)}_{f}$ and function evaluations $\vec{\boldsymbol{F}}^{(k)}_{I, f}$, $\vec{\boldsymbol{F}}^{(k)}_{E, f}$ from the previous MLSDC 
 iteration $(k)$ on the fine level.}
 \KwResult{Approximate solution $\vec{\boldsymbol{\Theta}}^{(k+1)}_{\ell}$ and function evaluations $\vec{\boldsymbol{F}}^{(k+1)}_{I, \ell}$, $\vec{\boldsymbol{F}}^{(k+1)}_{E, \ell}$ on all levels.}
 
 \vspace{0.4cm}

 \textit{\textbf{A)} Perform a fine sweep} \\
 $\vec{\boldsymbol{\Theta}}^{(k+1)}_{f}, \, \vec{\boldsymbol{F}}^{(k+1)}_{I, f}, \, \vec{\boldsymbol{F}}^{(k+1)}_{E, f} \longleftarrow \textbf{SweepFine}\big( \vec{\boldsymbol{\Theta}}^{(k)}_{f}, \, \vec{\boldsymbol{F}}^{(k)}_{I, f}, \, \vec{\boldsymbol{F}}^{(k)}_{E, f} \big) $ \\
 
 \vspace{0.4cm}

   \textit{\textbf{B)} Restrict, re-evaluate, and save restriction} \\ 
   \For{$m = 1, \dots, M_c$}{
     $\boldsymbol{\Theta}^{m,(k)}_{c} \longleftarrow \textbf{Restrict} \big( \boldsymbol{\Theta}^{m,(k+1)}_{f} \big)$ \\
     $\boldsymbol{F}^{m,(k)}_{I, c}, \, \boldsymbol{F}^{m,(k)}_{E, c} \longleftarrow \textbf{Evaluate\_F} \big( \boldsymbol{\Theta}^{m,(k)}_{c} \big)$ \\
     $\boldsymbol{\tilde{\Theta}}^{m,(k)}_{c} \longleftarrow \boldsymbol{\Theta}^{m,(k)}_{c}$ \\
     $\boldsymbol{\tilde{F}}^{m,(k)}_{I,c}, \, \boldsymbol{\tilde{F}}^{m,(k)}_{E,c}  \longleftarrow \boldsymbol{F}^{m,(k)}_{I,c}, \, \boldsymbol{F}^{m,(k)}_{E,c}$
   }
   \textit{\textbf{C)} Compute FAS correction and sweep} \\ 
   $\boldsymbol{\tau}_{c} \longleftarrow \text{FAS} \big( \vec{\boldsymbol{F}}^{(k)}_{I, f}, \, \vec{\boldsymbol{F}}^{(k)}_{E, f}, \, \vec{\boldsymbol{F}}^{(k)}_{I, c}, \, \vec{\boldsymbol{F}}^{(k)}_{E, c}, \,  \boldsymbol{\tau}_{f} \big)$ \\
   $\vec{\boldsymbol{\Theta}}^{(k+1)}_{c}, \, \vec{\boldsymbol{F}}^{(k+1)}_{I, c}, \, \vec{\boldsymbol{F}}^{(k+1)}_{E, c} \longleftarrow \textbf{SweepCoarse} \big( \vec{\boldsymbol{\Theta}}^{(k)}_{c}, \, \vec{\boldsymbol{F}}^{(k)}_{I, c}, \, \vec{\boldsymbol{F}}^{(k)}_{E, c}, \, \boldsymbol{\tau}_{c} \big)$



 \vspace{0.4cm}
 
 \textit{\textbf{D)} Return to finest level before next iteration} \\
 \For{$m = 1, \dots, M_f$}{
   $\boldsymbol{\Theta}^{m,(k+1)}_{f} \longleftarrow \boldsymbol{\Theta}^{m,(k+1)}_{f} + \textbf{Interpolate} \big( \boldsymbol{\Theta}^{m, (k+1)}_{c} - \boldsymbol{\tilde{\Theta}}^{m, (k)}_{c} \big)$ \\
   $\boldsymbol{F}^{m, (k+1)}_{I, f} \longleftarrow \boldsymbol{F}^{m, (k+1)}_{I, f} + \textbf{Interpolate} \big( \boldsymbol{F}^{m, (k+1)}_{I, c} - \boldsymbol{\tilde{F}}^{m, (k)}_{I, c}\big)$ \\ 
   $\boldsymbol{F}^{m, (k+1)}_{E, f} \longleftarrow\boldsymbol{F}^{m, (k+1)}_{E, f} +  \textbf{Interpolate} \big( \boldsymbol{F}^{m, (k+1)}_{E, c} - \boldsymbol{\tilde{F}}^{m, (k)}_{E, c}\big)$ 
 }

\end{algorithm}

The single-level SDC iteration only consists of the fine sweep of Step $\boldsymbol{A}$. Both schemes share the same
initialization procedure. Specifically, before the first iteration, for $k=0$, we initialize the algorithm described above 
by simply copying the initial data for the time step, denoted by $\boldsymbol{\Theta}^{n,0}_f$, to all the other
SDC nodes, that is, for $m \in \{ 0, \dots, M_f \}$:
\begin{equation}
\boldsymbol{\Theta}^{m,(k = 0)}_f := \boldsymbol{\Theta}^{n,0}_f.
\end{equation}
In Algorithm~\oldref{alg:mlsdc_iteration}, the procedure \textit{SweepFine} consists in applying 
\ref{update_equation_sdcq_discrete_form} once on the fine level. The procedure \textit{SweepCoarse}
involves applying the correction described by \ref{modified_update_equation_sdcq_discrete_form}  
on the coarse level. We found that for the numerical examples considered in this work,
doing multiple sweeps on the coarse level instead of one every time \textit{SweepCoarse} is
called does not improve the accuracy of MLSDC, but increases the computational cost. This
is why the procedure \textit{SweepCoarse} only involves one sweep per call.
The procedure \textit{Evaluate\_F} involves computing the implicit and explicit right-hand sides. We 
highlight that the last step of Algorithm~\oldref{alg:mlsdc_iteration} does not involve any 
function evaluation. Instead, when we return to the fine level, we interpolate the coarse 
solution update as well as the coarse right-hand side corrections to the fine level. By 
avoiding $M_f$ function evaluations, this reduces the computational cost of the algorithm 
without undermining the order of accuracy of the scheme, as shown with numerical examples 
in Section~\oldref{section_numerical_examples}. 
We now discuss two key determinants of the performance of MLSDC-SH, namely, the coarsening 
strategy and the solver for the implicit systems.

\subsubsection{\label{subsubsection_coarsening_strategy_and_transfer_functions}Coarsening strategy and transfer functions}

In this section, we describe the linear restriction and interpolation operators used in the MLSDC-SH algorithm to transfer the approximate solution from fine to 
coarse levels, and vice-versa. 
In this work, the spatial restriction and interpolation procedures are performed in spectral space and heavily rely on the decomposition 
\ref{spectral_space} resulting from the SH basis. As explained below, this approach is based on the truncation of high-frequency modes, and 
therefore avoids the generation of spurious modes in the approximate solution that would propagate in the spectrum due to nonlinear wave 
interweaving over one coarse sweep. 

We reiterate that $K_{\ell}$ denotes the number of spectral coefficients used in \ref{spectral_space} 
at level $\ell$ and $M_{\ell}+1$ denotes the number of SDC nodes at level $\ell$. Therefore, the space-time vector storing the state of the system at level 
$\ell$, denoted by $\vec{\boldsymbol{\Theta}}_{\ell}$, is in $\mathbb{C}^{(M_{\ell}+1)K_{\ell}}$.
The two-step restriction process from fine level $\ell = f$ to coarse level $\ell = c$ consists in applying a restriction operator in time, denoted by 
$(\boldsymbol{R}^{t})^{c}_{f}$, followed by a restriction operator in space, denoted by $(\boldsymbol{R}^{s})^{c}_{f}$, that is,
\begin{equation}
\vec{\boldsymbol{\Theta}}_{c} = \boldsymbol{R}^{c}_{f} \vec{\boldsymbol{\Theta}}_{f} 
                                 = ( \boldsymbol{R}^{\textit{s}} )^{c}_{f} ( \boldsymbol{R}^{\textit{t}} )^{c}_{f} \vec{\boldsymbol{\Theta}}_{f}.
\label{two_step_restriction}
\end{equation}
In \ref{two_step_restriction}, the restriction operator in time is defined using the Kronecker product as
\begin{equation}
( \boldsymbol{R}^{\textit{t}} )^{c}_{f} \equiv \boldsymbol{\Pi}^{c}_{f} \otimes \boldsymbol{I}_{K_{f}} \in \mathbb{R}^{(M_{c}+1)K_{f} \times (M_{f}+1)K_{f}},
\label{time_restriction}
\end{equation}
where $\boldsymbol{I}_{K_{f}} \in \mathbb{R}^{K_{f} \times K_{f}} $ is the identity matrix, and 
$\boldsymbol{\Pi}^{c}_{f} \in \mathbb{R}^{(M_{c}+1) \times (M_{f}+1)}$ is the rectangle matrix employed to interpolate a scalar function 
from the fine temporal discretization to the coarse temporal discretization. Using the Lagrange polynomials 
$L^{j}_{f}$ on the fine temporal discretization, this matrix reads
\begin{equation}
(\boldsymbol{\Pi}^{c}_{f})_{ij} = L^{j-1}_{f}(t^{i-1}_{c}),
\end{equation}
using the SDC node $i-1$ at the coarse level, denoted by $t^{i-1}_{c}$. We note that, in the special case of two, three, and five Gauss-Lobatto nodes,
applying this restriction operator in time amounts to performing pointwise injection. The restriction operator in space consists in truncating the 
spectral representation of the primary variables \ref{spectral_space} based on the SH transform to remove the high-frequency features from 
the approximate solution. This is achieved by applying the matrix
\begin{equation}
( \boldsymbol{R}^s )^{c}_{f} \equiv \boldsymbol{I}_{M_{c}+1} \otimes \boldsymbol{D}^{c}_{f} \in \mathbb{R}^{(M_{c}+1)K_{c} \times (M_{c} + 1)K_{f} }.
\label{space_restriction}
\end{equation}
In \ref{space_restriction}, $\boldsymbol{D}^{c}_{f} \in \mathbb{R}^{K_{c} \times K_{f}}$ is a rectangle truncation matrix defined as
\begin{equation}
(\boldsymbol{D}^{c}_{f})_{ij} = 
\left\{
\begin{array}{l l}
   1 & i = j \\[4pt]
   0 & \text{otherwise.} 
\end{array} \right. 
\end{equation}

In the interpolation procedure employed to transfer the approximate solution from the coarse level to the fine level, we start 
with the application of the interpolation operator in space, $(\boldsymbol{P}^s)^{f}_{c}$, followed by the application of the interpolation 
operator in time, $(\boldsymbol{P}^t)^{f}_{c}$,
\begin{equation}
\vec{\boldsymbol{\Theta}}_{f} \equiv \boldsymbol{P}^{f}_{c} \vec{\boldsymbol{\Theta}}_{c} = (\boldsymbol{P}^t)^{f}_{c} (\boldsymbol{P}^{s})^{f}_{c} \vec{\boldsymbol{\Theta}}_{c}.
\end{equation}
The interpolation operator in space consists in padding the spectral representation of the primary variables at the coarse level 
with $K_{f} - K_{c}$  zeros and can be defined as the transpose of the restriction operator in space (see \ref{space_restriction}), 
that is,
\begin{equation}
( \boldsymbol{P}^s )^{f}_{c} \equiv \big( ( \boldsymbol{R}^s )^{c}_{f} \big)^T \in \mathbb{R}^{(M_{c}+1)K_{f} \times (M_{c}+1)K_{c}}.
\end{equation}
Finally, the interpolation operator in time is analogous to \ref{time_restriction} and reads 
\begin{equation}
( \boldsymbol{P}^t )^{f}_{c} \equiv \boldsymbol{\Pi}^{f}_{c} \otimes \boldsymbol{I}_{K_{f}} \in \mathbb{R}^{(M_{f}+1)K_{f} \times (M_{c}+1)K_{f}},
\end{equation}
where the rectangle interpolation matrix $\boldsymbol{\Pi}^{f}_{c}$ is constructed with the Lagrange polynomials 
$L^j_{c}$ on the coarse temporal discretization. For two, three, and five Gauss-Lobatto nodes, this amounts to performing
pointwise injection at the fine nodes that correspond to the coarse nodes, and then polynomial interpolation to compute the 
solution at the remaining fine nodes. This completes the presentation of the MLSDC-SH algorithm for the time integration of the 
shallow-water equations on the rotating sphere. Next, we discuss the implicit solver used in this work.

\subsection{\label{implicit_solver_for_sdc_and_mlsdc}Solver for the implicit systems in SDC and MLSDC}

The time integration schemes of Sections~\oldref{subsection_implicit_explicit_spectral_deferred_correction} and 
\oldref{subsection_multi_level_spectral_deferred_correction} involve solving implicit linear systems in the form 
\begin{equation}
\boldsymbol{\Theta}^{m+1,(k+1)} - \Delta t \tilde{q}^I_{m+1,m+1} \boldsymbol{F}_I( \boldsymbol{\Theta}^{m+1,(k+1)} ) = \boldsymbol{b},
\label{implicit_systems}
\end{equation}
where $\boldsymbol{b}$ is obtained from \ref{update_equation_sdcq_discrete_form} or \ref{modified_update_equation_sdcq_discrete_form},
and where we have dropped the subscripts denoting the space-time levels for simplicity. The structure of the implicit linear systems results 
from the spatial discretization based on the SH transform, but also from the temporal splitting between implicit and explicit terms 
described in Section~\oldref{subsection_temporal_splitting}. The solution strategy for \ref{implicit_systems} is performed in spectral space 
and follows two steps briefly outlined below. 

First, we algebraically form a reduced linear system containing only the geopotential 
unknowns -- that is, $K/3$ degrees of freedom, where $K$ denotes the total number of spectral coefficients needed to represent
the three primary variables in \ref{spectral_space}. Given that the longitudinal and latitudinal coupling terms present in the Coriolis 
term and in the nonlinear operators are discretized explicitly, the $K/3$ geopotential degrees of freedom are fully decoupled from 
one another. The geopotential linear system is therefore diagonal and trivial to solve. 
Second, we have to solve for the remaining $2K/3$ vorticity and divergence degrees of freedom. This is again a trivial 
operation that does not require a linear solver, since we have to solve two diagonal linear systems to update the vorticity 
and divergence variables, respectively.

Therefore, solving the implicit system \ref{implicit_systems} is purely based on local operations during which the degrees of 
freedom are updated one at a time in spectral space. We refer to \cite{schreiber2018sph} for the detailed formulation of 
the geopotential, vorticity, and divergence diagonal linear systems.


\subsection{\label{subsection_computational_cost_of_sdc_and_mlsdc}Computational cost of SDC and MLSDC}

In this section, we compare the computational cost of the MLSDC-SH scheme
described in Section \oldref{subsection_multi_level_spectral_deferred_correction} to that of the single-level SDC scheme.
We refer to the single-level SDC scheme with $M_f+1$ temporal nodes and $N_S$ fine sweeps as SDC($M_f+1$,$N_S$).
We denote by MLSDC($M_f+1$, $M_c+1$, $N_{ML}$, $\alpha$) the MLSDC-SH scheme with $M_f+1$ nodes on the fine level, 
$M_c + 1$ nodes of the coarse level, $N_{\textit{ML}}$ iterations, and a spatial coarsening ratio, $\alpha$, defined using \ref{spectral_to_physical_space} as 
$\alpha = R_{c}/R_{f}$.
The parameters of the SDC and MLSDC-SH schemes are summarized in Tables~\oldref{tbl:sdc_parameters} and \oldref{tbl:mlsdc_overview},
respectively. 

\begin{table}[!ht]
\begin{center}
\begin{tabular}{|c|c|}
  \hline
  \multicolumn{2}{|c|}{SDC($M_f+1$, $N_{\textit{S}}$)} \\  
  \hline
  \textbf{Parameter} & \multicolumn{1}{c|}{\textbf{Description}} \\
  \hline
  $M_f+1$	     & SDC nodes on fine level\\
  \hline
  $N_{S}$          & Number of SDC iterations  \\
  \hline
\end{tabular}
\end{center}
\vspace{-0.5cm}
\caption{\label{tbl:sdc_parameters}Parameters for the SDC scheme. The SDC iteration only involves one sweep on the fine level.}
\end{table}

\begin{table}[!ht]
\begin{center}
\begin{tabular}{|c|c|}
  \hline
  \multicolumn{2}{|c|}{MLSDC($M_f+1$, $M_c+1$, $N_{\textit{ML}}$, $\alpha$)} \\  
  \hline
  \textbf{Parameter} & \multicolumn{1}{c|}{\textbf{Description}} \\
  \hline
  $M_f+1$	     & SDC nodes on fine level\\
  \hline
  $M_c+1$   	     & SDC nodes on coarse level\\
  \hline
  $N_{ML}$          & Number of MLSDC iterations  \\
  \hline
  $\alpha$           & Spatial coarsening ratio \\
  \hline
\end{tabular}
\end{center}
\vspace{-0.5cm}
\caption{\label{tbl:mlsdc_overview}Parameters for the MLSDC-SH scheme. The MLSDC-SH iteration is described in 
Alg.~\oldref{alg:mlsdc_iteration}, and involves one sweep on the fine level and one sweep on the coarse level.}
\end{table}

To evaluate the theoretical computational cost of the MLSDC-SH scheme, we count the number of 
function evaluations and the number of solves involved in a time step. We denote by $C^s_{\ell}$ 
the cost of a solve at level $\ell$, and by $C^{\textit{fi}}_{\ell}$ (respectively, $C^{\textit{fe}}_{\ell}$) 
the cost of an implicit (respectively, explicit) function evaluation at level $\ell$. We neglect the 
cost of computing the FAS correction. The quantities $C^s_{c}$, $C^{\textit{fi}}_{c}$, $C^{\textit{fe}}_{c}$ 
depend on the spatial coarsening ratio, $\alpha = R_c/R_f$. Here, $R_{\ell}$ represents the highest Fourier wavenumber in 
the east-west representation of \ref{spectral_to_physical_space} on level $\ell$.

The cost of a time step with the two-level MLSDC($M_f+1$,$M_c+1$,$N_{\textit{ML}}$,$\alpha$) is
\begin{align}
C^{\textit{MLSDC}(M_f+1,M_c+1,N_{\textit{ML}},\alpha)} &= N_{\textit{ML}} M_{f} ( C^s_{f} + C^{\textit{fi}}_{f} + C^{\textit{fe}}_{f} ) \nonumber \\
                                            &+ N_{\textit{ML}} M_{c} ( C^s_{c} + C^{\textit{fi}}_{c} + C^{\textit{fe}}_{c} )  \nonumber \\
                                            &+ N_{\textit{ML}} M_{c} ( C^{\textit{fi}}_{c} + C^{\textit{fe}}_{c}), 
\label{computational_cost_mlsdc}
\end{align}
where the term in the right-hand side of the first line represents the cost of the fine sweeps, the second term represents
the cost of the coarse sweeps, and the third term accounts for the cost of the evaluation of the right-hand 
sides at the coarse nodes after the restriction. This can be compared with the cost of a time 
step in SDC($M_f+1$,$N_{\textit{S}}$),  given by 
\begin{equation}
C^{\textit{SDC}(M_f+1,N_{\textit{S}})} = N_{\textit{S}} M_f ( C^s_{f} + C^{\textit{fi}}_{f} + C^{\textit{fe}}_{f} ).
\label{computational_cost_sdc}
\end{equation}
Furthermore, we assume that the cost of a linear solve is the same 
as the cost of evaluating the right-hand side, that is,
\begin{equation}
C^s_{f} = C^{\textit{fi}}_{f} = C^{\textit{fe}}_{f},
\end{equation}
with this assumption being motivated by Section~\oldref{implicit_solver_for_sdc_and_mlsdc}. 
In addition, we will also assume that the 
computational cost of the operators is proportional to the number of spectral coefficients in \ref{spectral_space}, denoted by $K_\ell$.
There are three primary variables, and each of them is represented in the triangular truncation framework with $R_{\ell}(R_{\ell}+1)/2$ 
spectral coefficients, where $R_{\ell}$ denotes the highest Fourier wavenumber in the east-west representation of \ref{spectral_to_physical_space}.
This yields $K_{\ell} = 3 R_{\ell}(R_{\ell}+1) / 2$. Using this notation and the definition $\alpha = R_c/R_f$, we can obtain an expression 
of $C^s_c$ as a function $\alpha$, $C^s_f$, and $R_f$, by writing
\begin{equation}
C^s_c = \frac{K_c}{K_f} C^s_f = \frac{R_c(R_c+1)}{R_f(R_f+1)} C^s_f = \alpha^2 \frac{R_f + 1/\alpha}{R_f + 1} C^s_f.
\end{equation}
Similarly, using the same assumptions, we obtain
\begin{equation}
C^{fi}_c = \alpha^2 \frac{R_f+1/\alpha}{R_f+1} C^{fi}_f, \qquad \qquad C^{fe}_c = \alpha^2 \frac{R_f+1/\alpha}{R_f+1} C^{fe}_f.
\end{equation}
Using these notations and assuming that MLSDC-SH and SDC achieve the same accuracy, the theoretical speedup obtained 
with MLSDC-SH, denoted by $\mathcal{S}^{\textit{theo}}$, reads
\begin{equation}
\mathcal{S}^{\textit{theo}} = \frac{C^{\textit{SDC}(M_f+1,N_{\textit{S}})}}{C^{\textit{MLSDC}(M_f+1,M_c+1,N_{\textit{ML}},\alpha)}} 
            = \frac{N_S}{N_{ML}} \times \frac{1}{1 + \displaystyle  \alpha^2  \frac{5 (R_{f} + 1/\alpha) M_c}{3 (R_{f}+1) M_f }}.
\label{theoretical_speedup}
\end{equation}
That is, MLSDC(3,2,2,1/2), based on two fine sweeps and two coarse sweeps, yields a theoretical
speedup $\mathcal{S}^{\textit{theo}} \approx 1.66$ compared to SDC(3,4), which uses four fine sweeps.
This corresponds to a reduction of 40 \% in the wall-clock time.
MLSDC(5,3,4,1/2), based on four fine sweeps and four coarse sweeps, also results in a theoretical 
speedup $\mathcal{S}^{\textit{theo}} \approx 1.66$ compared to SDC(5,8), which relies on eight fine sweeps. 
This reasoning assumes that one MLSDC-SH iteration can replace two single-level SDC iterations
and still achieve the same accuracy. This point is investigated in the next section using numerical examples.

We conclude this section by comparing the computational cost of MLSDC-SH to that of the fully implicit single-level
SDC scheme based on the Spectral Element Method (SEM) presented in \cite{jia2013spectral} and referred to as SDC-SEM
in the remainder of this paper. The fully implicit SDC-SEM iteration entails solving a large nonlinear system
on the fine problem to update all the degrees of freedom simultaneously. 
This is achieved using the Jacobian-free Newton-Krylov (JFNK) method. As stated by the authors, the computational
cost of this nonlinear solve can be large and heavily depends on the availability of a scalable preconditioner
for the linear systems. Instead, the MLSDC-SH iteration involves trivial diagonal linear solves that can easily be 
parallelized. In addition, our multi-level time integration framework relies on a hierarchy of space-time levels to shift
a significant fraction of the computational work to the coarser representation of the problem. This reduces the number
of fine sweeps in the algorithm and therefore further reduces the cost of a time step.

As a result, we expect IMEX MLSDC-SH to be significantly less expensive than the fully implicit SDC-SEM on a per-timestep basis for moderate resolutions. 
Assuming that the linear systems can be efficiently preconditioned -- as in \cite{lott2015algorithmically} -- the key to the performance
of SDC-SEM lies in its ability to take much larger stable time steps than MLSDC-SH to compensate for its relatively
high cost on a per-timestep basis. Exploration of this trade-off requires a careful analysis that will be presented in future work.

\section{\label{section_numerical_examples}Numerical examples}

We assess the performance of MLSDC-SH with state-of-the-art test cases for the development of dynamical cores. All the 
test cases are nonlinear. They are selected to focus on particular challenges that arise with MLSDC-SH. The 
first test case in Section~\oldref{subsection_steady_zonal_jet} targets geostrophically balanced modes. It evaluates 
the effects of multi-level mode truncation and the relation to the diffusion used in the simulations. The 
second test case in Section~\oldref{subsection_nonlinear_propagation_of_gaussian_dome} studies the observed order of 
convergence of MLSDC-SH  upon refinement in time for waves propagating on the rotating sphere. While the first two test cases 
are mainly dominated by the linear parts, the following benchmarks assess the performance of MLSDC-SH in the presence 
of stronger nonlinear interactions. The Rossby-Haurwitz benchmark in Section~\oldref{subsection_rossby_haurwitz_wave} 
studies the advection of a wave that propagates around the sphere without changing shape.
This is followed by the unstable barotropic wave benchmark in Section~\oldref{subsection_galewsky} with an initially 
linear balanced flow perturbed by the introduction of a Gaussian bump in the geopotential field. All these benchmarks 
provide a key insight into the numerical properties of MLSDC-SH in the context of atmospheric simulations.


\subsection{\label{subsection_steady_zonal_jet}Steady zonal jet}

\begin{figure}[ht!]
\centering
\begin{tikzpicture}
\node[anchor=south west,inner sep=0] at (0,0){\includegraphics[scale=0.27]{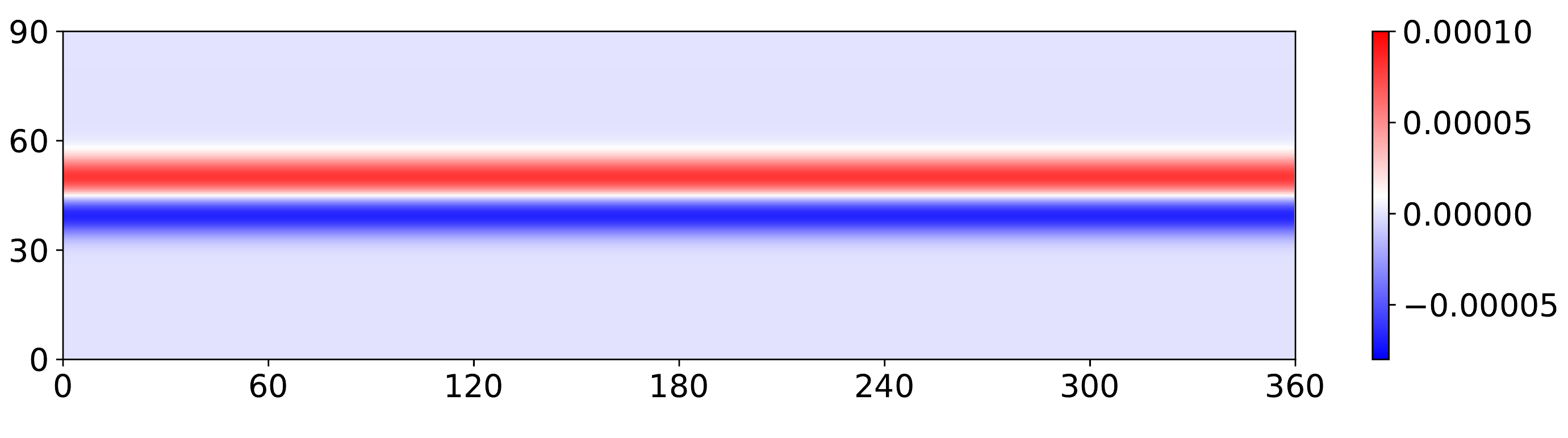}};
\node (ib_1) at (3.4,-0.05) {\scriptsize Longitude (degrees)};
\node[rotate=90] (ib_1) at (-0.2,1.2) {\scriptsize Latitude (degrees)};
\end{tikzpicture}
\vspace{-0.4cm}
\caption{\label{fig:galewsky_test_case_without_localized_bump_vorticity_field} 
Steady zonal jet: vorticity field with a resolution of $R_f = S_f = 256$ after 144 
hours. This solution is obtained with the single-level SDC(5,8). The diffusion coefficient is 
$\nu_{\mathcal{B}} = 1.0 \times 10^{5} \, \text{m}^2.\text{s}^{-1}$.
}
\end{figure}

\begin{figure}[ht!]
\centering
\begin{tikzpicture}
\node[anchor=south west,inner sep=0] at (0,0){\includegraphics[scale=0.365]{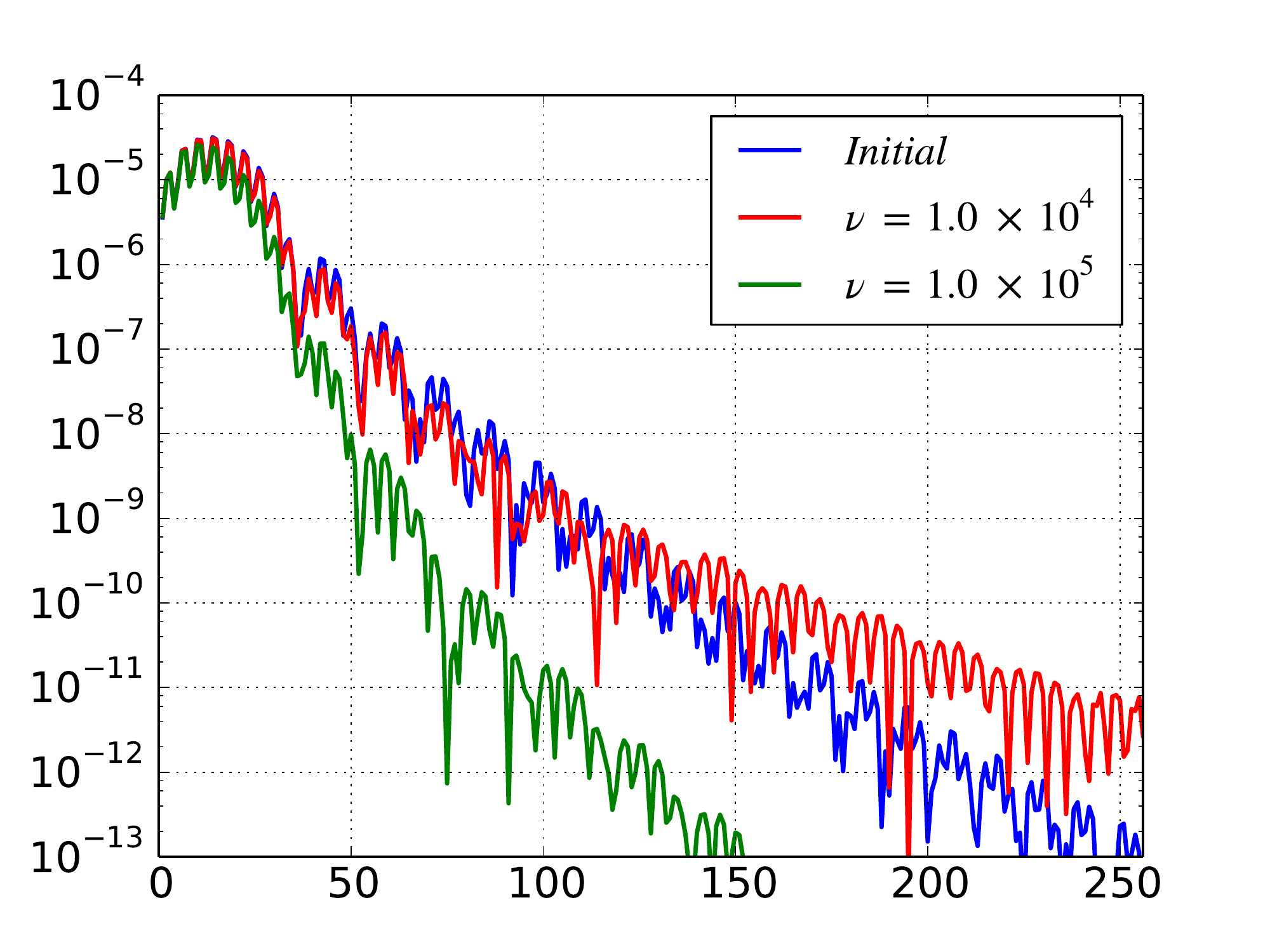}};
\node (ib_1) at (3.8,-0.05) {$n_0$};
\node[rotate=90] (ib_1) at (0,2.7) {$|\zeta_{n_0}|$};



\node (ib_1) at (1.65,5.14) {};
\node (ib_2) at (1.65,0.45) {};
\path [draw=black!100, thick, dashed] (ib_1) -- (ib_2);

\node (ib_1) at (2.35,5.14) {};
\node (ib_2) at (2.35,0.45) {};
\path [draw=black!100, thick, dashed] (ib_1) -- (ib_2);

\node (ib_1) at (3.75,5.14) {};
\node (ib_2) at (3.75,0.45) {};
\path [draw=black!100, thick, dashed] (ib_1) -- (ib_2);

\node[rotate=-90] (ib_1) at (1.8,1.425) {\scriptsize $\alpha = 1/8$};
\node[rotate=-90] (ib_2) at (2.5,4.38)  {\scriptsize $\alpha = 1/4$};
\node[rotate=-90] (ib_3) at (3.9,3.42) {\scriptsize $\alpha = 1/2$}; 

\end{tikzpicture}
\vspace{-0.4cm}
\caption{\label{fig:galewsky_test_case_no_localized_bump_vorticity_spectrum} 
Steady zonal jet: max-spectrum of the vorticity field at the beginning of the simulation and after 144 hours for different 
values of the diffusion coefficient. The quantity on the $y$-axis is defined as $|\zeta_{n_0}| = \max_{r} |\zeta^r_{n_0}|$.
}
\end{figure}
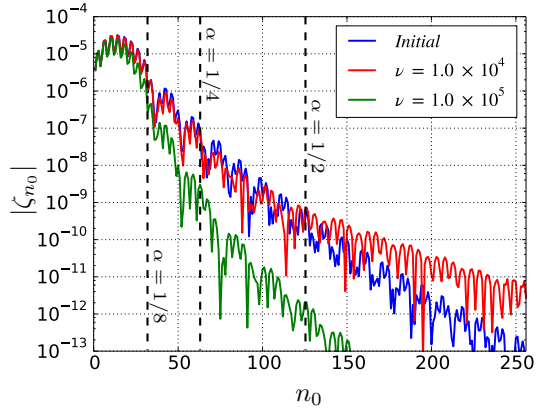

We first study the behavior of the multi-level SDC scheme on a steady test case derived from \cite{galewsky2004initial}. 
This test case consists in the simulation of a steady, analytically specified mid-latitude jet with an unperturbed, balanced height field. 
This test assesses the ability of the numerical schemes to maintain this balanced state for 144 hours. The vorticity field obtained 
with the single-level SDC(5,8) with a modal resolution of $R_f = S_f = 256$ is in 
Fig.~\oldref{fig:galewsky_test_case_without_localized_bump_vorticity_field}, along with the corresponding vorticity spectrum in 
Fig.~\oldref{fig:galewsky_test_case_no_localized_bump_vorticity_spectrum}. 

This steady numerical test is used to illustrate the order of convergence of MLSDC-SH upon refinement in time. We will consider 
the computational cost of the multi-level scheme in subsequent examples. We highlight that we do not use the steady geostrophic 
balance test case of \cite{williamson1992standard} here because it is based on an initial vorticity field that can be represented 
with a few modes only. Instead, the vorticity field of Fig.~\oldref{fig:galewsky_test_case_no_localized_bump_vorticity_spectrum} 
has a spectrum that spans a larger number of modes. This is key for our analysis because it allows us to better study the impact 
of the coarsening strategy based on spectral coefficient truncation (see Section~\oldref{subsubsection_coarsening_strategy_and_transfer_functions})
on the convergence rate upon temporal refinement. To our best knowledge, this is the first time that a study of the effect of the 
coarsening strategy on the observed order of convergence of MLSDC is conducted.

\begin{figure}[ht!]
\centering
\subfigure[]{
\begin{tikzpicture}
\node[anchor=south west,inner sep=0] at (0,0){\includegraphics[scale=0.365]{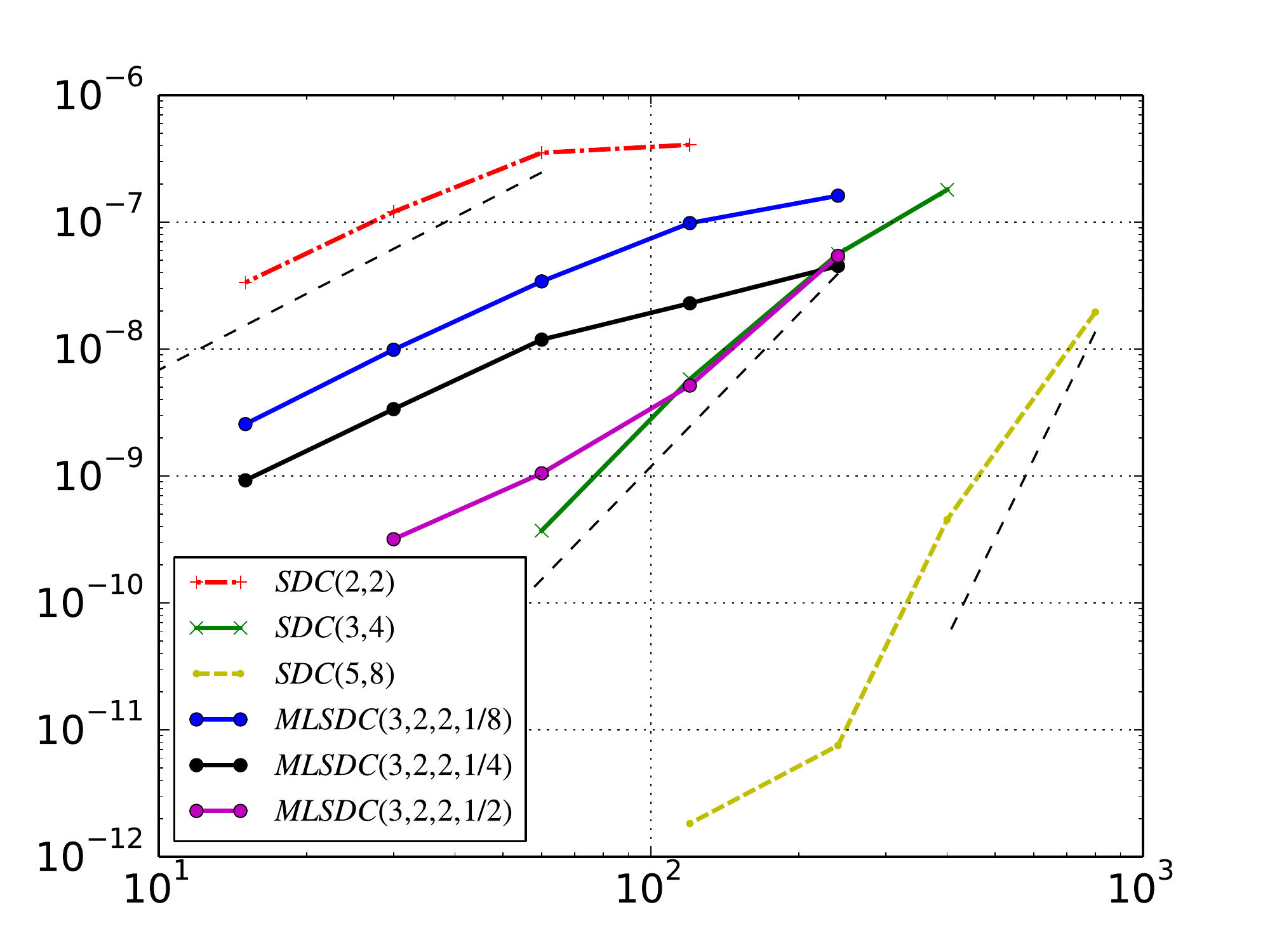}};
\node (ib_1) at (3.8,5.3) {\small $\nu_{\mathcal{A}} = 1.0 \times 10^{4} \, \text{m}^2.\text{s}^{-1}$};
\node (ib_1) at (3.8,-0.05) {$\Delta t$};
\node[rotate=90] (ib_1) at (0,2.8) {$L_{\infty}$-norm of the error};
\end{tikzpicture}
\label{fig:galewsky_test_case_no_localized_bump_accuracy_vorticity_A} 
}
\subfigure[]{
\begin{tikzpicture}
\node[anchor=south west,inner sep=0] at (0,0){\includegraphics[scale=0.365]{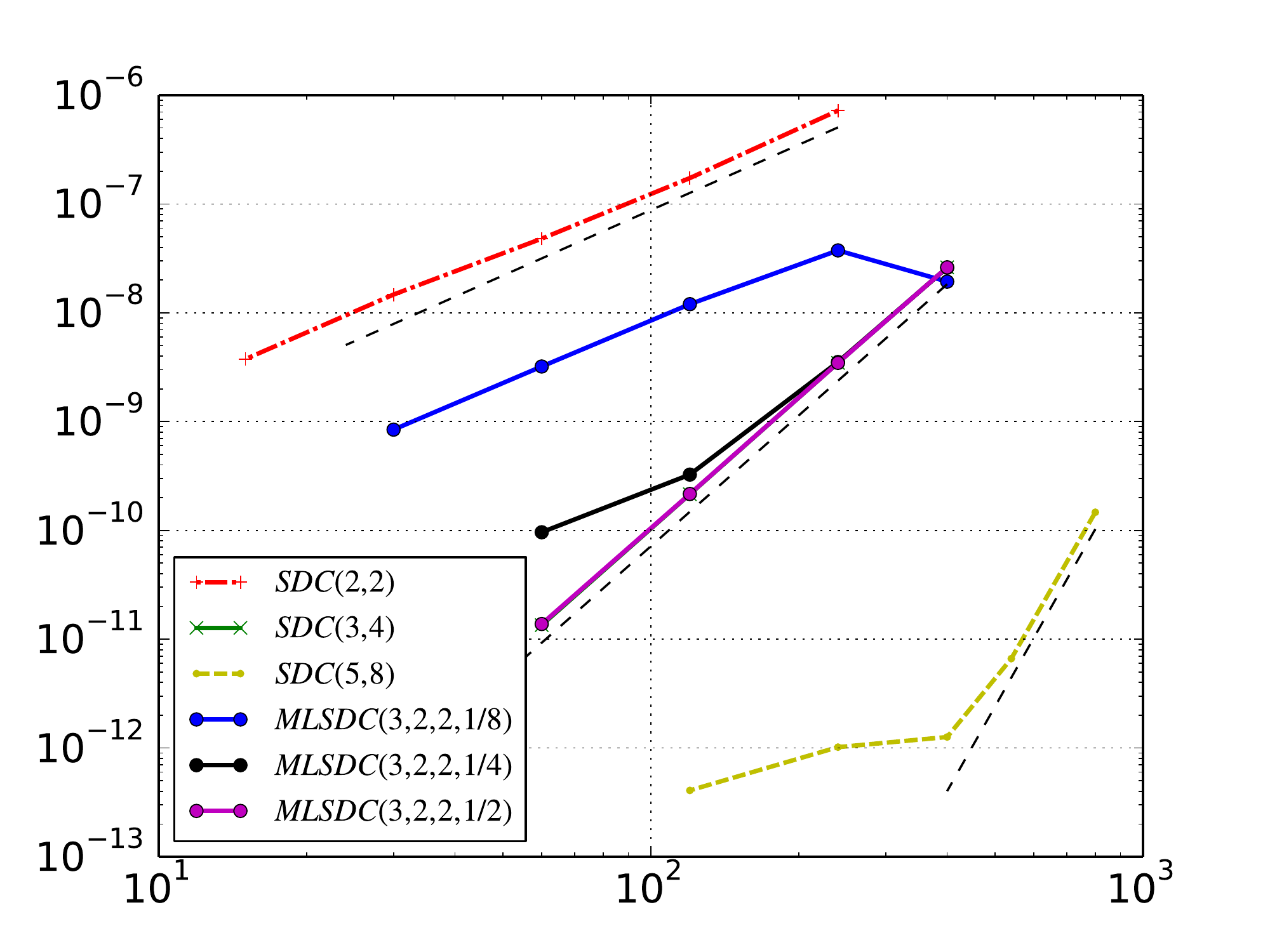}};
\node (ib_1) at (3.8,5.3) {\small $\nu_{\mathcal{B}} = 1.0 \times 10^{5} \, \text{m}^2.\text{s}^{-1}$};
\node (ib_1) at (3.8,-0.05) {$\Delta t$};
\node[rotate=90] (ib_1) at (0,2.8) {$L_{\infty}$-norm of the error};
\end{tikzpicture}
\label{fig:galewsky_test_case_no_localized_bump_accuracy_vorticity_B} 
}
\vspace{-0.5cm}
\caption{\label{fig:galewsky_test_case_no_localized_bump_accuracy_vorticity} 
Steady zonal jet: $L_{\infty}$-norm of the error in the vorticity field with respect to the reference solution as a function 
of time step size. When the norm of the error is smaller than the magnitude of the spectral terms truncated during spatial 
coarsening (given in Fig.~\oldref{fig:galewsky_test_case_no_localized_bump_vorticity_spectrum}), MLSDC(3,2,2,$\alpha$) 
exhibits only second-order convergence upon refinement in time. Above this threshold, we observe fourth-order convergence.
}
\end{figure}

We perform a refinement study in time to assess the impact of the spatial coarsening ratio, $\alpha$, on the observed order of convergence of 
the MLSDC-SH scheme. The study is done with a fixed fine resolution of $R_f = S_f = 256$. We consider two configurations, denoted by 
$\mathcal{A}$ and $\mathcal{B}$, in which the diffusion coefficient is set to $\nu_{\mathcal{A}} = 1.0 \times 10^4 \, \text{m}^2. \text{s}^{-1}$ 
and $\nu_{\mathcal{B}} = 1.0 \times 10^5 \, \text{m}^2. \text{s}^{-1}$, respectively. In each configuration, the reference solution is obtained with 
SDC(5,8) using a time step size of $\Delta t_{\textit{ref}} = 90 \, \text{s}$. Fig.~\oldref{fig:galewsky_test_case_no_localized_bump_accuracy_vorticity} 
shows the norm of the error in the vorticity field with respect to the reference solution as a function of the time step size. 
We use the $L_{\infty}$-norm to illustrate the connection between the convergence rate of MLSDC-SH upon refinement in time and the magnitude 
of the spectral coefficients of the vorticity that are truncated during the spatial restriction from the fine level to the coarse level. 
For this test case, we focus on MLSDC(3,2,2,$\alpha$), which relies on three fine temporal nodes, two coarse temporal nodes, and uses two 
iterations (with one fine sweep and one coarse sweep per iteration). In both configurations, the observed order of convergence of MLSDC(3,2,2,$\alpha$) 
varies significantly as a function of the spatial coarsening ratio, $\alpha$. 

In configuration $\mathcal{A}$, MLSDC(3,2,2,1/2) achieves fourth-order convergence upon refinement in time for stable time steps larger than
$120 \, \text{s}$, but exhibits only second-order convergence for shorter time steps. The reduction in the observed 
order of convergence can be explained by considering the vorticity spectrum of Fig.~\oldref{fig:galewsky_test_case_no_localized_bump_vorticity_spectrum}. 
Specifically, we note that for the time step size range defined by $\Delta t \leq 90 \, \text{s}$, the MLSDC(3,2,2,1/2) scheme reaches 
an $L_{\infty}$-norm of the error smaller than $10^{-9}$. We see in Fig.~\oldref{fig:galewsky_test_case_no_localized_bump_vorticity_spectrum} 
that this threshold corresponds to the order of magnitude of the truncated terms during the restriction to the coarse level when 
$R_c = S_c = 128$. 
In MLSDC(3,2,2,1/4) and MLSDC(3,2,2,1/8), the truncated coefficients in the vorticity spectrum are relatively large which causes the observed 
convergence rate upon refinement in time to be reduced to second order in the entire range of stable time step sizes. Conversely, 
MLSDC(3,2,2,4/5) achieves the same convergence rate as SDC(3,4) over the full time step range (not shown here for brevity).

In configuration $\mathcal{B}$, the use of a larger diffusion coefficient significantly reduces the magnitude of the spectral coefficients 
associated with the high-frequency modes. Therefore, the MLSDC(3,2,2,1/2) scheme achieves fourth-order convergence in the entire time step range 
considered here. MLSDC(3,2,2,1/4) achieves fourth-order convergence in a larger fraction of the range of stable time step sizes, but still exhibits 
a reduction of its observed order of convergence when the norm of the error reaches the magnitude of the terms that are truncated during the restriction 
procedure. MLSDC(3,2,2,1/8) is still limited to second-order convergence. For this test case, the order of convergence of the MLSDC-SH scheme 
for the geopotential and divergence variables, not shown here, is similar to that observed for the vorticity. 

The key insight of this section is that the accuracy of MLSDC-SH as the time step is reduced depends 
on the interplay between two key factors, namely the spectrum of the fine solution and the magnitude of the spatial coarsening 
ratio. They determine the range of scales of the fine solution that can be captured by the coarse correction, and as a result have a 
strong impact on the observed order of convergence of MLSDC-SH upon refinement in time. In particular, the presence of large
high-frequency modes in the fine solution imposes of lower limit on the spatial coarsening ratio to preserve the high-order
convergence of the multi-level scheme. Next, we study the computational 
cost of the MLSDC-SH scheme using unsteady test cases, starting with Gaussian dome propagation.

\subsection{\label{subsection_nonlinear_propagation_of_gaussian_dome}Propagation of a Gaussian dome}

\begin{figure}[ht!]
\centering
\subfigure[]{
\begin{tikzpicture}
\node[anchor=south west,inner sep=0] at (0,0){\includegraphics[scale=0.267]{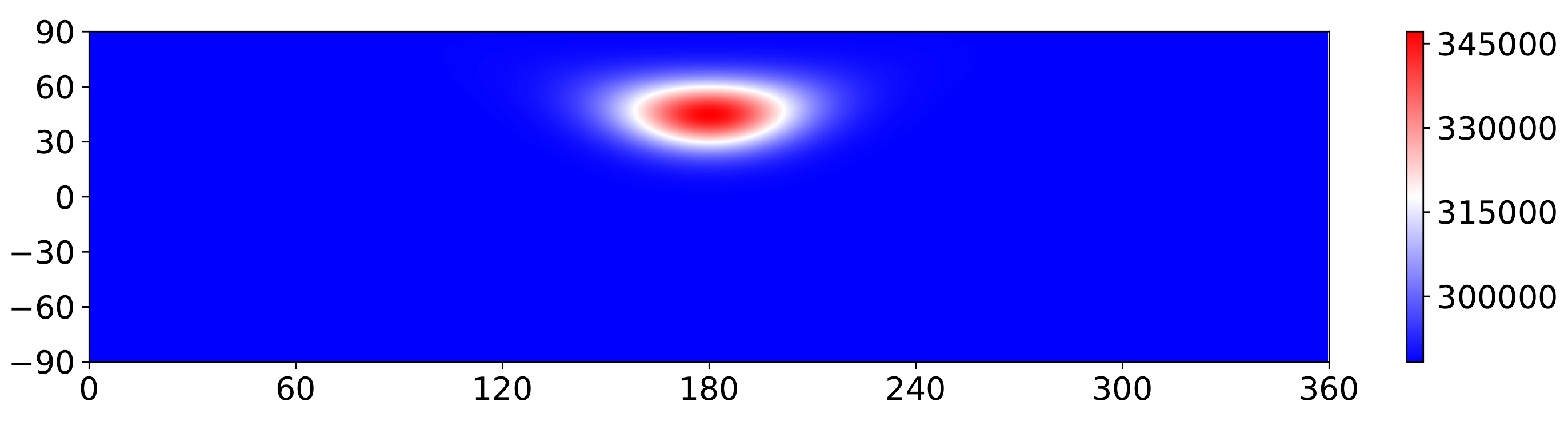}};
\node (ib_1) at (3.4,-0.05) {\scriptsize Longitude (degrees)};
\node[rotate=90] (ib_1) at (-0.2,1.2) {\scriptsize Latitude (degrees)};
\end{tikzpicture}
\label{fig:gaussian_bump_test_case_geopotential_field_0}
}
\hspace{-0.5cm}
\subfigure[]{
\begin{tikzpicture}
\node[anchor=south west,inner sep=0] at (0,0){\includegraphics[scale=0.267]{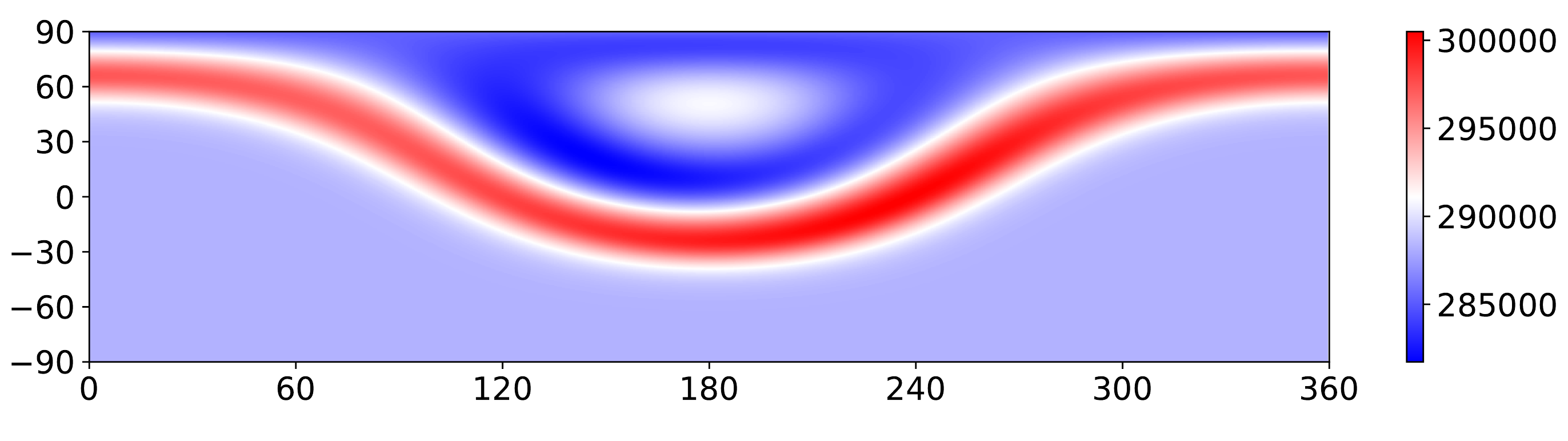}};
\node (ib_1) at (3.4,-0.05) {\scriptsize Longitude (degrees)};
\end{tikzpicture}
\label{fig:gaussian_bump_test_case_geopotential_field_2}
}  \\\vspace{-0.3cm}
\subfigure[]{
\begin{tikzpicture}
\node[anchor=south west,inner sep=0] at (0,0){\includegraphics[scale=0.267]{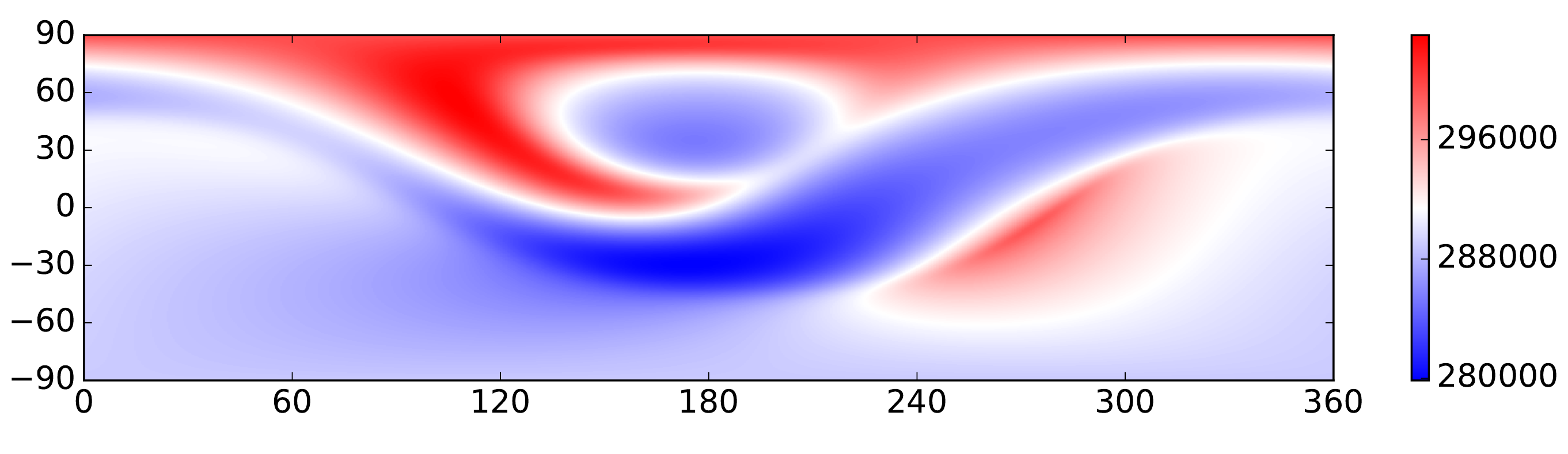}};
\node (ib_1) at (3.4,-0.05) {\scriptsize Longitude (degrees)};
\node[rotate=90] (ib_1) at (-0.2,1.2) {\scriptsize Latitude (degrees)};
\end{tikzpicture}
\label{fig:gaussian_bump_test_case_geopotential_field_4}
} 

\vspace{-0.45cm}
\caption{\label{fig:gaussian_bump_test_case_geopotential_field} 
  Gaussian dome: geopotential field with a resolution of $R_f = S_f = 256$ at the start of the simulation
  in~\oldref{fig:gaussian_bump_test_case_geopotential_field_0}, after $9000 \, \text{s}$
  in~\oldref{fig:gaussian_bump_test_case_geopotential_field_2},
  and after one day in \oldref{fig:gaussian_bump_test_case_geopotential_field_4}. 
  This solution is obtained with the single-level SDC(5,8). The diffusion coefficient
  is $\nu = 1.0 \times 10^{5} \, \text{m}^2.\text{s}^{-1}.$
}
\end{figure}

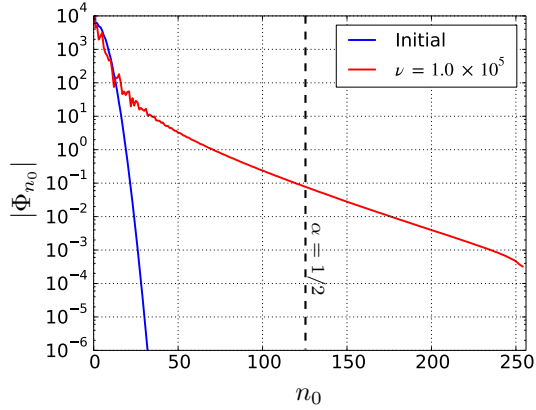
\begin{figure}[ht!]
\centering
\begin{tikzpicture}
\node[anchor=south west,inner sep=0] at (0,0){\includegraphics[scale=0.365]{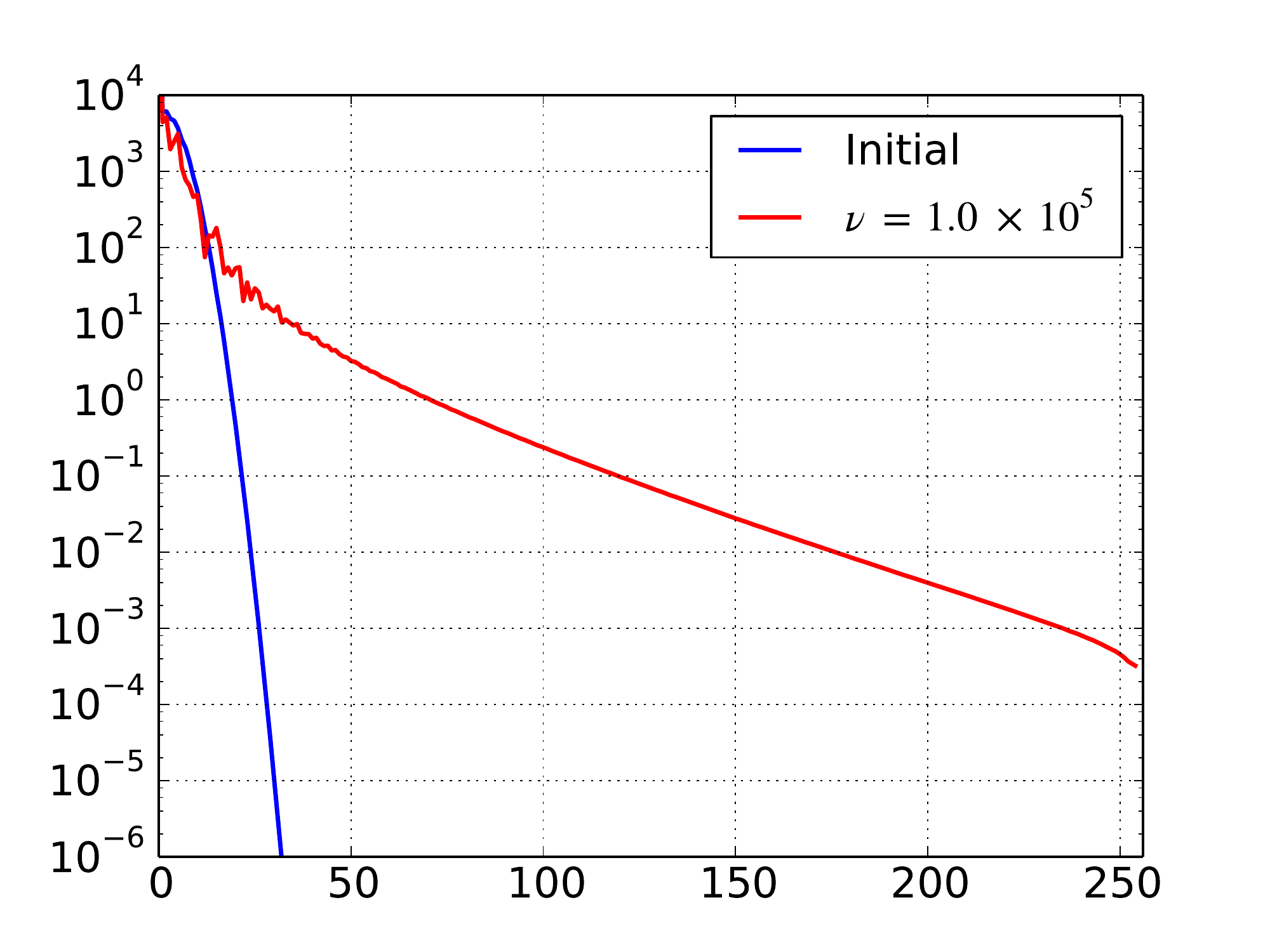}};
\node (ib_1) at (3.8,-0.05) {$n_0$};
\node[rotate=90] (ib_1) at (0,2.7) {$|\Phi_{n_0}|$};

\node (ib_1) at (3.75,5.14) {};
\node (ib_2) at (3.75,0.45) {};
\path [draw=black!100, thick, dashed] (ib_1) -- (ib_2);

\node[rotate=-90] (ib_3) at (3.9,1.78) {\scriptsize $\alpha = 1/2$};

\end{tikzpicture}
\vspace{-0.4cm}
\caption{\label{fig:gaussian_test_case_spectra} 
Gaussian dome: max-spectrum of the geopotential field at the beginning of the simulation and after one day. 
We use a diffusion coefficient $\nu = 1.0 \times 10^{5} \, \text{m}^2.\text{s}^{-1}$. The quantity on the 
$y$-axis is defined as $|\Phi_{n_0}| = \max_{r} |\Phi^r_{n_0}|$.
}
\end{figure}

We now consider an initial condition derived from the third numerical experiment of \cite{swarztrauber2004shallow}. The velocities are initially 
equal to zero ($u = v = 0$). We place a Gaussian dome in the initial geopotential field, such that
\begin{equation}
h( \lambda, \phi ) = \bar{h} + A \text{e}^{-\alpha (d / a )^2}, \label{gaussian_bump}
\end{equation}
where $a$ denotes the Earth radius. The distance $d$ is defined as
\begin{equation}
d = \sqrt{x^2 + y^2 + z^2},
\end{equation}
with 
\begin{align}
x &= a \big( \cos(\lambda) \cos(\phi) - \cos(\lambda_c) \cos(\phi_c) \big), \\
y &= a \big( \sin(\lambda) \cos(\phi) - \sin(\lambda_c) \cos(\phi_c) \big), \\
z &= a \big( \sin(\phi) - \sin(\phi_c) \big).
\end{align} 
This corresponds to a Gaussian dome centered at $\lambda_c = \pi$ and $\phi_c = \pi/4$. 
We use realistic values for the Earth radius, the gravitational acceleration,
and the angular rate of rotation $\Omega$ involved in the Coriolis force.
We set $\bar{h} = 29400 \, \text{m}$ and $A = 6000 \, \text{m}$, which is about ten times larger than in the original test case. 
The simulation of the collapsing dome is run for one day to study the behavior of MLSDC-SH on an advection-dominated test case.  
The geopotential field at different times is in Fig.~\oldref{fig:gaussian_bump_test_case_geopotential_field}, and the spectrum is in 
Fig.~\oldref{fig:gaussian_test_case_spectra}.

\begin{figure}[ht!]
\centering
\subfigure[]{
\begin{tikzpicture}
\node[anchor=south west,inner sep=0] at (0,0){\includegraphics[scale=0.365]{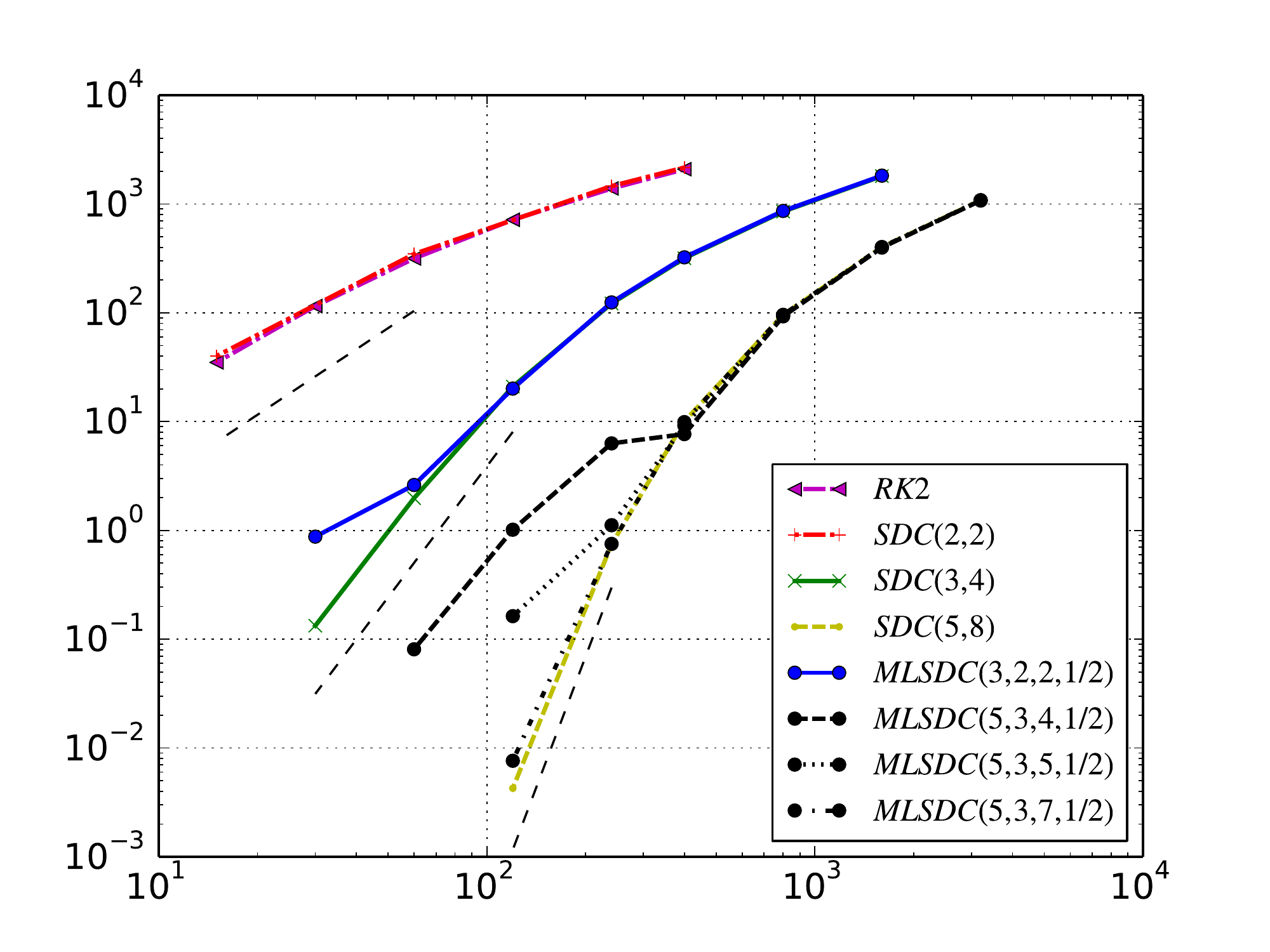}};
\node (ib_1) at (3.8,-0.05) {$\Delta t$};
\node[rotate=90] (ib_1) at (0,2.8) {$L_{\infty}$-norm of the error};
\end{tikzpicture}
\label{fig:gaussian_bump_test_case_accuracy_geopotential} 
}
\subfigure[]{
\begin{tikzpicture}
\node[anchor=south west,inner sep=0] at (0,0){\includegraphics[scale=0.365]{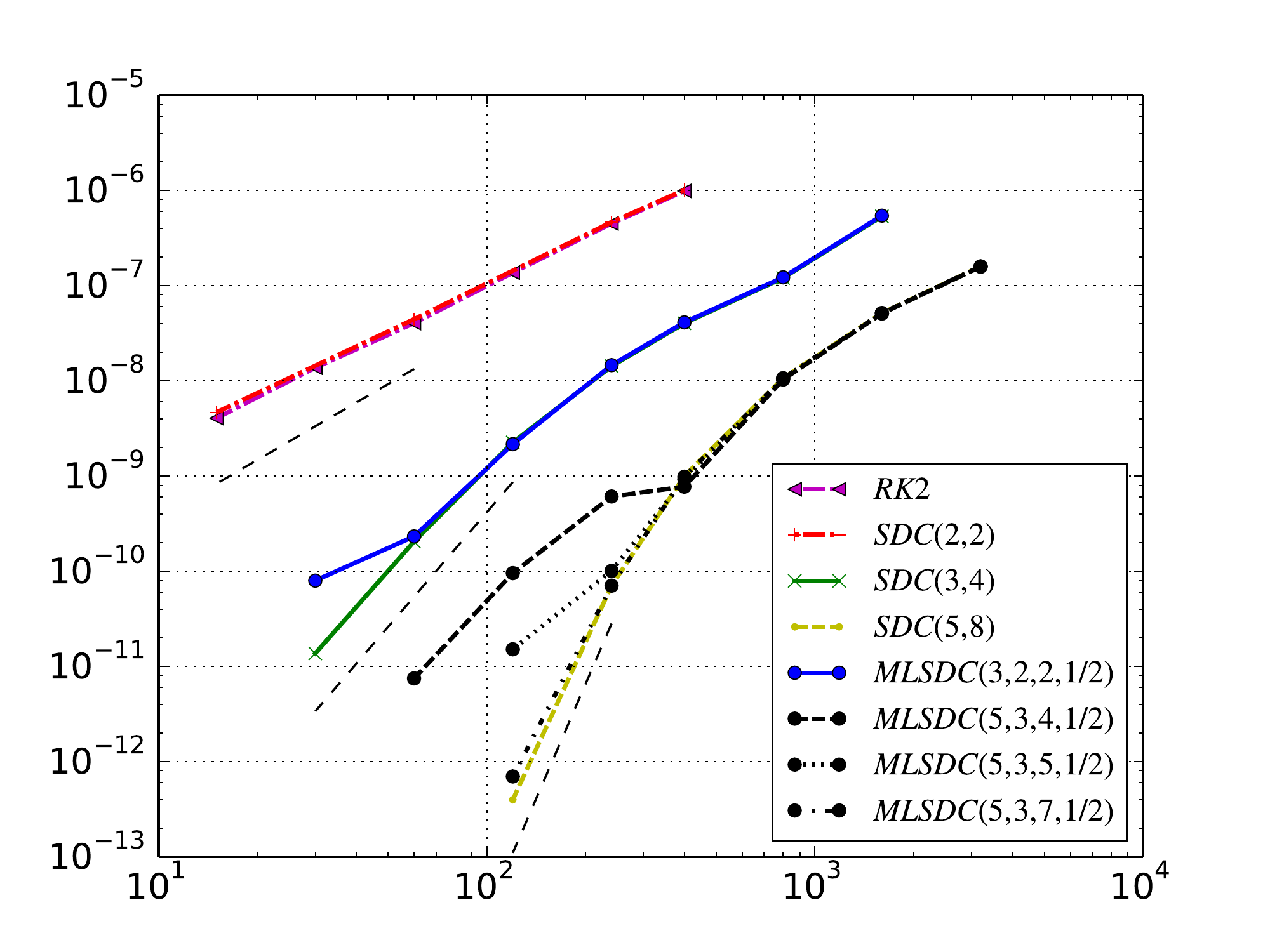}};
\node (ib_1) at (3.8,-0.05) {$\Delta t$};
\node[rotate=90] (ib_1) at (0.,2.8) {$L_{\infty}$-norm of the error};
\end{tikzpicture}
\label{fig:gaussian_bump_test_case_accuracy_vorticity} 
}
\vspace{-0.5cm}
\caption{\label{fig:gaussian_bump_test_case_accuracy} 
  Gaussian dome: $L_{\infty}$-norm of the error in the geopotential field
  in~\oldref{fig:gaussian_bump_test_case_accuracy_geopotential} and the vorticity field 
  in~\oldref{fig:gaussian_bump_test_case_accuracy_vorticity} with respect to the reference
  solution as a function of time step size. MLSDC(3,2,2,1/2) (respectively, MLSDC(5,3,4,1/2))
  achieves the same order of convergence as the SDC(3,4) (respectively, SDC(5,8))
  for larger time steps. For smaller time steps, more MLSDC-SH iterations are necessary
  to achieve the same order of convergence as the corresponding single-level schemes.
}
\end{figure}

To assess the accuracy of the MLSDC-SH scheme, we perform a refinement study in time for a fixed spatial resolution
($R_f = S_f = 256$) over one day. The reference solution is obtained with the single-level SDC(5,8) and a time step
size of $\Delta t_{\textit{ref}} \ = 60 \, \text{s}$. The  diffusion coefficient is set to
$\nu = 1.0 \times 10^5 \, \text{m}^2. \text{s}^{-1}$. For this test case, we focus again on the geopotential 
and vorticity fields, because the results for the divergence field are qualitatively similar. For completeness,
the study includes the results obtained with a second-order implicit-explicit Runge-Kutta scheme based on the
temporal splitting \ref{system_of_odes} to \ref{system_of_odes_explicit_part}. The $L_{\infty}$-norm of the error
with respect to the reference solution, as a function of time step size, is shown in Fig.~\oldref{fig:gaussian_bump_test_case_accuracy}. 

We see that MLSDC(3,2,2,1/2) achieves fourth-order convergence and an error of the same magnitude as
that obtained with the single-level SDC(3,4) for time step sizes such that $\Delta t \geq 60 \, \text{s}$. 
For a smaller time step size of $30 \, \text{s}$, the error induced by the truncation of high-frequency
modes during coarsening reduces the observed convergence of MLSDC(3,2,2,1/2) to second order. As in the
previous numerical example, the reduction in the observed convergence rate occurs when the $L_{\infty}$-norm of
the error reaches the magnitude of the spectral coefficients truncated during spatial coarsening -- about $10^{-1}$ 
for the geopotential according to Fig.~\oldref{fig:gaussian_test_case_spectra}.
Numerical results not included for brevity indicate that this reduction in accuracy caused by spatial 
coarsening persists even in the absence of temporal coarsening ($M_f = M_c$).
Still, in Fig.~\oldref{fig:gaussian_bump_test_case_accuracy}, the magnitude of the error obtained with MLSDC(3,2,2,1/2) in this range of small time steps
remains significantly smaller than that obtained with the single-level second-order SDC(2,2).
Fig.~\oldref{fig:gaussian_bump_test_case_accuracy} also shows that MLSDC(5,3,4,1/2) achieves the same 
convergence rate as the single-level SDC(5,8) whenever $\Delta t \geq 400 \, \text{s}$. For smaller
time step sizes, the observed convergence of MLSDC(5,3,4,1/2) is reduced to fourth order. But,
Fig.~\oldref{fig:gaussian_bump_test_case_accuracy} demonstrates that doing more iterations with 
MLSDC(5,3,7,1/2) is sufficient to recover eighth-order convergence in the asymptotic range. 

To interpret these results, we distinguish two regimes in the temporal refinement study of
Fig.~\oldref{fig:gaussian_bump_test_case_accuracy}. For very large time steps, large errors
are caused by the fact that the fine and coarse corrections do not resolve the large temporal scales 
accurately and overstep the small scales present in the problem. Reducing the time step size in this
regime reduces the errors associated with the large temporal scales. These scales can be resolved on
both the coarse level and the fine level which explains why MLSDC-SH and SDC converge at the same rate.
The second regime starts for smaller time steps once these large scales have been resolved accurately.
The error is then dominated by small-scale features that can be resolved by the fine correction but
cannot be captured by the coarse correction due to its lower spatial resolution of the 
coarse problem. This undermines the observed order of convergence of MLSDC-SH as the time step size
becomes very small.

\begin{figure}[ht!]
\centering
\subfigure[]{
\begin{tikzpicture}
\node[anchor=south west,inner sep=0] at (0,0){\includegraphics[scale=0.365]{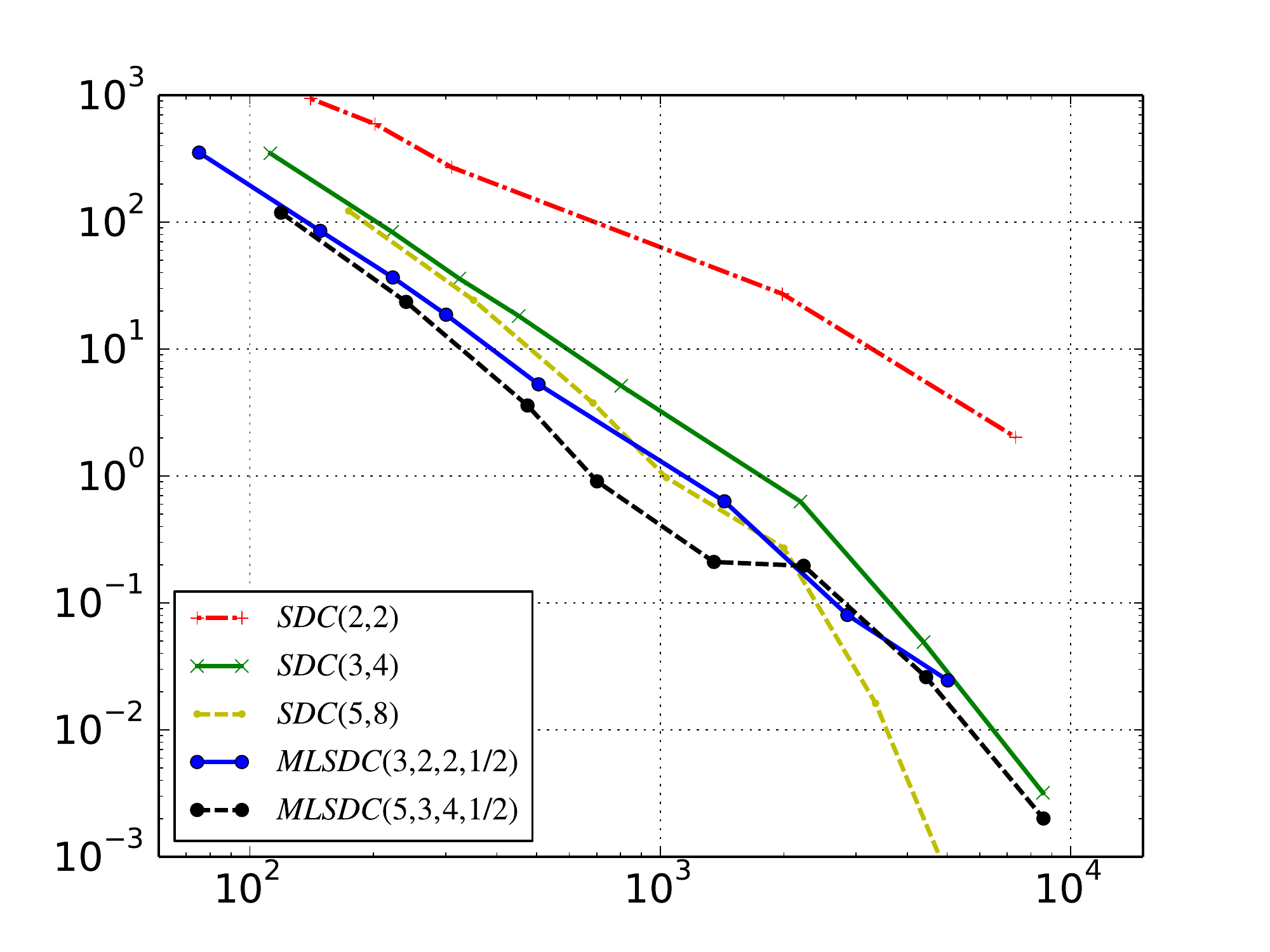}};
\node (ib_1) at (3.8,-0.05) {Wall-clock time};
\node[rotate=90] (ib_1) at (0.0,2.8) {$L_{2}$-norm of error};
\end{tikzpicture}
\label{fig:gaussian_bump_test_case_computational_cost_geopotential}
}
\subfigure[]{
\begin{tikzpicture}
\node[anchor=south west,inner sep=0] at (0,0){\includegraphics[scale=0.365]{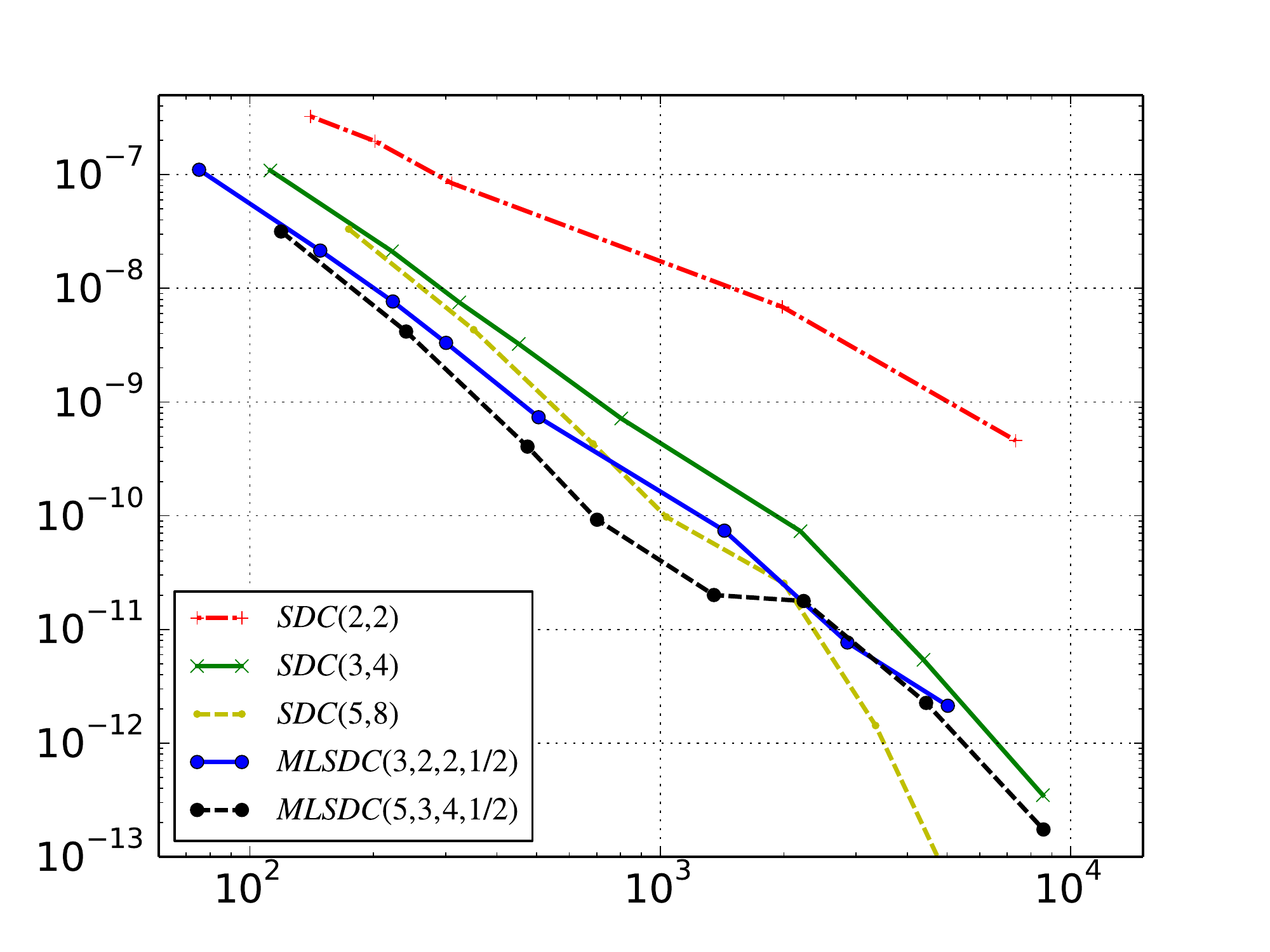}};
\node (ib_1) at (3.8,-0.05) {Wall-clock time};
\node[rotate=90] (ib_1) at (0.,2.8) {$L_{2}$-norm of error};
\end{tikzpicture}
\label{fig:gaussian_bump_test_case_computational_cost_vorticity}
}
\vspace{-0.5cm}
\caption{\label{fig:gaussian_bump_test_case_computational_cost}
  Gaussian dome: $L_{2}$-norm of the error in the geopotential field
  in~\oldref{fig:gaussian_bump_test_case_computational_cost_geopotential} 
  and the vorticity field in~\oldref{fig:gaussian_bump_test_case_computational_cost_vorticity}
  with respect to the reference solution as a function of the computational cost for the
  Gaussian dome test case. MLSDC(3,2,2,1/2) is more efficient than SDC(3,4) in the full 
  time step range. MLSDC(5,3,4,1/2) is more efficient than SDC(5,8) when the norm of the error is above $10^{-1}$ for the geopotential,
  and above $10^{-11}$ for the vorticity. Below these thresholds, MLSDC(5,3,4,1/2) only performs as efficiently
  as MLSDC(3,2,2,1/2) because of the reduction in its observed order of convergence upon refinement
  in time. 
}
\end{figure}

In Fig.~\oldref{fig:gaussian_bump_test_case_computational_cost}, we investigate the computational
cost the MLSDC-SH scheme by measuring the wall-clock time of the simulations.
MLSDC(3,2,2,1/2) is more efficient than the SDC(2,2) and SDC(3,4) schemes for the time step
sizes considered here despite the reduced observed order of convergence when $\Delta t \leq 60 \, \text{s}$.
MLSDC(5,3,4,1/2) is more efficient than SDC(5,8) whenever $\Delta t \geq 120 \, \text{s}$.
Below this time step size, its efficiency deteriorates slightly. Since the MLSDC-SH scheme
does not necessarily match the convergence rate of SDC in the simulations, we compute an observed speedup, $\mathcal{S}^{\textit{obs}}$, 
as the ratio of the computational cost of SDC over that of MLSDC-SH for a given error norm in
Fig.~\oldref{fig:gaussian_bump_test_case_computational_cost}.
For an error norm of $10^0$ in the geopotential field, MLSDC(3,2,2,1/2) achieves an observed speedup 
$\mathcal{S}^{\textit{obs}} \approx 1.58$ -- i.e., a reduction of 37 \% in wall-clock time -- compared to SDC(3,4). 
Considering that the cost of the FAS correction has been neglected in \ref{theoretical_speedup}, 
this is close to the theoretical speedup $\mathcal{S}^{\textit{theo}} \approx 1.66$ computed in 
Section~\oldref{section_temporal_discretization}. 
For the same magnitude of the error norm, MLSDC(5,3,4,1/2) achieves an observed speedup 
$\mathcal{S}^{\textit{obs}} \approx 1.50$ compared to SDC(5,8), which represents a reduction
of 33 \% in wall-clock time. This is again relatively close
to the theoretical speedup $\mathcal{S}^{\textit{theo}} \approx 1.66$. Although this is not
included in the figure for clarity, we point out that MLSDC(5,3,5,1/2) and MLSDC(5,3,7,1/2)
are more accurate, but also more expensive than MLSDC(5,3,4,1/2) in the range of time step
sizes considered here.

\subsection{\label{subsection_rossby_haurwitz_wave}Rossby-Haurwitz wave}

In this section, we apply MLSDC-SH to the Rossby-Haurwitz wave test case included in \cite{williamson1992standard} and also 
considered in \cite{jia2013spectral}. The initial analytically specified velocity field is non-divergent, and is computed with
wavenumber 4. The initial geopotential field is obtained by solving the balance equation. The resulting Haurwitz pattern moves
from east to west. We consider a fine resolution defined by $R_f = S_f = 256$ and we use a
diffusion coefficient $\nu = 1.0 \times 10^{5} \, \text{m}^2.\text{s}^{-1}$. In Fig.~\oldref{fig:rossby_haurwitz_wave_field}
(respectively, Fig.~\oldref{fig:rossby_haurwitz_wave_vorticity_spectrum}), we show the solution (respectively, the vorticity spectrum)
obtained with SDC(5,8) after one day.

\begin{figure}[ht!]
\centering
\subfigure[]{
\begin{tikzpicture}
\node[anchor=south west,inner sep=0] at (0,0){\includegraphics[scale=0.27]{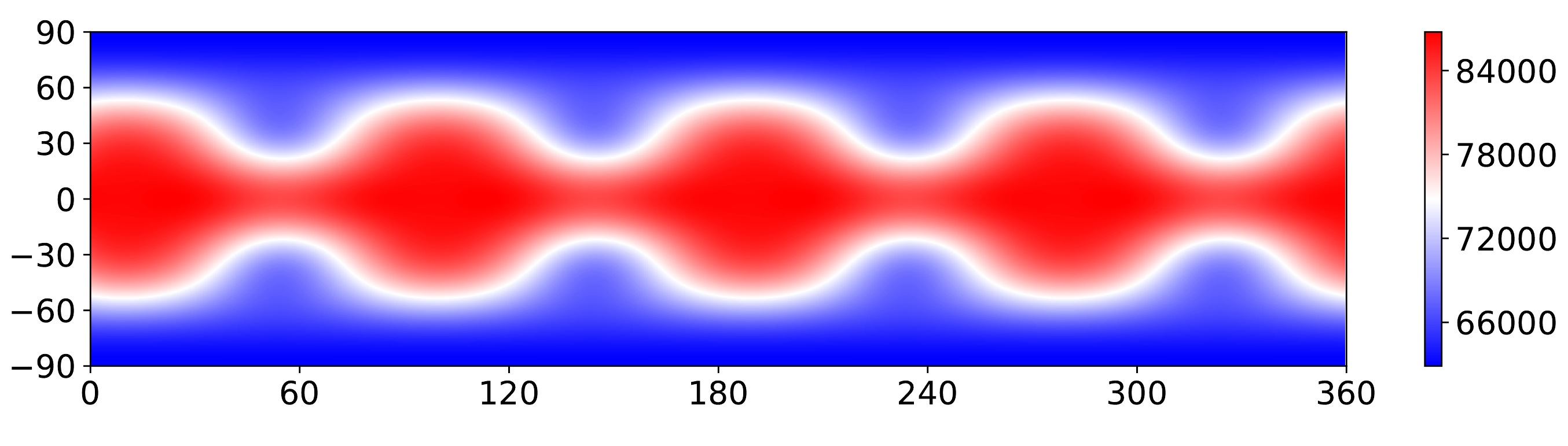}};
\node (ib_1) at (3.4,-0.05) {\scriptsize Longitude (degrees)};
\node[rotate=90] (ib_1) at (-0.2,1.2) {\scriptsize Latitude (degrees)};
\end{tikzpicture}
\label{fig:rossby_haurwitz_wave_geopotential_field}
}
\hspace{-0.5cm}
\subfigure[]{
\begin{tikzpicture}
\node[anchor=south west,inner sep=0] at (0,0){\includegraphics[scale=0.27]{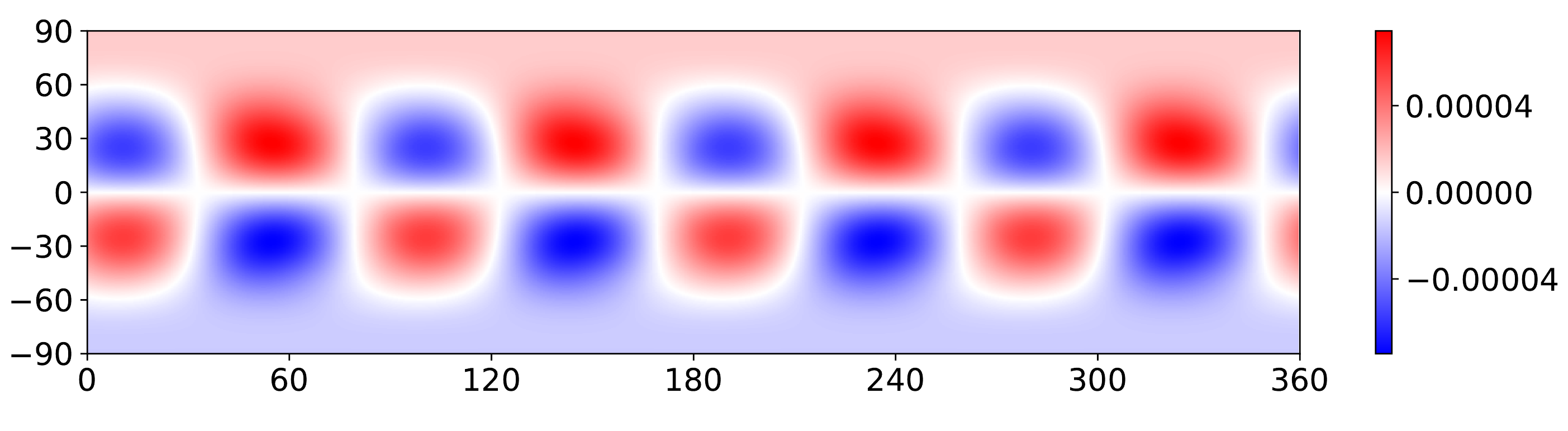}};
\node (ib_1) at (3.4,-0.05) {\scriptsize Longitude (degrees)};
\end{tikzpicture}
\label{fig:rossby_haurwitz_wave_vorticity_field}
} 
\vspace{-0.45cm}
\caption{\label{fig:rossby_haurwitz_wave_field} 
Rossby-Haurwitz wave: geopotential field in~\oldref{fig:rossby_haurwitz_wave_geopotential_field} and vorticity field 
in~\oldref{fig:rossby_haurwitz_wave_vorticity_field} 
with a resolution of $R_f = S_f = 256$ after one day. This solution is obtained with the single-level SDC(5,8). 
The diffusion coefficient is $\nu = 1.0 \times 10^{5} \, \text{m}^2.\text{s}^{-1}$. 
}
\end{figure}
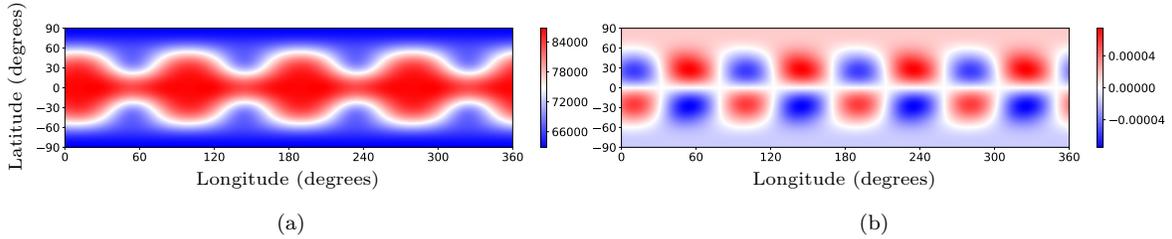
\begin{figure}[ht!]
\centering
\begin{tikzpicture}
\node[anchor=south west,inner sep=0] at (0,0){\includegraphics[scale=0.365]{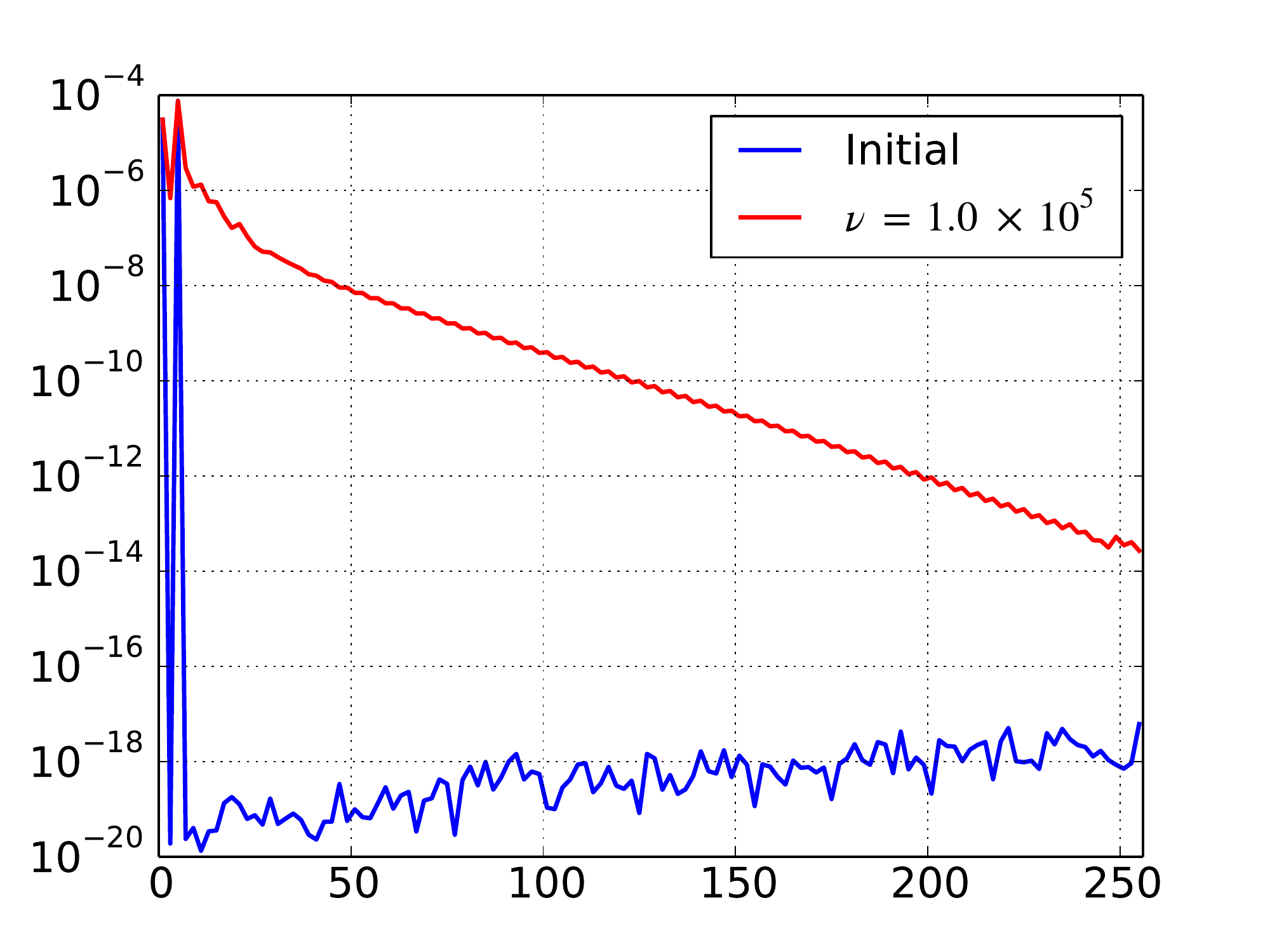}};
\node (ib_1) at (3.8,-0.05) {$n_0$};
\node[rotate=90] (ib_1) at (-0.05,2.7) {$|\zeta_{n_0}|$};

\node (ib_1) at (3.75,5.14) {};
\node (ib_2) at (3.75,0.45) {};
\path [draw=black!100, thick, dashed] (ib_1) -- (ib_2);

\node[rotate=-90] (ib_3) at (3.9,2.12) {\scriptsize $\alpha = 1/2$}; 

\end{tikzpicture}
\vspace{-0.4cm}
\caption{\label{fig:rossby_haurwitz_wave_vorticity_spectrum} 
Rossby-Haurwitz wave: max-spectrum of the vorticity field for the Rossby-Haurwitz wave test case at the beginning of the simulation
and after one day. The quantity on the $y$-axis is defined as $|\zeta_{n_0}| = \max_{r} |\zeta^r_{n_0}|$. To 
simplify the figure, we only show $|\zeta_{n_0}|$ for the even $n_0$. The odd $n_0$ correspond to negligible
$|\zeta_{n_0}|$ and are therefore omitted.
}
\end{figure}
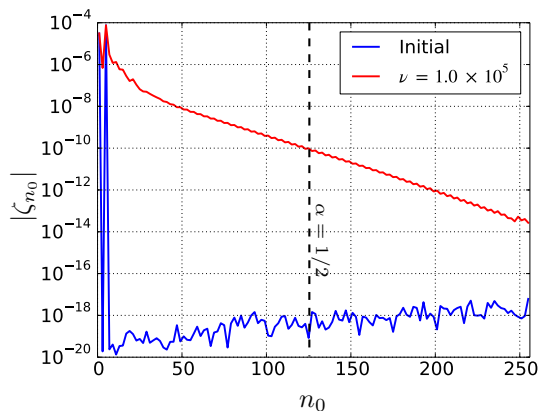

\begin{figure}[ht!]
\centering
\subfigure[]{
\begin{tikzpicture}
\node[anchor=south west,inner sep=0] at (0,0){\includegraphics[scale=0.365]{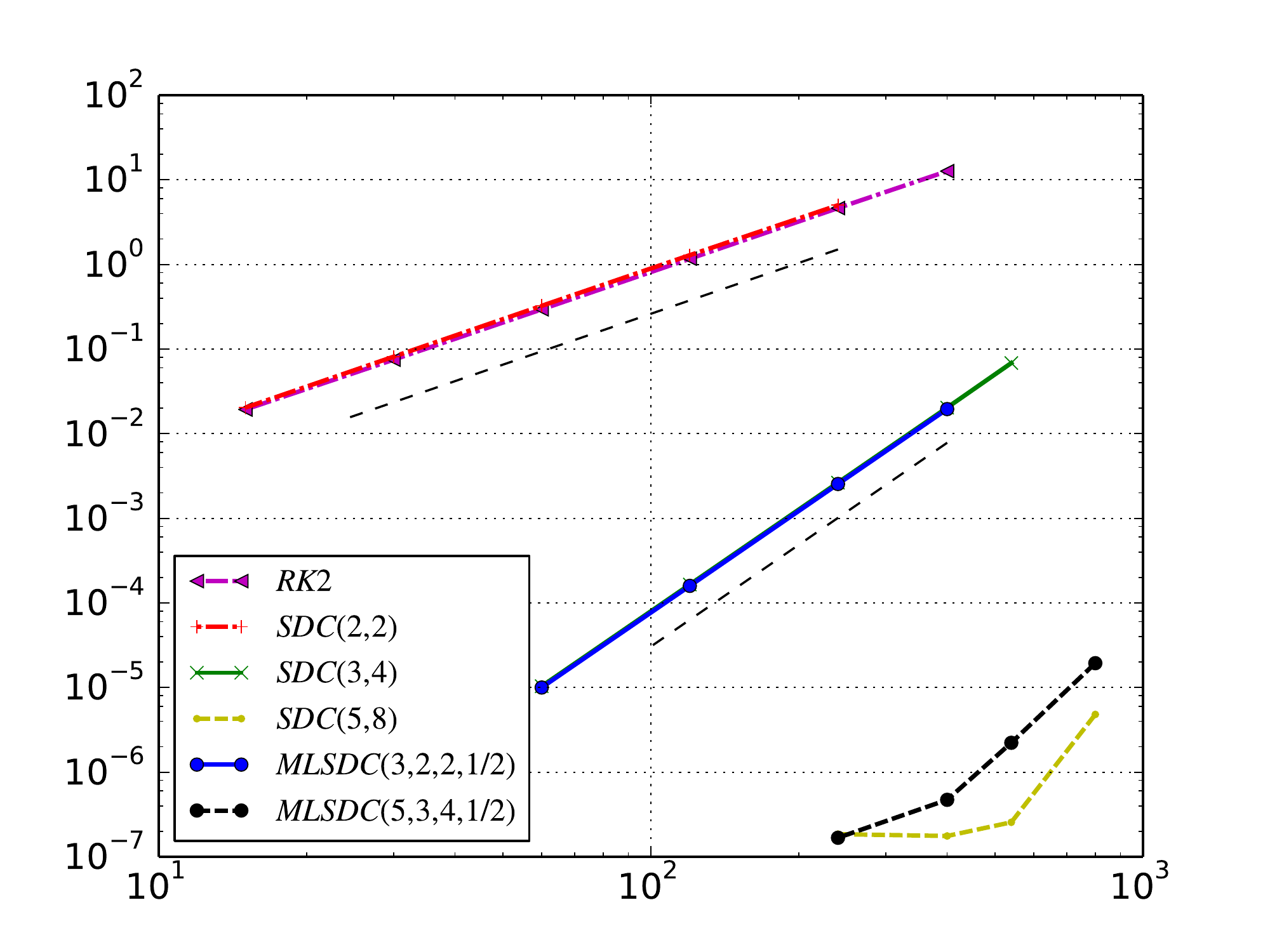}};
\node (ib_1) at (3.8,-0.05) {$\Delta t$};
\node[rotate=90] (ib_1) at (0,2.8) {$L_{\infty}$-norm of the error};
\end{tikzpicture}
\label{fig:rossby_haurwitz_wave_test_case_accuracy_geopotential} 
}
\subfigure[]{
\begin{tikzpicture}
\node[anchor=south west,inner sep=0] at (0,0){\includegraphics[scale=0.365]{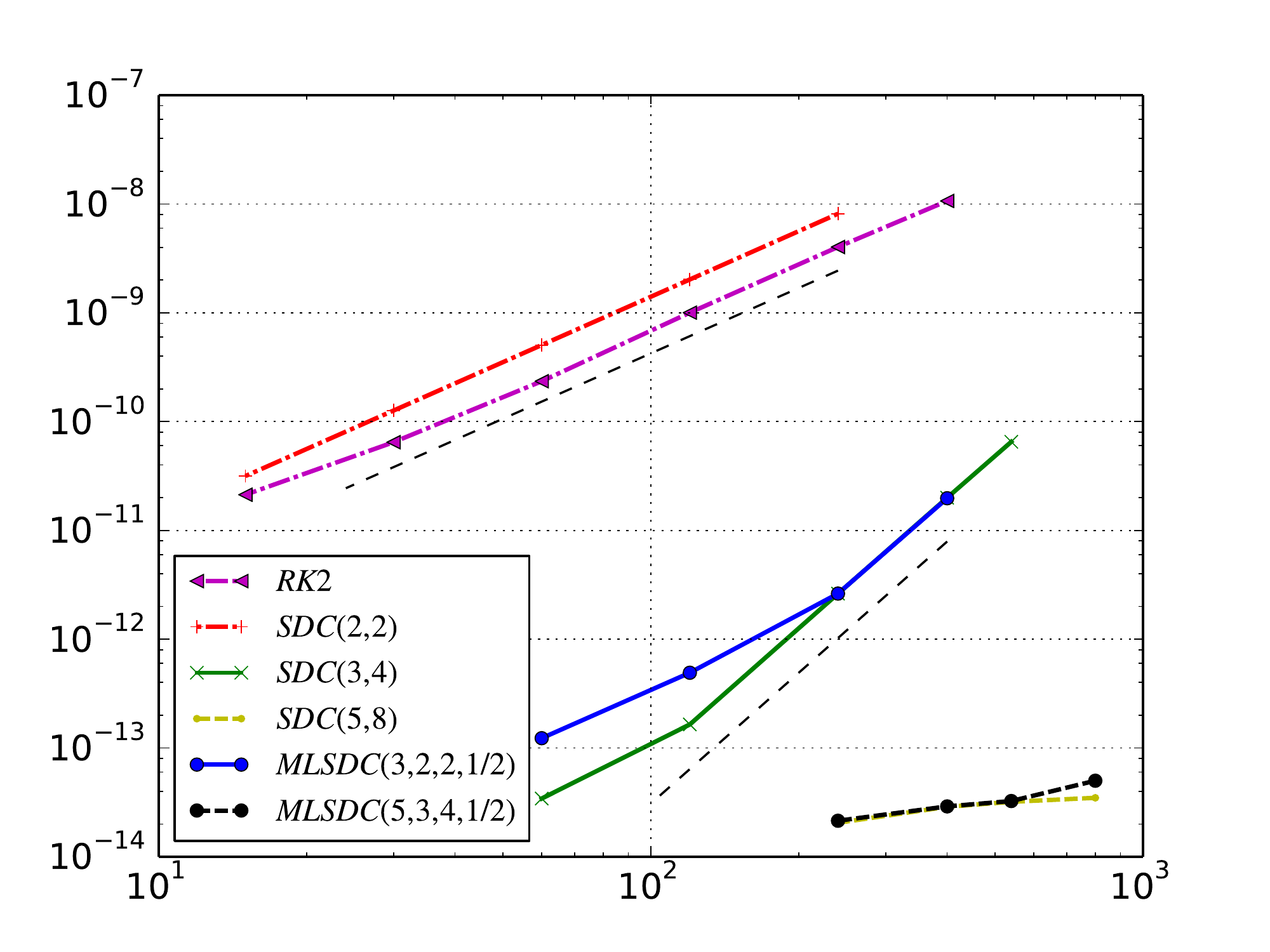}};
\node (ib_1) at (3.8,-0.05) {$\Delta t$};
\node[rotate=90] (ib_1) at (0.,2.8) {$L_{\infty}$-norm of the error};
\end{tikzpicture}
\label{fig:rossby_haurwitz_wave_test_case_accuracy_vorticity} 
}
\vspace{-0.45cm}
\caption{\label{fig:rossby_haurwitz_wave_test_case_accuracy} 
  Rossby-Haurwitz wave: $L_{\infty}$-norm of the error in the geopotential field
  in~\oldref{fig:rossby_haurwitz_wave_test_case_accuracy_geopotential}   and the vorticity field
  in~\oldref{fig:rossby_haurwitz_wave_test_case_accuracy_vorticity} with respect to the reference
  solution as a function of time step size. For this example, MLSDC(3,2,2,1/2) (respectively,
  MLSDC(5,3,4,1/2)) achieves the same observed order of convergence as the single-level SDC(3,4)
  (respectively, SDC(5,8)).
}
\end{figure}

As in the previous sections, we carry out a refinement study in time using a reference solution obtained with SDC(5,8)
over one day using a time step size $\Delta t_{\textit{ref}} = 120 \, \text{s}$. The results, shown in 
Fig.~\oldref{fig:rossby_haurwitz_wave_test_case_accuracy}, differ between the 
geopotential variable and the vorticity variable. Specifically, for the former, MLSDC(3,2,2,1/2) achieves fourth-order 
convergence upon refinement in time and the same error magnitude as SDC(3,4) for the full time step size range considered here.
MLSDC(5,3,4,1/2) is also more accurate than in the previous examples and reaches fifth-order convergence. But,
for the vorticity variable, MLSDC(3,2,2,1/2) exhibits a reduction in its convergence rate when $\Delta t \leq 120 \, \text{s}$,
which is slightly earlier than SDC(3,4). Here again, this reduction is caused by the truncation of high-frequency modes during
coarsening (see the vorticity spectrum in Fig.~\oldref{fig:rossby_haurwitz_wave_vorticity_spectrum}). MLSDC(5,3,4,1/2) is 
not in the asymptotic range and has already converged for the range of time step sizes considered here, which explains the 
flat line in Fig.~\oldref{fig:rossby_haurwitz_wave_test_case_accuracy_vorticity}.

\begin{figure}[ht!]
\centering
\subfigure[]{
\begin{tikzpicture}
\node[anchor=south west,inner sep=0] at (0,0){\includegraphics[scale=0.365]{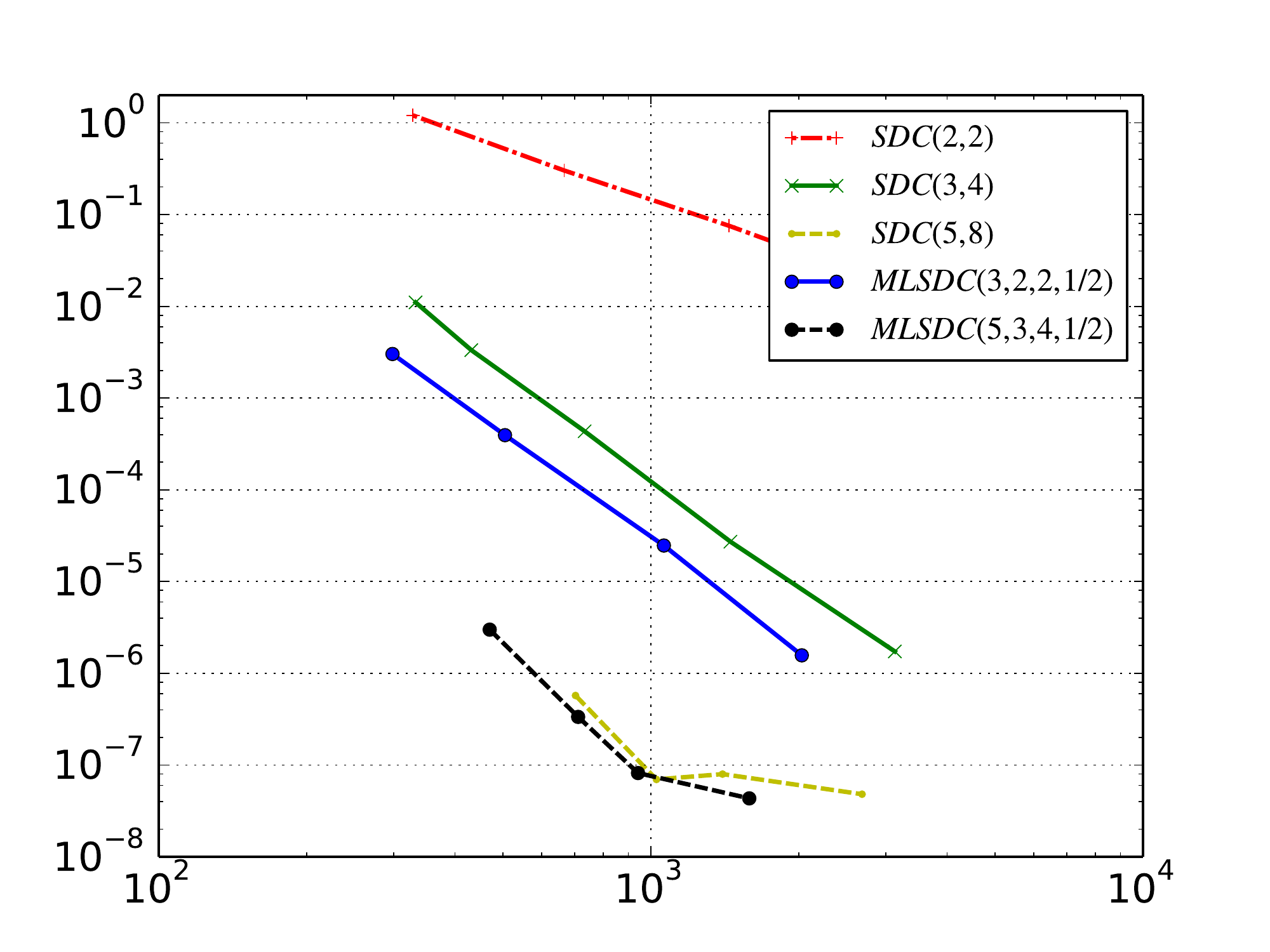}};
\node (ib_1) at (3.8,-0.05) {Wall-clock time};
\node[rotate=90] (ib_1) at (0.0,2.8) {$L_{2}$-norm of error};
\end{tikzpicture}
\label{fig:rossby_haurwitz_wave_test_case_computational_cost_geopotential}
}
\subfigure[]{
\begin{tikzpicture}
\node[anchor=south west,inner sep=0] at (0,0){\includegraphics[scale=0.365]{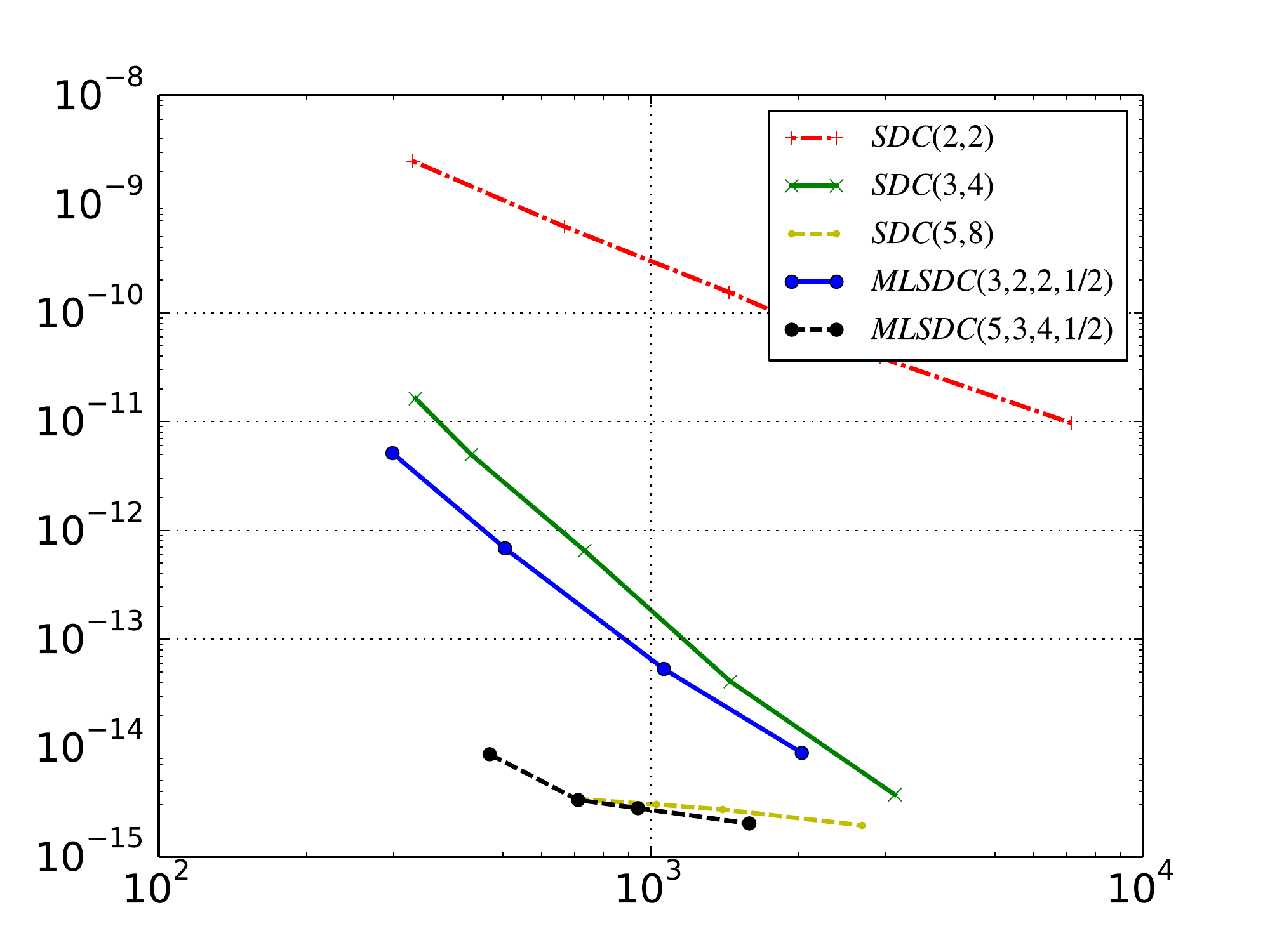}};
\node (ib_1) at (3.8,-0.05) {Wall-clock time};
\node[rotate=90] (ib_1) at (0.,2.8) {$L_{2}$-norm of error};
\end{tikzpicture}
\label{fig:rossby_haurwitz_wave_test_case_computational_cost_vorticity}
}
\vspace{-0.45cm}
\caption{\label{fig:rossby_haurwitz_wave_test_case_computational_cost}
  Rossby-Haurwitz wave: $L_{2}$-norm of the error in the geopotential field
  in~\oldref{fig:rossby_haurwitz_wave_test_case_computational_cost_geopotential} and the vorticity field
  in~\oldref{fig:rossby_haurwitz_wave_test_case_computational_cost_vorticity} with respect to the reference 
  solution as a function of the computational cost. The most efficient scheme is MLSDC(5,3,4,1/2). We also
  note that MLSDC(3,2,2,1/2) is more efficient than the single-level SDC(3,4) for the range of time step
  sizes considered here.
}
\end{figure}

In terms of wall-clock time, the most efficient scheme for both variables is MLSDC(5,3,4,1/2), as shown in 
Fig.~\oldref{fig:rossby_haurwitz_wave_test_case_computational_cost}. This multi-level scheme achieves a very low error norm for both 
variables while performing a large portion of the computations on the coarse level. MLSDC(3,2,2,1/2) is less efficient than 
MLSDC(5,3,4,1/2) but still more efficient than SDC(3,4) on the full range of time steps.
In particular, for an error norm of $10^{-4}$ in the geopotential field, MLSDC(3,2,2,1/2) achieves an observed speedup
$\mathcal{S}^{\textit{obs}} \approx 1.50$ compared to SDC(3,4), that is, a reduction in wall-clock time of 35 \%. This is 
in good agreement with the theoretical speedup $\mathcal{S}^{\textit{theo}} \approx 1.66$ computed
in Section~\oldref{subsection_computational_cost_of_sdc_and_mlsdc}. The MLSDC(3,2,2,1/2) performance
deteriorates for the vorticity variable for smaller time steps and the speedup compared to SDC(3,4) decreases because of the 
reduction in its observed order of convergence. For an error norm of $10^{-12}$ in the vorticity field, the observed speedup is
$S^{\textit{obs}} \approx 1.50$, but it is reduced to $S^{\textit{obs}} \approx 1.13$ for an error norm of $10^{-14}$. Next, we 
conclude the analysis of MLSDC-SH with a challenging unsteady test case representative of atmospheric flows.






\newpage

\subsection{\label{subsection_galewsky}Nonlinear evolution of an unstable barotropic wave}

In this section, we consider the barotropic instability test case proposed in \cite{galewsky2004initial}.
This is done by introducing a localized bump in the height field to perturb the balanced state described in
Section~\oldref{subsection_steady_zonal_jet}. The perturbation first triggers the development of gravity waves and 
then leads to the formation of complex vortical dynamics. These processes operate on multiple time scales
and are representative of the horizontal features of atmospheric flows. We run the simulations using two configurations, $\mathcal{B}$ and $\mathcal{C}$, 
based on a diffusion coefficient $\nu_{\mathcal{B}} = 1.0 \times 10^5 \, \text{m}^2.\text{s}^{-1}$ -- as in \cite{galewsky2004initial} --
and $\nu_{\mathcal{C}} = 2.0 \times 10^5 \, \text{m}^2.\text{s}^{-1}$, respectively. The reference solutions for the refinement 
studies detailed below are obtained with SDC(5,8) with a time step size $\Delta t_{\textit{ref}} = 60 \, \text{s}$. 
In \cite{jia2013spectral}, the largest time step size used in the single-level SDC scheme 
combined with the Spectral Element Method (SEM) based on 24 elements along each cube edge and a polynomial basis of degree seven is $1200 \, \text{s}$,
which is the same order of magnitude as the largest time step size used here.
The vorticity fields after 122 hours and 144 hours for configuration $\mathcal{B}$ are shown in 
Fig.~\oldref{fig:galewsky_test_case_with_localized_bump_vorticity_field}. 
Fig.~\oldref{fig:galewsky_test_case_with_localized_bump_vorticity_spectrum} presents the corresponding spectrum of the vorticity field 
at the same times. 

\begin{figure}[ht!]
\centering
\subfigure[]{
\begin{tikzpicture}
\node[anchor=south west,inner sep=0] at (0,0){\includegraphics[scale=0.27]{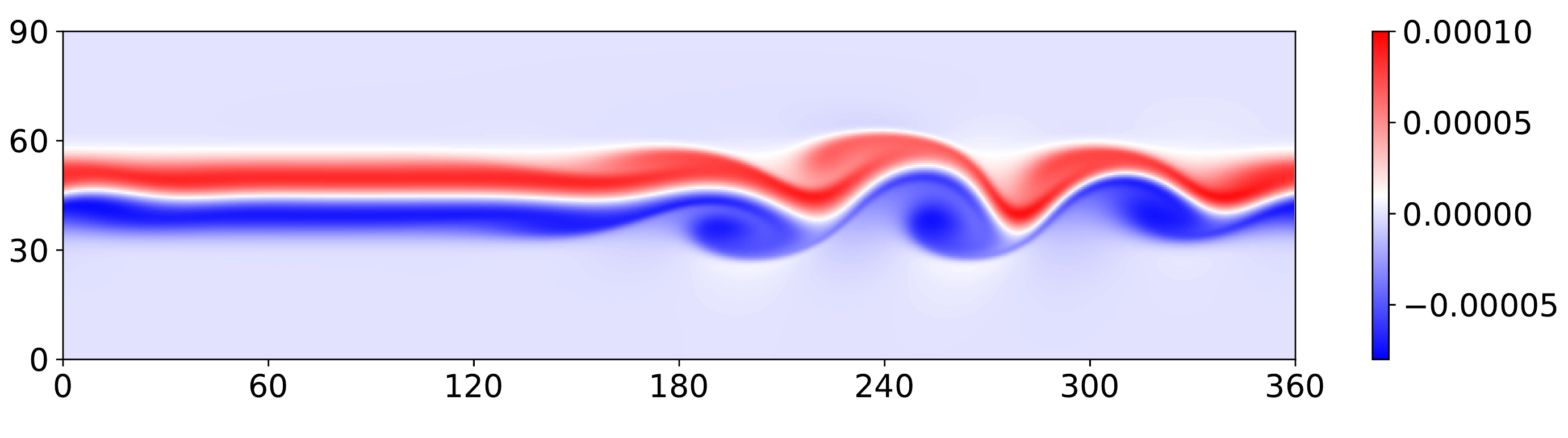}};
\node (ib_1) at (3.4,-0.05) {\scriptsize Longitude (degrees)};
\node[rotate=90] (ib_1) at (-0.2,1.2) {\scriptsize Latitude (degrees)};
\end{tikzpicture}
\label{fig:galewsky_test_case_with_localized_bump_vorticity_field_120}
}
\hspace{-0.5cm}
\subfigure[]{
\begin{tikzpicture}
\node[anchor=south west,inner sep=0] at (0,0){\includegraphics[scale=0.27]{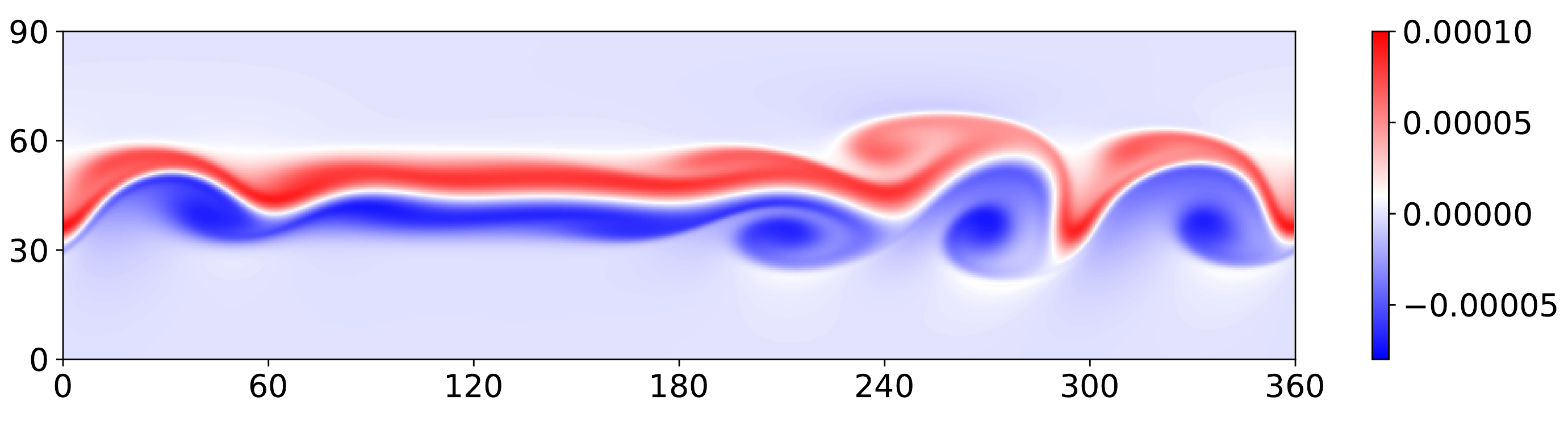}};
\node (ib_1) at (3.4,-0.05) {\scriptsize Longitude (degrees)};
\end{tikzpicture}
\label{fig:galewsky_test_case_with_localized_bump_vorticity_field_144}
} 
\vspace{-0.45cm}
\caption{\label{fig:galewsky_test_case_with_localized_bump_vorticity_field} 
Unstable barotropic wave: vorticity field with a resolution of $R_f = S_f = 256$ after 120 hours 
in~\oldref{fig:galewsky_test_case_with_localized_bump_vorticity_field_120} and 144 hours 
in~\oldref{fig:galewsky_test_case_with_localized_bump_vorticity_field_144}. This solution is obtained with the 
single-level SDC(5,8). The diffusion coefficient is 
$\nu_{\mathcal{B}} = 1.0 \times 10^{5} \, \text{m}^2.\text{s}^{-1}$. 
}
\end{figure}
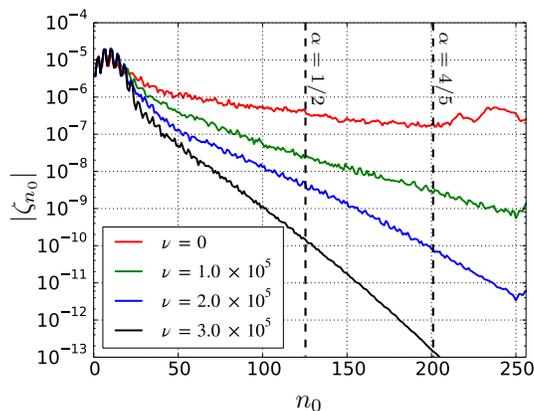
\begin{figure}[ht!]
\centering
\begin{tikzpicture}
\node[anchor=south west,inner sep=0] at (0,0){\includegraphics[scale=0.365]{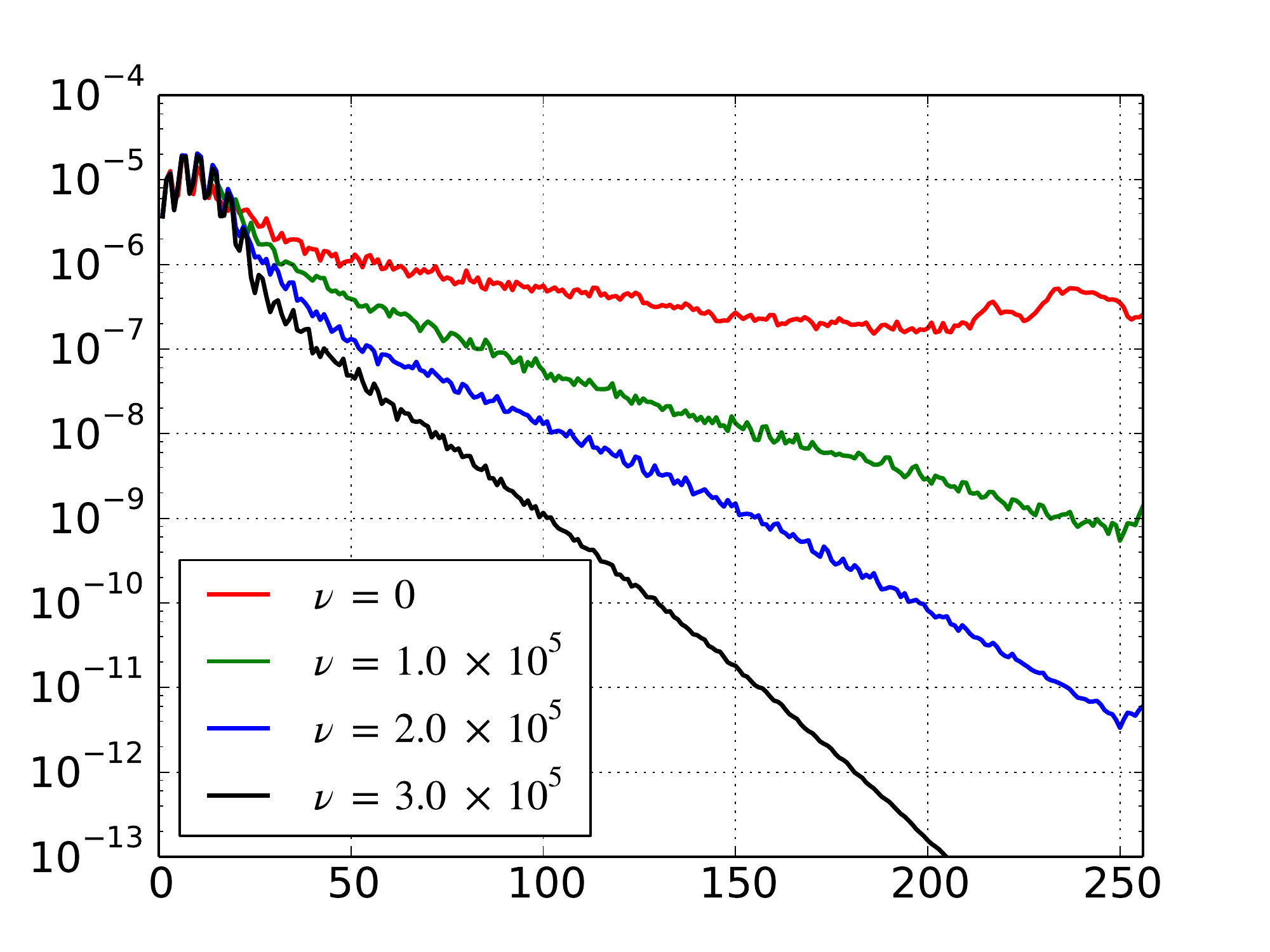}};
\node (ib_1) at (3.8,-0.05) {$n_0$};
\node[rotate=90] (ib_1) at (0,2.7) {$|\zeta_{n_0}|$};

\node (ib_1) at (3.75,5.14) {};
\node (ib_2) at (3.75,0.45) {};
\path [draw=black!100, thick, dashed] (ib_1) -- (ib_2);

\node (ib_1) at (5.4475,5.14) {};
\node (ib_2) at (5.4475,0.45) {};
\path [draw=black!100, thick, dashed] (ib_1) -- (ib_2);

\node[rotate=-90] (ib_3) at (3.9,4.4) {\scriptsize $\alpha = 1/2$}; 

\node[rotate=-90] (ib_3) at (5.5975,4.4) {\scriptsize $\alpha = 4/5$};

\end{tikzpicture}
\vspace{-0.4cm}
\caption{\label{fig:galewsky_test_case_with_localized_bump_vorticity_spectrum} 
Unstable barotropic wave: max-spectrum of the vorticity field after 144 hours 
for various values of the diffusion coefficient. The 
quantity on the $y$-axis is defined as $|\zeta_{n_0}| = \max_{r} |\zeta^r_{n_0}|$.
We observe a strong damping of the high-frequency modes when the diffusion coefficient
is large. This has a significant impact on the observed order of convergence of MLSDC-SH as 
shown in Fig.~\oldref{fig:galewsky_test_case_with_localized_bump_vorticity_accuracy_l_infty_norm}.
}
\end{figure}
\begin{figure}[ht!]
\centering
\begin{tikzpicture}
\node[anchor=south west,inner sep=0] at (0,0){\includegraphics[scale=0.365]{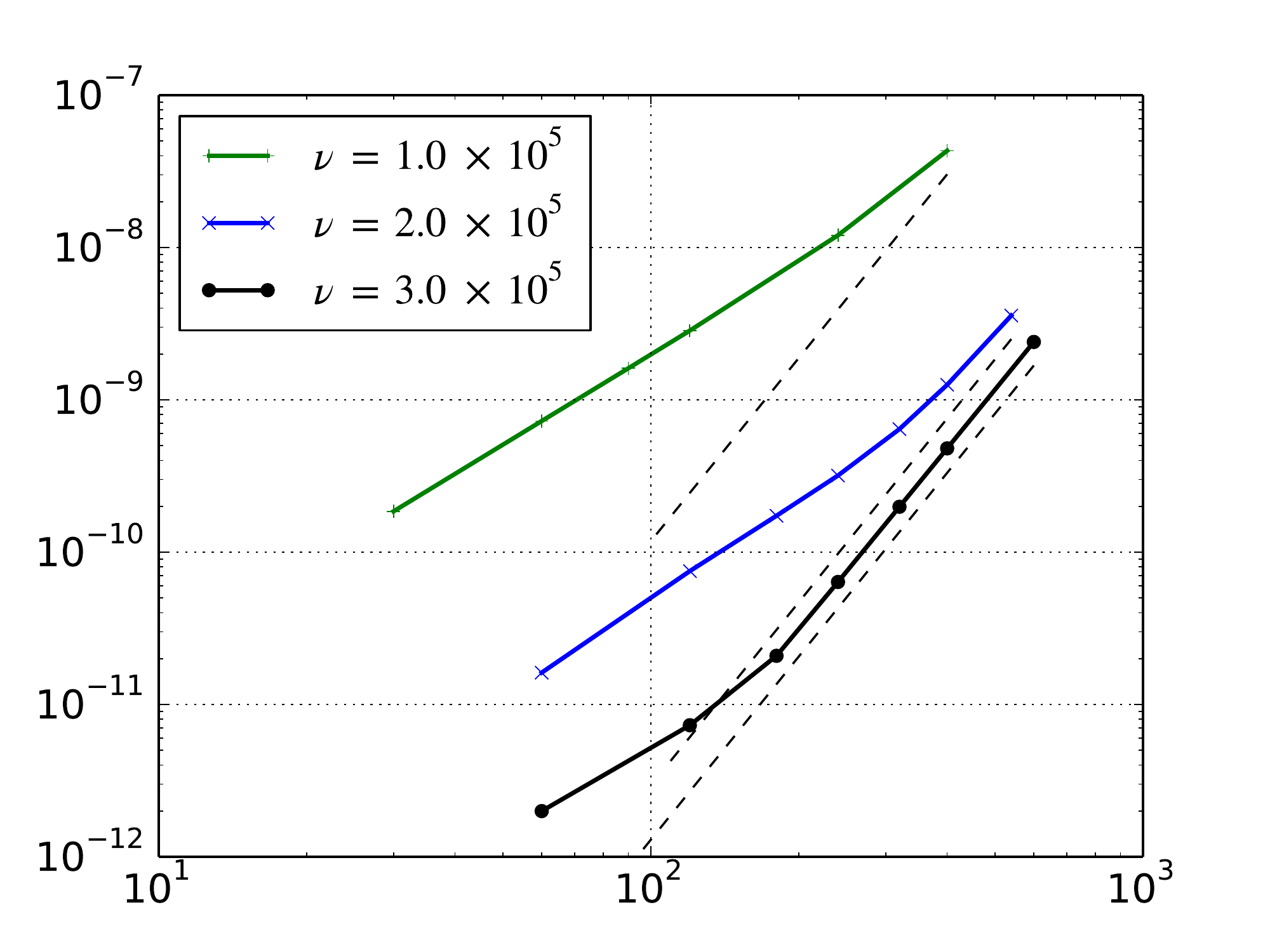}};
\node (ib_1) at (3.8,-0.05) {$\Delta t$};
\node[rotate=90] (ib_1) at (0,2.8) {$L_{\infty}$-norm of the error};
\end{tikzpicture}
\vspace{-0.4cm}
\caption{\label{fig:galewsky_test_case_with_localized_bump_vorticity_accuracy_l_infty_norm} 
Unstable barotropic wave: $L_{\infty}$-norm of the error in the vorticity field with respect to the reference 
solution as a function of the time step size. These results are obtained with MLSDC(3,2,2,1/2). We note again 
that a reduction in the observed order of convergence of MLSDC-SH occurs when the norm of the error is smaller 
than the magnitude of the spectral coefficients that are truncated during spatial coarsening 
(see Fig.~\oldref{fig:galewsky_test_case_with_localized_bump_vorticity_spectrum}).
}
\end{figure}

As in Section~\oldref{subsection_steady_zonal_jet}, we first highlight the connection between the spectrum of the vorticity 
field and the observed order of convergence of the MLSDC-SH scheme upon refinement in time. This is done with MLSDC(3,2,2,1/2) -- 
that is, MLSDC-SH with three nodes on the fine level, two nodes on the coarse level, two iterations, and $R_c = S_c = 128$ -- 
in the refinement study in time shown in Fig.~\oldref{fig:galewsky_test_case_with_localized_bump_vorticity_accuracy_l_infty_norm}.
When $\nu_{\mathcal{B}} = 1.0 \times 10^5 \, \text{m}^2.\text{s}^{-1}$, the magnitude of the truncated terms in the vorticity spectrum is 
of the order of $10^{-8}$ (see Fig.~\oldref{fig:galewsky_test_case_with_localized_bump_vorticity_spectrum}). Since this is also the 
order of the $L_{\infty}$-norm of the error for the largest stable time step ($\Delta t = 400 \, \text{s}$), MLSDC(3,2,2,1/2)  
achieves only second-order convergence for the range of time step sizes considered here. With $\nu_{\mathcal{C}} = 2.0 \times 10^5 \, \text{m}^2.\text{s}^{-1}$, 
the magnitude of the truncated terms is of the order of $10^{-9}$, and MLSDC(3,2,2,1/2) reaches fourth-order convergence until this threshold 
is reached for $\Delta t = 320 \, \text{s}$. For completeness, we have run the same test with $\nu = 3.0 \times 10^5 \, \text{m}^2.\text{s}^{-1}$, 
in which case this threshold is lower, which allows MLSDC(3,2,2,1/2) to exhibit fourth-order convergence for $\Delta t \geq 120 \, \text{s}$.
In the following paragraphs, we show that this reduction in the observed order of convergence can be overcome by doing additional MLSDC-SH iterations.
For instance, we demonstrate that MLSDC(3,2,3,1/2)  recovers fourth-order convergence in 
configurations $\mathcal{B}$ and $\mathcal{C}$.

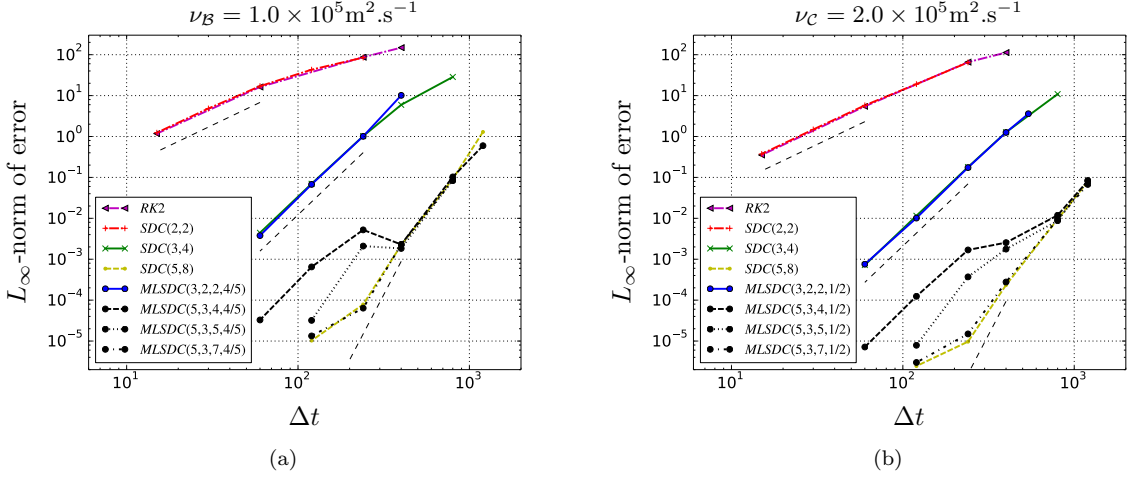
\begin{figure}[ht!]
\centering
\subfigure[]{
\begin{tikzpicture}
\node[anchor=south west,inner sep=0] at (0,0){\includegraphics[scale=0.365]{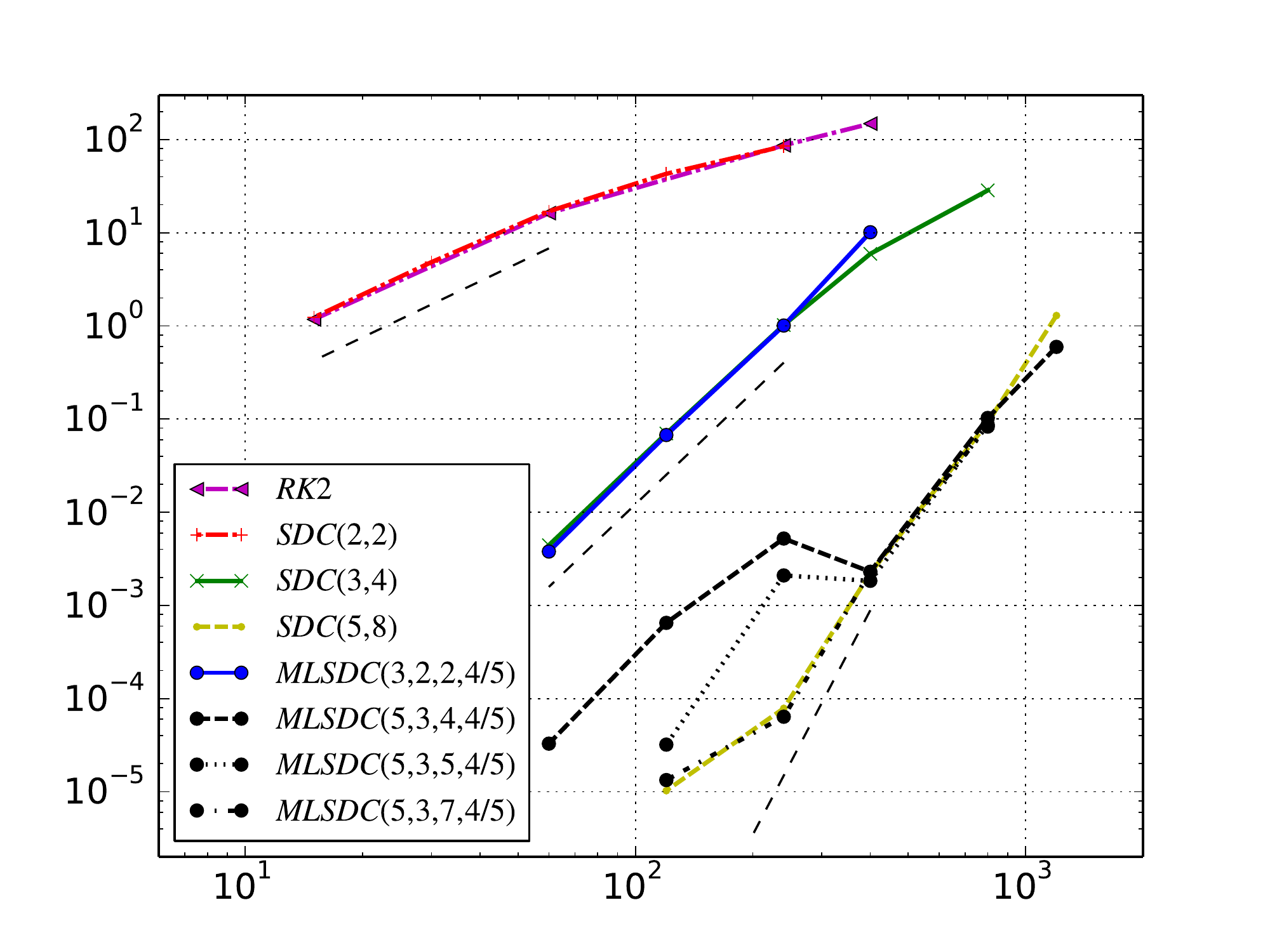}};
\node (ib_1) at (3.8,-0.05) {$\Delta t$};
\node[rotate=90] (ib_1) at (0.0,2.8) {$L_{\infty}$-norm of error};
\node (ib_1) at (3.8,5.3) {\small $\nu_{\mathcal{B}} = 1.0 \times 10^{5} \text{m}^2.\text{s}^{-1}$};
\end{tikzpicture}
\label{fig:galewsky_test_case_with_localized_bump_accuracy_geopotential_A}
}
\subfigure[]{
\begin{tikzpicture}
\node[anchor=south west,inner sep=0] at (0,0){\includegraphics[scale=0.365]{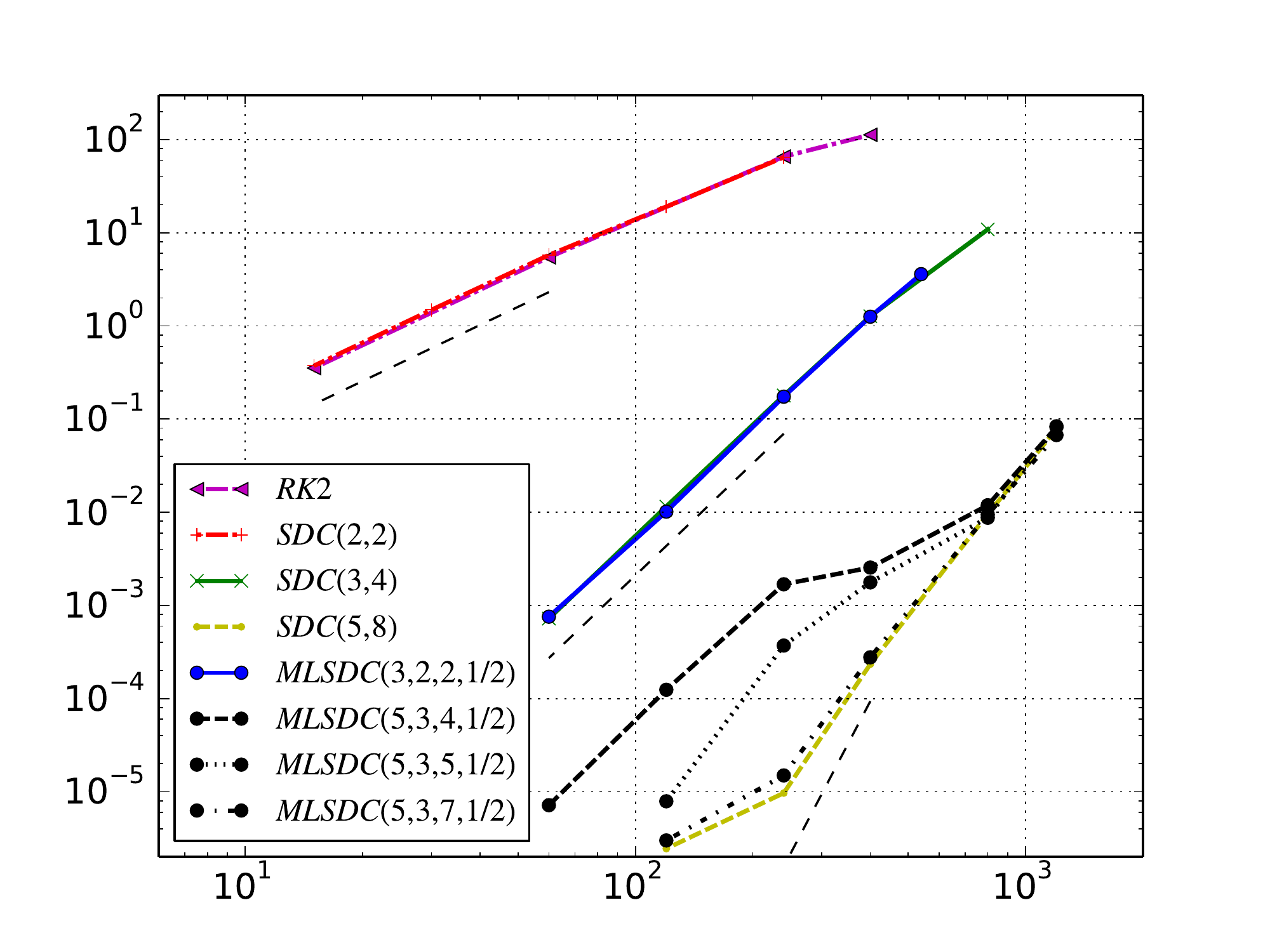}};
\node (ib_1) at (3.8,-0.05) {$\Delta t$};
\node[rotate=90] (ib_1) at (0.0,2.8) {$L_{\infty}$-norm of error};
\node (ib_1) at (3.8,5.3) {\small $\nu_{\mathcal{C}} = 2.0 \times 10^{5} \text{m}^2.\text{s}^{-1}$};
\end{tikzpicture}
\label{fig:galewsky_test_case_with_localized_bump_accuracy_geopotential_B}
}
\vspace{-0.45cm}
\caption{\label{fig:galewsky_test_case_with_localized_bump_accuracy_geopotential} 
Unstable barotropic wave: $L_{\infty}$-norm of the error in the geopotential field with respect to the reference 
solution as a function of the time step size. MLSDC(3,2,2,$\alpha$) achieves the same order of convergence
as SDC(3,4) in both configurations. MLSDC(5,3,4,$\alpha$) is also more accurate than SDC(3,4), but achieves the same 
accuracy as SDC(5,8) for larger time steps only. Doing additional MLSDC-SH iterations -- in this case, seven 
iterations with MLSDC(5,3,7,$\alpha$) -- is needed to match the accuracy of SDC(5,8) in the full range of time step sizes.
}
\end{figure}

We now use this knowledge of the spectrum of the vorticity field to motivate our choice of the spatial 
coarsening ratio in each configuration. The goal is to make the coarse sweeps as inexpensive as possible 
without undermining the observed order of convergence of the MLSDC-SH scheme upon refinement in time. In configuration 
$\mathcal{B}$, we choose a relatively modest spatial coarsening ratio $\alpha_{\mathcal{B}} = 4/5$ to account 
for the presence of large spectral coefficients associated with the high-frequency modes. At the coarse level, this 
choice yields $R^{\mathcal{B}}_c = S^{\mathcal{B}}_c = 204$. In configuration $\mathcal{C}$, we can choose a more 
aggressive coarsening strategy with $\alpha_{\mathcal{C}} = 1/2$, leading to $R^{\mathcal{C}}_c = S^{\mathcal{C}}_c = 128$ 
as in the previous test cases. These choices are such that the magnitude of the truncated spectral vorticity coefficients 
have the same order of magnitude, that is $|\zeta^{\mathcal{B}}_{204}| \approx |\zeta^{\mathcal{C}}_{128}| \approx 3 \times 10^{-9}$ in 
Fig.~\oldref{fig:galewsky_test_case_with_localized_bump_vorticity_spectrum}. The results of the refinement study in $L_{\infty}$-norm are 
shown in Fig.~\oldref{fig:galewsky_test_case_with_localized_bump_accuracy_geopotential} for the geopotential. The asymptotic rates observed 
for the divergence and the vorticity are qualitatively similar to those of the geopotential and are therefore omitted for brevity. 

In Fig.~\oldref{fig:galewsky_test_case_with_localized_bump_accuracy_geopotential}, 
MLSDC(3,2,2,$\alpha$) exhibits the same observed order of convergence and error magnitude as 
SDC(3,4) for both diffusion configurations. MLSDC(5,3,4,$\alpha$) also converges at a fourth-order 
rate in the asymptotic range, but achieves a significantly smaller error magnitude than MLSDC(3,2,2,$\alpha$). 
MLSDC(5,3,4,$\alpha$) is as accurate as SDC(5,8) for larger time step sizes. But, to achieve the same observed order of convergence 
as SDC(5,8) in the entire range of time step sizes, seven iterations -- with MLSDC(5,3,7,$\alpha$) -- are needed. 
Finally, we note that numerical examples not shown here for brevity confirm that the observed order 
of convergence of MLSDC-SH increases significantly for a larger coefficient $\alpha$ -- i.e., less aggressive spatial coarsening --,
but this also drastically increases the computational cost.

\begin{figure}[ht!]
\centering
\subfigure[]{
\begin{tikzpicture}
\node[anchor=south west,inner sep=0] at (0,0){\includegraphics[scale=0.365]{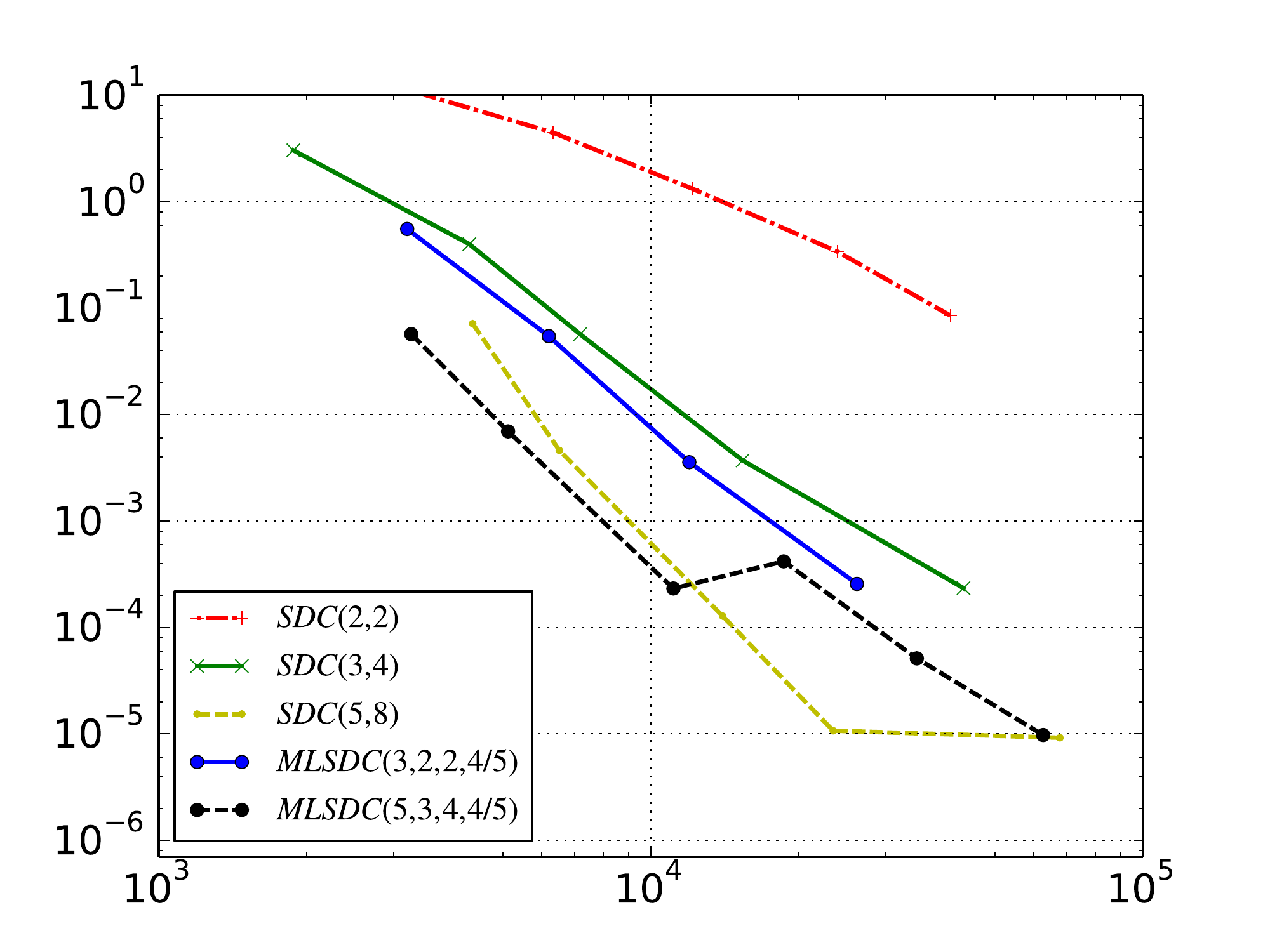}};
\node (ib_1) at (3.8,-0.05) {Wall-clock time};
\node[rotate=90] (ib_1) at (0.0,2.8) {$L_2$-norm of error};
\node (ib_1) at (3.8,5.3) {\small $\nu_{\mathcal{B}} = 1.0 \times 10^{5} \text{m}^2.\text{s}^{-1}$};
\end{tikzpicture}
\label{fig:galewsky_test_case_with_localized_bump_computational_cost_geopotential_A}
}
\subfigure[]{
\begin{tikzpicture}
\node[anchor=south west,inner sep=0] at (0,0){\includegraphics[scale=0.365]{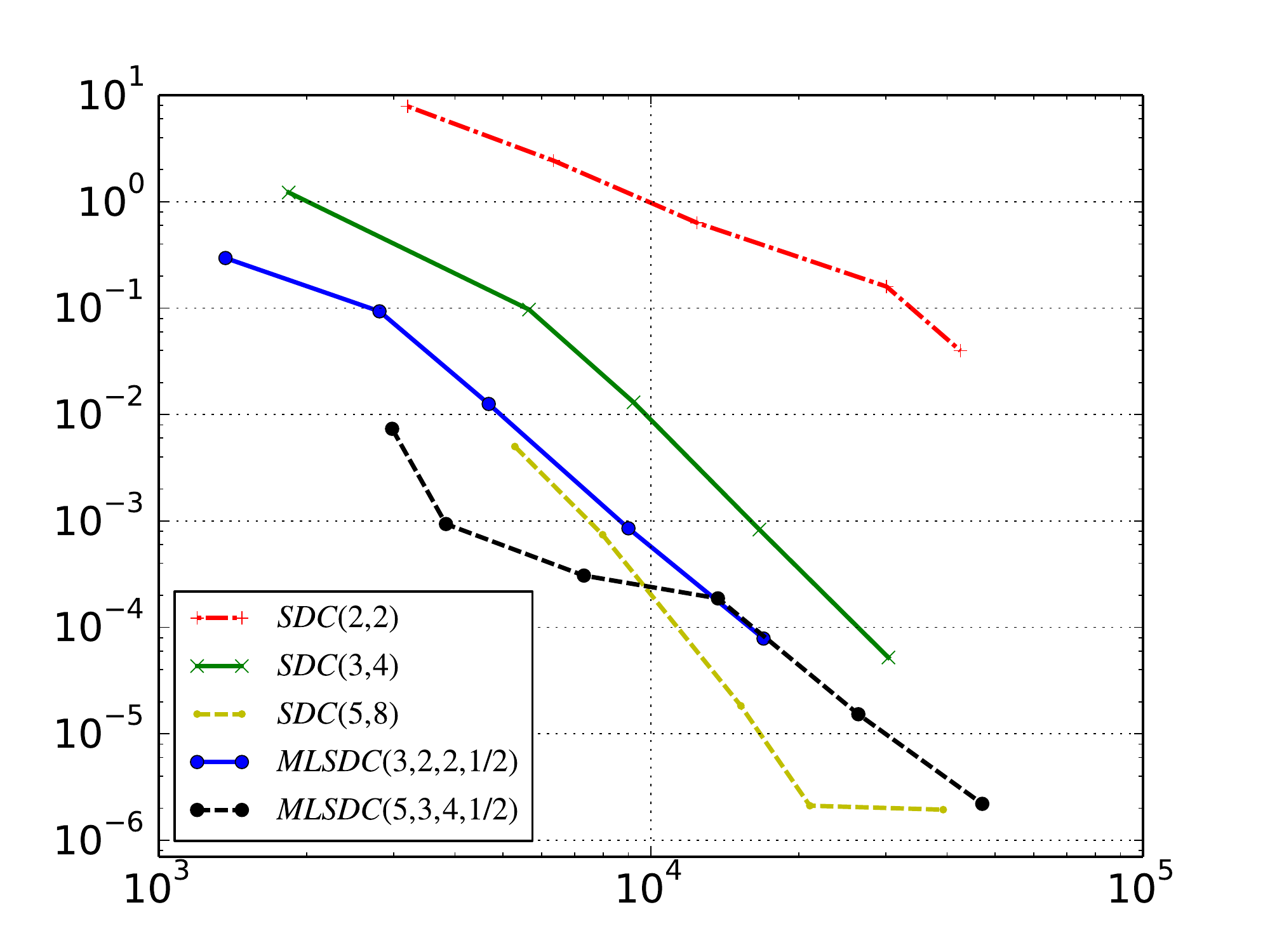}};
\node (ib_1) at (3.8,-0.05) {Wall-clock time};
\node[rotate=90] (ib_1) at (0.0,2.8) {$L_2$-norm of error};
\node (ib_1) at (3.8,5.3) {\small $\nu_{\mathcal{C}} = 2.0 \times 10^{5} \text{m}^2.\text{s}^{-1}$};
\end{tikzpicture}
\label{fig:galewsky_test_case_with_localized_bump_computational_cost_geopotential_B}
}
\vspace{-0.45cm}
\caption{\label{fig:galewsky_test_case_with_localized_bump_computational_cost_geopotential} 
Unstable barotropic wave: $L_2$-norm of the error in the geopotential field with respect to the reference 
solution as a function of the computational cost. MLSDC(3,2,2,$\alpha$) is more efficient than SDC(3,4). 
The cost reduction is larger when the spatial coarsening is more aggressive ($\alpha = 1/2$). 
MLSDC(5,3,4,$\alpha$) is more efficient that SDC(5,8) for larger error magnitudes. In the range of 
relatively larger errors, MLSDC(5,3,4,$\alpha$) is the most efficient scheme among those considered here. 
}
\end{figure}
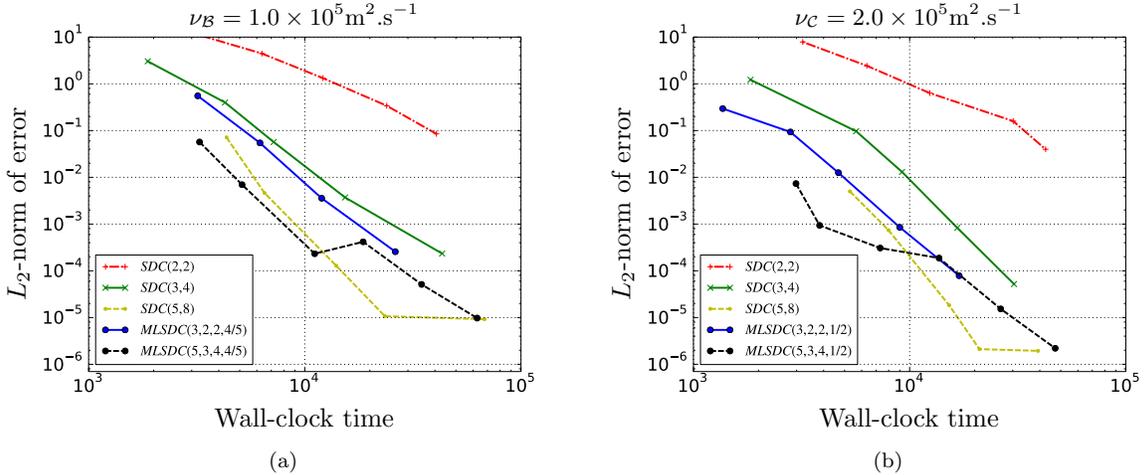

In Fig.~\oldref{fig:galewsky_test_case_with_localized_bump_computational_cost_geopotential}, we show the $L_2$-norm 
of the error as a function of the wall-clock time of the simulations for the geopotential. We see that MLSDC(3,2,2,$\alpha$) is significantly less expensive 
than SDC(3,4) in the full range of time step sizes. This cost reduction is larger with MLSDC(3,2,2,$\alpha_{\mathcal{C}}$) 
since configuration $\mathcal{C}$ allows for a more aggressive spatial coarsening strategy than configuration $\mathcal{B}$. 
In both configurations, the observed speedup of MLSDC(3,2,2,$\alpha$) compared to SDC(3,4) is close to the theoretical speedup. 
Specifically, in configuration $\mathcal{B}$, the observed speedup is $\mathcal{S}_{\mathcal{B}}^{\textit{obs}} \approx 1.28$ 
for an error norm of $3 \times 10^{-3}$ in the geopotential field whereas the theoretical speedup -- obtained with \ref{theoretical_speedup} 
evaluated with $\alpha_{\mathcal{B}} = 4/5$ -- is $\mathcal{S}_{\mathcal{B}}^{\textit{theo}} \approx 1.30$. In configuration $\mathcal{C}$, 
MLSDC(3,2,2,1/2) achieves $\mathcal{S}^{\textit{obs}}_{\mathcal{C}} \approx 1.56$ for an error norm of $8 \times 10^{-5}$ in the geopotential,
for a theoretical speedup $\mathcal{S}^{\textit{theo}}_{\mathcal{C}} \approx 1.66$. 
MLSDC(5,3,4,$\alpha$) is the most efficient scheme for relatively large error magnitudes and also achieves observed speedups
 close to the theoretical speedup. But, the performance of MLSDC(5,3,4,$\alpha$)  deteriorates for lower error magnitudes. 
We found that doing additional iterations, for instance with MLSDC(3,2,3,$\alpha$) or MLSDC(5,3,5,$\alpha$), does not improve 
the efficiency of MLSDC-SH.

\section{\label{section_conclusion}Conclusions and future work}


We have studied a high-order implicit-explicit iterative multi-level time integration scheme for the 
nonlinear shallow-water equations on the rotating sphere. Our algorithm relies on the Multi-Level Spectral 
Deferred Corrections (MLSDC) scheme of \cite{emmett2012toward,speck2015multi} combined with a spatial discretization 
performed with the global Spherical Harmonics (SH) transform. MLSDC-SH applies a sequence of updates distributed 
on a hierarchy of space-time levels obtained by coarsening the problem in space and in time. This approach 
makes it possible to shift a significant portion of the computational work to the coarse representation of 
the problem to reduce the time-to-solution while preserving accuracy.

We have discussed the requirements of consistent inter-level transfer operators which play a crucial role in MLSDC-SH. 
Our approach consists in exploiting the canonical basis of the multi-level scheme.
This SH-based algorithm leads to restriction and interpolation procedures performed in spectral space to transfer 
the solution between different spatio-temporal levels. The proposed restriction and interpolation methods do not 
introduce spurious modes that would, driven by nonlinear interactions, propagate across the spectrum. Our results 
show that this is one the key features needed to obtain an efficient MLSDC-SH scheme. The development of restriction 
and interpolation operators for other non-global spatial discretization schemes is left for future work.

We have shown that MLSDC-SH is efficient for the nonlinear wave-propagation-dominated problems arising from the 
discretized Shallow-Water Equations (SWE) on the rotating sphere. Our numerical studies are based on challenging 
test cases that are representative of the horizontal effects present in the full atmospheric dynamics. With a 
steady zonal jet test case,  we have first examined the impact of the coarsening strategy on the observed 
accuracy of MLSDC-SH upon refinement in time. Then, using unsteady numerical examples, we have shown 
that MLSDC-SH can achieve up to eighth-order convergence upon refinement in time, and that MLSDC-SH can take 
stable time steps that are as large as those of the single-level SDC schemes. We have also demonstrated 
that MLSDC-SH is more efficient than the single-level SDC schemes, and in particular requires fewer 
function evaluations. Our results show that MLSDC-SH can reduce the wall-clock time of the simulations by 
up to 37\% compared to single-level SDC schemes.


As a final note, we mention here that MLSDC is one of the key building blocks of the Parallel 
Full Approximation Scheme in Space and in Time (PFASST). The present work therefore lays the foundations 
of a parallel-in-time integration of the full shallow-water equations on the sphere with PFASST.

\section{Code availability}

The code used to generate the simulations presented here is publicly available \citep{schreiber2018sweet}.

\section{Acknowledgements}
The work of Fran\c cois Hamon and Michael Minion was supported by the U.S. Department of Energy, Office of Science,
Office of Advanced Scientific Computing Research, Applied Mathematics program under contract number DE-AC02005CH11231. Part of the simulations were performed using resources of the National Energy Research Scientific Computing Center (NERSC), a DOE Office of Science User Facility supported by the Office of Science of the U.S. Department of Energy under Contract No. DE-AC02-05CH11231.
Martin Schreiber gratefully acknowledges support from the Computational and Information Systems Laboratory (CISL) Visitor Program at the National Center for Atmospheric Research in Boulder, CO.
We thank Andreas Kreienbuehl for discussions on his work on MLSDC and PFASST using the Spectral Element Method.

\bibliography{biblio}

\end{document}